\documentclass[12pt]{amsart}
\usepackage{amssymb}
\usepackage{amsfonts}
\usepackage{latexsym}
\usepackage{amscd}
\usepackage[mathscr]{euscript}
\usepackage{xy} \xyoption{all}

\usepackage{graphicx}

\usepackage{amsmath}
\usepackage{amscd}
\usepackage{amssymb}
\usepackage{latexsym,nicefrac}

\newtheorem{definition}{Definition}[section]
\newtheorem{theorem}[definition]{Theorem}
\newtheorem{lemma}[definition]{Lemma}
\newtheorem{corollary}[definition]{Corollary}
\newtheorem{remark}[definition]{Remark}
\newtheorem{example}[definition]{Example}
\newtheorem{conjecture}[definition]{Conjecture}

\newtheorem{notation}[definition]{Notation}

\newtheorem{proposition}[definition]{Proposition}

\input xy
\xyoption{all}

\newcommand{\uloopr}[1]{\ar@'{@+{[0,0]+(-4,5)}@+{[0,0]+(0,10)}@+{[0,0] +(4,5)}}^{#1}}
\newcommand{\uloopd}[1]{\ar@'{@+{[0,0]+(5,4)}@+{[0,0]+(10,0)}@+{[0,0]+ (5,-4)}}^{#1}}
\newcommand{\dloopr}[1]{\ar@'{@+{[0,0]+(-4,-5)}@+{[0,0]+(0,-10)}@+{[0, 0]+(4,-5)}}_{#1}}
\newcommand{\dloopd}[1]{\ar@'{@+{[0,0]+(-5,4)}@+{[0,0]+(-10,0)}@+{[0,0 ]+(-5,-4)}}_{#1}}

\newcommand{\luloop}[1]{\ar@'{@+{[0,0]+(-8,2)}@+{[0,0]+(-10,10)}@+{[0, 0]+(2,2)}}^{#1}}

\newcommand{\V}{\mathcal V}
\newcommand{\Z}{\mathbb{Z}}
\newcommand{\N}{\mathbb{N}}

\newcommand{\M}{{ \rm M}}
\newcommand{\C}{\mathbb{C}}

\begin{document}

\title{Leavitt path algebras: \ the first decade}



\author{Gene Abrams
}





\begin{abstract}
The algebraic structures known as {\it Leavitt path algebras} were initially  developed  in 2004  by Ara, Moreno and Pardo, and almost simultaneously (using  a different approach) by the author and Aranda Pino.   
  During the intervening decade, these algebras have attracted significant interest and attention, not only from ring theorists, but from analysts working in C$^*$-algebras, group theorists,  and symbolic dynamicists as well.   The goal of this article is threefold:  to introduce the notion of Leavitt path algebras to the general mathematical community; to present some of the important results in the subject; and to describe some of the field's currently  unresolved questions.  
\end{abstract}

\maketitle

\medskip

Our goal in writing this article is threefold:  first,  to provide a history and overall viewpoint of the ideas which comprise the subject of Leavitt path algebras; second, to give the reader a general sense  of the results which have been achieved in the field; and finally, to give a broad picture of some of the research lines which are currently being pursued.     The history and overall viewpoint (Section \ref{intro}) are presented with a completely general mathematical audience in mind; the writing style here will be more chatty than formal. Our description of the results in the field has been split into two pieces:  we  describe  Leavitt path algebras of   row-finite graphs (Sections \ref{rowfinitesection}, \ref{section:InternalStructure}, and \ref{ModuleSection}), and  subsequently discuss various generalizations of these  (Section \ref{arbitrarygraphssection}).    Our intent and hope in ordering the presentation this way  is to allow the non-expert to appreciate the key ideas of the subject,  without getting ensnarled in the at-first-glance formidable constructs which drive the generalizations. We close with Section \ref{section:CurrentLines}, in which we describe some of the current lines of investigation in the subject.  In part, our hope here is  to attract mathematicians from a wide variety of fields to join in the research effort.  

The exhilarating  increase in the level of interest in Leavitt path algebras during the first decade since their introduction  has resulted in the publication of scores of articles on these and related structures.  Certainly it is {\it not} the goal of the current article   to review the entirety of the literature in the subject.  Rather, we have tried to strike a balance between presenting enough information to make clear the beauty and diversity of the subject on the one hand, while avoiding ``information overload" on the other.  Apologies are issued in advance to those authors whose work in the field has consequently not been included herein.     

In keeping with our goal of making this article accessible to a broad audience,  we will offer either a complete, formal {\it Proof}, or an intuitive, informal {\it Sketch of Proof},  only for specific results for which such proofs are particularly illuminating.     In other situations we will simply present statements without proof.   Appropriate  references are provided for all key results.  

\smallskip

{\it Acknowledgments}:  
This work was partially supported by Simons Foundation Collaboration Grant \#20894.  The author thanks P.N. \'{A}nh, Pere Ara, Zachary Mesyan, Paul Muhly, Enrique Pardo, and Mark Tomforde for helping to clarify some of the historical aspects of the subject.  The author is extremely grateful to S. Paul Smith for providing  a summary of germane ideas related to noncommutative algebraic geometry, and for a number of suggestions which clarified the presentation. Finally, the author warmly thanks the referee for providing significant, valuable feedback on the original version of this article. 


\section{ \ History and overview}
\label{intro}

\smallskip


{\bf 
\subsection{ \ History and overview:  Module type, and Leavitt's Theorem}\label{section:LeavittAlgebras}
}

The fundamental examples of rings that are encountered during one's  algebraic pubescence  (e.g., fields $K$, $\Z$, $K[x]$, $K[x,x^{-1}]$, $ {\rm M}_n(K)$) all have the following property.

\begin{definition}\label{def:IBN}   The unital ring $R$ has the {\it Invariant Basis Number} (IBN) property in case,   for each pair $i,j \in \N$, if the left $R$-modules $R^i$ and $R^j$ are isomorphic, then $i=j$.   
\end{definition}

A wide class of rings can easily be shown to have the IBN property, including rings possessing   any sort of reasonable chain condition on one-sided ideals, as well as commutative rings.   But there are naturally occurring examples of algebras which are not IBN. 

 Let $V$ be a countably infinite dimensional vector space over a field $K$, and let $R$ denote ${\rm End}_K(V)$, the algebra of linear transformations from $V$ to itself.  It is not hard to see that  $R^i \cong R^j$ as left (or right) $R$-modules for any pair $i,j\in \N$, as follows.    One starts by viewing $R$  as the $K$-algebra ${\rm RFM}_\N(K)$ consisting of those $\N \times \N$ matrices $M$ having the  property that each row of $M$ contains at most finitely many nonzero entries.   (In this context we view transformations as acting on the right, and define composition of transformations by setting $f\circ g$ to mean ``first $f$, then $g$".  Of course, depending on the reader's tastes,  the same analysis can be performed  by considering the analogous algebra ${\rm CFM}_\N(K)$ of column-finite matrices.)     Then a left-module isomorphism $ {}_RR \rightarrow {}_RR^2$ is easy to establish, by considering  the map that associates with any row-finite matrix $M$ the pair of row-finite matrices $(M_1,M_2)$, where $M_1$ is built from the odd-indexed columns of $M$, and $M_2$ from the even-indexed columns.     
 Once such an isomorphism ${}_RR \cong {}_RR^2$ is guaranteed, then by using the obvious generalization of the observation ${}_RR^3 \cong {}_RR^2 \oplus {}_RR \cong {}_RR \oplus {}_RR \cong {}_RR,$ we see that ${}_RR^i \cong {}_RR^j$ for all $i,j\in \N$.

The following is easy to prove.  
 
 \begin{lemma}\label{lem:RisoR^n}
 Let $R$ be a unital ring, and let $n\in \N$.  Then ${}_RR \cong {}_RR^n$ as left $R$-modules if and only if there exist elements $x_1, \hdots, x_n, y_1, \hdots, y_n$ of $R$ for which 
 \begin{equation} \tag{$\dagger$}
y_i x_j = \delta_{i,j} 1, \ \ \mbox{and} \ \ \sum_{i=1}^n x_iy_i = 1.
\end{equation}
 \end{lemma}

In effect, the ring $R = {\rm End}_K(V)$ 
lies on the complete opposite end of the spectrum from the IBN property, in that {\it every} pair of finitely generated free left $R$-modules are isomorphic;  such a ring is said to have  the Single Basis Number (SBN) property.   The question posed (and answered completely) by William G. Leavitt in the  early 1960's regards the existence of a  middle ground between IBN and SBN:  do there exist rings for which ${}_RR^i \cong {}_RR^j$ for some, but not all, pairs $i,j \in \N$?    If we assume an isomorphism exists between ${}_RR^i$ and ${}_RR^j$ for some pair $i\neq j$, then clearly by appending $k$ copies of $R$  to this isomorphism we get that ${}_RR^{i+k} \cong {}_RR^{j+k}$.  With this idea in mind, and using only basic properties of  the semigroup $\N \times \N$, it is easy to prove the following.    

\begin{lemma}\label{lem:Leavittnecessary} 
Let $R$ be a unital ring.      
Assume $R$ is not IBN, i.e., that there exist  $i\neq j \in \N$ with ${}_RR^i \cong {}_RR^j$.  Let $m$ be the least integer for which ${}_RR^m \cong {}_RR^j$ for some $j \neq m$.  For this $m$, let $n$ denote be the least integer for which ${}_RR^m \cong {}_RR^n$ and $n>m$.  Let $k$ denote $n-m \in \N$.    Then for any pair $i,j \in \N$,
$${}_RR^i \cong {}_RR^j \ \ \Longleftrightarrow   \ \  i,j \geq m \ \mbox{ and } \  i \equiv j \ (mod k).$$
\end{lemma}

We call the pair $(m,n)$ the {\it module type of} $R$. (We caution that some authors, including Leavitt, instead use the phrase {\it module type} to denote the pair $(m,k)$.)   In particular, the ring $R = {\rm End}_K(V)$ has module type $(1,2)$. 

With Lemma \ref{lem:Leavittnecessary} in mind, Leavitt proved the following ``anything-that-{\it can}-happen-actually-{\it does}-happen" result.

\begin{theorem}{\rm (}{\bf Leavitt's Theorem}{\rm )} \cite[Theorem 8]{Leav62}\label{LeavittsTheorem}  
Let $m,n \in \N$ with $n>m$, and let $K$ be any field.     Then there exists a  $K$-algebra $L_K(m,n)$ having module type $(m,n)$.   Additionally:

(1)  $L_K(m,n)$ is universal, in the sense that if $S$ is any $K$-algebra having module type $(m,n)$, then there exists a nonzero $K$-algebra homomorphism $\varphi: L_K(m,n) \rightarrow S$.

(2)  $L_K(m,n)$ is simple (i.e., has no nontrivial two-sided ideals) if and only if $m=1$.  (This was shown in \cite[Theorems 2 and 3]{Leav64}.)  In this case, for each $0\neq x \in L_K(1,n)$, there exists $a,b\in L_K(1,n)$ for which $axb = 1$.  

(3)  $L_K(m,n)$ is  explicitly described in terms of generators and relations.

\end{theorem}
We will refer to  $L_K(m,n)$ as the {\it Leavitt algebra of type} $(m,n)$.   

Since the ring $R =  {\rm RFM}_\N(K) \cong {\rm End}_K(V) $ has module type $(1,2)$, by Lemma \ref{lem:RisoR^n}  there necessarily exists a set of  four elements in $R$ which satisfy the appropriate relations; for clarity, we note that one such set is given by 
$$Y_1 = (\delta_{i,2i-1}), \ \ Y_2 = (\delta_{i,2i}), \ \ X_1 = (\delta_{2i-1,i}), \ \mbox{and}  \ \ X_2 = (\delta_{2i,i}).$$
Moreover, the subalgebra of ${\rm RFM}_\N(K)$ generated by these four matrices is isomorphic to $L_K(1,2)$.  

In the specific case when $m=1$, the explicit description of the algebra $L_K(1,n)$ given in \cite{Leav62}  yields that $L_K(1,n)$ is the free associative $K$-algebra  $K\langle  x_1, \hdots, x_n, y_1, \hdots, y_n\rangle $ modulo the relations given in $(\dagger)$.  
So, with Lemma \ref{lem:RisoR^n} in mind, we may view $L_K(1,n)$ as essentially the ``smallest" algebra of type $(1,n)$.   

\smallskip

In the following subsection we will rediscover the algebras $L_K(1,n)$ from a different point of view.


{\bf 
\subsection{ \ History and overview:  The $\mathcal{V}$-monoid,  and Bergman's Theorem}
\label{section:introtoV-monoid}
}

We begin this subsection with a well-known idea, perhaps cast more formally than is typical.  Let $P$ be a finitely generated projective left $R$-module.   (The notions of ``finitely generated" and ``projective" are categorical, and make sense even in case $R$ is nonunital.) 

\smallskip

\begin{definition}\label{V(R)definition}
For a ring $R$, let $\mathcal{V}(R)$ denote the set of  isomorphism classes of finitely generated projective left $R$-modules, and define the obvious  binary operation $\oplus$ on $\mathcal{V}(R)$ by setting  $[P] \oplus [Q] = [P \oplus Q]$.    Then $(\mathcal{V}(R),\oplus)$  is easily seen to be a commutative monoid (with neutral element $[0]$). 
\end{definition}
 
For any idempotent $e\in R$,   $[Re] \in \V(R)$; indeed, elements of $\V(R)$ of this form will play a central role in the subject.  If $R$ is a division ring, then $\V (R) \cong \mathbb{Z}^+$; the same is true for $R = \mathbb{Z}$, as well as for various additional classes of rings.    The wide range of  monoids which can arise as $\V(R)$ will be demonstrated in Theorem \ref{BergmansTheorem}.  For an arbitrary ring $R$,   it's fair to say that an explicit description of $\V(R)$ is typically hard to come by.      The well-studied {\it Grothendieck group} $K_0(R)$ {\it of} $R$ is precisely the universal abelian group corresponding to the commutative monoid $\mathcal{V}(R)$.  
 
 When $R$ is unital then $[R] \in \V (R)$.   Key information about $R$ may be provided  by the pair $(\V (R) , [R])$.    If $R$ and $S$ are isomorphic as rings then there exists an isomorphism of monoids $\varphi: \V(R) \rightarrow \V(S)$ for which $\varphi([R]) = [S]$.  (In this situation we write $(\V (R) , [R]) \cong (\V (S) , [S])$.) 
 More generally,   if $R$ and $S$ are Morita equivalent  (denoted $R \sim_M S$), then  there exists an isomorphism of monoids $\psi: \V (R) \rightarrow \V (S)$.     However, such an induced isomorphism $\psi$  need not have the property $\psi([R]) = [S]$.  For instance, if $S = \M_2(K)$ for a field $K$, then $S \sim_M K$, and $\V (S) \cong \V(K) \cong \Z^+$.  But $(\V(K),[K]) = (\Z^+, 1)$, while $(\V(S), [S]) =  (\Z^+, 2)$, and clearly no automorphism of the monoid $\Z^+$ can take $1$ to $2$.  

In addition to the commutativity of $(\V (R), \oplus)$,  the monoid $\V (R)$  has the following two easy-to-see properties.    First,  $\V (R)$ is {\it conical}:   if $x, y \in \V(R)$ have $x \oplus y = 0$, then $x = y = 0$.  Second (for $R$ unital),  $\V (R)$ contains a {\it distinguished element}  $d$:   for each $x\in \V(R)$, there exists $y \in \V(R)$ and $n\in \N$ having $x \oplus y = nd$ (specifically, $d = [R]$).     In 1974, George Bergman  established the following remarkable result.

\begin{theorem} \rm{(}{\bf Bergman's Theorem})\label{BergmansTheorem} (\cite[Theorem 6.2]{Berg74}) 
Let $M$ be a finitely generated commutative conical monoid with distinguished element $d \neq 0$, and let $K$ be any field.   Then there exists a $K$-algebra $B = B(M,d)$ for which $(\V(B), [B]) \cong (M,d).$   Additionally:

1) $B = B(M,d)$ is universal, in the sense that if $C$ is any unital $K$-algebra for which there exists a monoid homomorphism $\psi: \V(B) \rightarrow  \V(C)$ having $\psi([B]) = [C]$, then there exists a (not necessarily unique) $K$-algebra homomorphism $\Psi: B \rightarrow C$ for which  the induced map $\overline{\Psi}:  C \otimes_B\underbar{\hspace{.15in}}: \V(B) \rightarrow \V(C)$ is precisely $\psi$.

2)   $B = B(M,d)$ is left and right hereditary (i.e., every left ideal and every right ideal of $B$ is projective). 

3)   The construction of $B(M,d)$ depends on the specific representation of $M$ as $\mathcal{F}/\langle \mathcal{R}\rangle$, where $\mathcal{F}$ is a finitely generated free abelian monoid, and $\mathcal{R}$ is a given (finite) set of relations in $\mathcal{F}$.   With $\mathcal{F}$ and $\mathcal{R}$ viewed as starting data, the algebra $B = B(M,d) = B(\mathcal{F} / \langle \mathcal{R} \rangle ,d)$ is constructed explicitly via a finite sequence of steps, where each step consists  of adjoining elements satisfying explicitly specified relations (provided by $\mathcal{R}$) to an explicitly described algebra. 

\smallskip

We will refer to $B = B(M,d) =  B(\mathcal{F} / \langle \mathcal{R}\rangle, d)$ as the {\it Bergman algebra} of $(\mathcal{F} / \langle \mathcal{R} \rangle,d)$.    

\end{theorem}

\begin{example}{\bf Three important examples.}\label{FundamentalExamples}
{\rm 

\medskip

(1)  
Perhaps not surprisingly, when $(M,d) =   (\Z^+ ,1) =  (\Z^+ / \langle \emptyset \rangle,1)$, then $B(\Z^+ / \langle \emptyset \rangle,1) = K$.   (In this situation we view $\Z^+$ as the free abelian monoid on one generator.)  

\smallskip

(2)  
Somewhat more subtly, we consider the same pair $(\Z^+,1)$, but this time represent the monoid $\Z^+$ as  $\Z^+ / \langle 1 = 1 \rangle$.    Then $B(\Z^+ / \langle 1 = 1 \rangle, 1) = K[x,x^{-1}]$, the {\it Laurent polynomial algebra} with coefficients in $K$.   

\smallskip

(3)    Let $n\in \N, n \geq 2$.  Let $V_n$ denote the free abelian monoid having a single generator $x$, subject to the relation $nx=x$.  So $V_n  = \{ 0,x,2x, \hdots , (n-1)x\}$, $ |M_n| =n$,  and $(M_n,x)$ clearly satisfies the hypotheses of Bergman's Theorem.  (In \cite{Berg74}, the semigroup $V_n$ is denoted $V_{1,n}$.)  
In this situation, Bergman's explicit construction yields that  $B(V_n,x) = K\langle x_1,  \hdots,  x_n, y_1, \hdots, y_n\rangle$, with relations given by exactly the same defining relations as given in $(\dagger)$ above, namely, 
$$y_i x_j = \delta_{i,j} 1, \ \ \mbox{and} \ \ \sum_{i=1}^n x_iy_i = 1.$$
Consequently, as observed in \cite[Theorem 6.1]{Berg74}, 
the Bergman algebra $B(V_n,x)$ is precisely the Leavitt algebra $L_K(1,n)$.  
}
\end{example}

\bigskip

\noindent

%
%
%
%

{\bf 
\subsection{ \ History and overview:  Graph C$^*$-algebras}
}

Because of the central role they play in both the genesis and the ongoing development Leavitt path algebras, no history of the subject would be complete without a discussion of {\it graph  C}$^*${\it-algebras}.  We present here only the most basic description of these algebras, just enough so that even the reader who is completely unfamiliar with them can get a sense of their connection to Leavitt path algebras.

Throughout this subsection all algebras are assumed to be unital algebras over the complex numbers $\C$ (but most of these ideas can be cast significantly more generally).     
 The algebra $A$ is a {\it $*$-algebra} in case there is a map $\ast : A \rightarrow A$ which has: \  $(x+y)^* = x^* + y^*$; \ $(xy)^* = y^*x^*$; \ $1^* = 1$; \ $(x^*)^* = x$; \ and $(\alpha \cdot x)^* = \overline{\alpha} x^*$ for all $x,y\in A$ and $\alpha \in \C$, where $\overline{\alpha}$ denotes the complex conjugate of $\alpha$.     Standard examples of $*$-algebras include matrix rings ${\rm M}_n(\C)$ (where $\ast$ is `conjugate transpose'), and the ring $C(\mathbb{T})$ of continuous functions from the unit circle $ \mathbb{T} = \{ z \in \C \ | \ |z| = 1 \}$ to $\C$ (where $\ast$ is defined by setting $f^*(z) = \overline{f(z)}$ for $z\in \mathbb{T}$).

A C$^*$-{\it norm} on a  $*$-algebra $A$ is a function $\|
\cdot \| \colon A\to \mathbb R^+$ for which:  

\noindent
 $\| ab\| \le \|a\|\cdot  \| b\| $;  \ 
 $\|a+b\| \le \| a\| + \| b\| $; \ 
 $\| aa^*\| =\| a\| ^2=\| a^*\| ^2 $;  \ $\| a\| =0$ $\Leftrightarrow$ $a=0$; \ and 
 $\| \lambda a \| =|\lambda | \|a\| $ for all $a, b \in A$ and $\lambda \in \mathbb
C$.     For $A = {\rm M}_n(\C)$, a $C^*$-norm on $A$ is given by operator norm, where we view elements of ${\rm M}_n(\C)$ as operators $\C^n \rightarrow \C^n$, with the Euclidean norm on $\C^n$.  (This operator norm assigns to $M \in {\rm M}_n(\C)$  the square root of the largest eigenvalue of the matrix $M^*M$.)  A C$^*$-norm on $C(\mathbb{T})$  is also given by an operator norm.  

A C$^*$-norm on a $\ast$-algebra $A$ induces a topology on $A$ in the usual way, by defining the $\epsilon$-ball around an element $a\in A$ to be $\{ b\in A \ | \ \| b - a \| < \epsilon \}$.  

\medskip

 \begin{definition}\label{def:Cstaralgebra}
 {\rm  A {\it $C^*$-algebra} is a
 $*$-algebra $A$ endowed with a $C^*$-norm $\|\cdot \|$, for which $A$ is complete with respect to the topology induced by $\|\cdot \|$. 
}
\end{definition}


A second description of a C$^*$-algebra, from an operator-theoretic point of view, is given here.  Let $\mathcal{H}$ be a Hilbert space, and let $B(\mathcal{H})$ denote the continuous linear operators on $\mathcal{H}$.  A C$^*$-algebra is an adjoint-closed subalgebra of $B(\mathcal{H})$ which is closed with respect to the norm topology on $B(\mathcal{H})$.     In general, and especially relevant in the current context, one often builds a  C$^*$-algebra   by starting with  a given set of elements in $B(\mathcal{H})$, and then forming the smallest C$^*$-subalgebra of $B(\mathcal{H})$ which contains that set. 

A {\it partial isometry} is an element $x$ in a C$^*$-algebra $A$ for which $y = x^*x$ is a self-adjoint idempotent; that is, in case $y^* = y$ and $y^2 = y$.   Such elements are characterized as those elements $z$ of $A$ for which $zz^*z = z$ in $A$.   For instance, in ${\rm M}_n(\C)$, any element which is the sum of distinct matrix units $e_{i,i}$ ($1\leq i \leq n$) is a partial isometry (indeed, a projection);  there are other partial isometries in ${\rm M}_n(\C)$ as well.  Since the only idempotents in $C(\mathbb{T})$ are the constant functions $0$ and $1$, it is not hard to show that the set of partial isometries in $C(\mathbb{T})$ consists of $\{0\} \cup \{ f \in C(\mathbb{T}) \ | \ f(\mathbb{T}) \subseteq \mathbb{T}\}$.   

The study of C$^*$-algebras has its roots in the early development of quantum mechanics; these  were used to model algebras of physical observables.    Various questions about the structure of C$^*$-algebras arose over the years.  One of the most important of these questions, the explicit description of a separable simple infinite C$^*$-algebra,  was resolved in 1977 by Cuntz (\cite[Theorem 1.12]{Cuntz}).   A C$^*$-algebra is {\it simple} in case it contains no nontrivial closed two-sided ideals.  (It can be shown that this is equivalent to the algebra containing no nontrivial two-sided ideals, closed or not.)   A C$^*$-algebra is {\it infinite} in case it contains an element $x$ for which $xx^* = 1$ and $x^*x \neq 1$.  

\medskip

{\bf Cuntz' Theorem}  Let $n \in \N$.  Consider a Hilbert space $\mathcal{H}$, and a set $\{S_i\}_{i=1}^n$ of isometries (i.e., $S_i^* S_i = 1$) on $\mathcal{H}$.  Assume that $\sum_{i=1}^n S_i S_i^* = 1$.  Let $\mathcal{O}_n$ denote $C^*(S_1, \dots, S_n)$,  the C$^*$-algebra generated by $\{S_i\}_{i=1}^n$.   Then the infinite separable C$^*$-algebra $\mathcal{O}_n$ is simple.

\medskip

Indeed, Cuntz proves much more in \cite[Theorem 1.12]{Cuntz} than we have stated here.   Additionally,  it is shown in \cite[Theorem 1.13]{Cuntz} that if $X$ is any nonzero element in $\mathcal{O}_n$, then there exist $A,B \in \mathcal{O}_n$ for which $AXB = 1$.  

  Cuntz notes that the condition  $\sum_{i=1}^n S_i S_i^* = 1$ implies that the $S_i S_i^*$ are pairwise orthogonal.    So the C$^*$-algebra $\mathcal{O}_n$ is the C$^*$-completion of a $\C$-subalgebra $T_n$ of $B(\mathcal{H})$, where $T_n$ as a $\C$-algebra is generated by isometries  $\{S_i\}_{i=1}^n$, for which $\sum_{i=1}^n S_i S_i^* = 1$.   Since a C$^*$-algebra is adjoint-closed, we see that $\mathcal{O}_n$ may also be viewed as the C$^*$-completion of a $\C$-subalgebra $L_n$ of $B(\mathcal{H})$  generated by isometries  $\{S_i\}_{i=1}^n$ {\it together with} $\{S_i^*\}_{i=1}^n$, for which $\sum_{i=1}^n S_i S_i^* = 1$.     

\smallskip

In retrospect, such a $\C$-algebra $L_n$ is seen to be isomorphic to $L_\C(1,n)$.    

\smallskip

Subsequent to the appearance of \cite{Cuntz}, a number of researchers in  operator algebras investigated natural generalizations of the Cuntz C$^*$-algebras $\mathcal{O}_n$; see especially \cite{CuntzKrieger}.     In the early 1980's, various constructions  of C$^*$-algebras corresponding to  directed graphs were studied by Watatani and others (e.g., \cite{Watatani82}).  Even though, via this approach, the Cuntz algebra $\mathcal{O}_n$ could be realized as the C$^*$-algebra  corresponding to the 
graph $R_n$ (see Example \ref{example:Leavitt} below), this  methodology  did not gain much traction at the time.  Instead, the study of these C$^*$-algebras from a different point of view (arising from matrices with non-negative integer entries, or arising from groupoids) became more  the vogue.     But then, in the fundamental article \cite{KPRR97} (in which groupoids are still in the picture, and the corresponding graphs could not have sinks), and the subsequent followup articles \cite{KPR98} and \cite{BPRS00},  
the power of constructing a C$^*$-algebra based on the data provided by a directed graph became clear.   

\begin{definition}\label{def:finitegraph}
{\rm  A (directed) {\it graph}  is a quadruple $E = (E^0, E^1, s,r)$, where $E^0$ and $E^1$ are sets (the {\it vertices} and {\it edges} of $E$, respectively), and $s$ and $r$ are functions from $E^1$ to $E^0$ (the {\it source} and {\it range} functions of $E$, respectively).   A {\it sink} is an element $v\in E^0$ for which $s^{-1}(v) = \emptyset$.  $E$ is {\it finite} in case both $E^0$ and $E^1$ are finite sets.
}
\end{definition}

\begin{definition}\label{graphCstardefinition}
{\rm 
Let $E$ be a finite graph.  Let $C^*(E)$ denote the universal C$^*$-algebra generated by a collection of mutually orthogonal projections $\{p_v \ | \ v\in E^0\}$ together with partial isometries $\{ s_e \ | \ e\in E^1\}$ which satisfy the {\it Cuntz-Krieger relations}:

(CK1)  $s_e^* s_e = p_{r(e)}$ for all $e\in E^1$, and

(CK2)  $p_v = \sum_{\{e\in E^1  | s(e) = v\}} s_e s_e^*$ for each non-sink $v\in E^0$.
}
\end{definition}


For example, in \cite{KPR98} the authors were able to identify those finite graphs $E$ for which $C^*(E)$ 
 is simple, and those for which $C^*(E)$ is  purely infinite simple.  (The germane graph-theoretic terms will be described in Notations \ref{graphnotation} and \ref{notation:heredandsat} below.  A unital C$^*$-algebra $A$ is {\it purely infinite simple} in case $A \not\cong \C$, and for each $0\neq x \in A$ there exist $a,b\in A$ with $axb = 1$.)    
 
 \begin{theorem}\label{Cstarsimplicity} {\bf (Simplicity and Purely Infinite Simplicity Theorems for graph C$^*$-algebras)}
 Let $E$ be a finite graph.   Then $C^*(E)$ is simple   if and only if the only hereditary saturated subsets of $E$ are trivial, and every cycle in $E$ has an exit.  Moreover, $C^*(E)$ is purely infinite simple if and only if C$^*(E)$ is simple, and $E$ contains at least one cycle.  
 \end{theorem}

   Subsequently, in  \cite{BPRS00},  results were clarified, sharpened, and extended; and the groupoid techniques were eliminated from the arguments.

During the same timeframe,  Kirchberg (unpublished) and Phillips (\cite{Phillips})  independently proved a beautiful, deep result which classifies up to isomorphism a class of C$^*$-algebras satisfying various properties.  Although the now-so-called {\it Kirchberg Phillips Theorem} covers a wide class of C$^*$-algebras, it manifests in the particular case of purely infinite simple graph C$^*$-algebras as follows.

\begin{theorem}\label{KPforgraphalgebras} {\bf (The Kirchberg Phillips Theorem for graph C$^*$-algebras)}   Let $E$ and $F$ be  finite graphs.  Suppose $C^*(E)$ and $C^*(F)$ are purely infinite simple.  Suppose there is an isomorphism $K_0(C^*(E)) \cong K_0(C^*(F))$ for which $[1_{C^*(E)}] \mapsto [1_{C^*(F)}].$   Then $C^*(E) \cong C^*(F)$.
\end{theorem}

The work described in \cite{BPRS00} became the basis of a newly-energized research program in the C$^*$-algebra community, a program which continues to flourish to this day.  For additional information about graph C$^*$-algebras, see \cite{RaeburnBook};  for a more complete description of the history of graph C$^*$-algebras, see \cite[Appendix B]{MarksHistoryBook}.


%
%
%
%

\bigskip

 {\bf
\subsection{ \ History and overview: The confluence of many ideas leads to the definition of Leavitt path algebras}
}

With the overview of Leavitt algebras, Bergman algebras, and graph C$^*$-algebras now in place, we are in position to describe the genesis of Leavitt path algebras.

\medskip

There are two plot lines to the history.    

\bigskip

 {\bf
\subsubsection{ \ Historical Plot Line \#1: Graph algebras as Bergman algebras}
}  \label{confluencesubsection}

The first historical plot line begins with an investigation into the algebraic notion of purely infinite simple rings, begun by Ara, Goodearl, and Pardo   (each of whom has significant expertise in both ring theory and C$^*$-algebras) in \cite{AGPpis}.  In it, the authors   ``... extend the notion of a purely infinite simple C$^*$-algebra to the context of unital rings, and study its basic properties, especially those related to $K$-theory".

   The authors note in the introduction of \cite{AGPpis} that ``The Cuntz algebra $\mathcal{O}_n$ is the C$^*$-completion of the Leavitt algebra $V_{1,n}$ over the field of complex numbers."    Although this connection between the Cuntz and Leavitt algebras is now viewed as almost obvious, it was not until the early 2000's that such a connection was first noted in the literature.  (A somewhat earlier mention of this connection appears in  \cite{AbAn02}; the observation in \cite{AbAn02} was included at the request of an anonymous referee.)

With the notion of purely infinite simple rings so introduced, the same three authors (together with Gonz\'{a}lez-Barroso)  set out to find large classes of explicit examples of such rings.   With the purely infinite simple graph C$^*$-algebras as motivation, the four authors in \cite{AGGP04} introduced the  ``algebraic Cuntz-Krieger (CK) algebras."   (Retrospectively, these are seen to be the Leavitt path algebras corresponding to finite graphs having neither sources nor sinks, and which do not consist of a disjoint union of cycles.)     These algebraic Cuntz-Krieger algebras arose as  specific examples of  {\it fractional skew monoid rings}, and the germane ones were shown to be purely infinite simple by using techniques which applied to the more general class.   

With  the $K$-theory of the corresponding graph C$^*$-algebras in mind, it was then natural to ask analogous $K$-theoretic questions about the algebraic CK algebras.   In addition, earlier  work by Ara, Goodearl, O'Meara and Pardo \cite{AGOP98} regarding semigroup-theoretic properties of $\V$ (e.g., separativity and refinement) for various classes of rings provided the motivation to ask similar questions about $\V(A)$ for these algebras.


Once various specific examples had been completely worked out, it became clear to Ara and Pardo that much of the information about the $\V$-monoid of the algebraic CK  algebras could be seen directly in terms of relations between vertices and edges in an associated graph $E$.   Indeed, these relations between vertices and edges could be codified as information which could then be used to generate a monoid in a natural way, defined here.  


\bigskip

\begin{definition}\label{def:graphmonoid}
Let $E$ be a finite graph, with $E^0 = \{v_1, v_2, \hdots, v_n\}$.   The {\it graph monoid} $M_E$ of $E$ is the free abelian monoid on a generating set $\{a_{v_1}, a_{v_2}, \hdots, a_{v_n}\}$, modulo the relations $$a_{v_i} = \sum_{\{e \in E^1 | s(e) = v_i\}} a_{r(e)} \ \ \ \  \mbox{for each non-sink} \ v_i \ .$$
\end{definition}

In a private communication to the author, Enrique Pardo wrote that, with all this information and background as context,  
\medskip

\begin{tabular}{@{\hspace{3ex}}p{29em}}
{\it  ... at some moment [early in 2004] one of us suggested that probably Bergman's coproduct construction would be a good manner of solving the computation and prove that both monoids coincide.}
\end{tabular}
\medskip

Once some additional necessary machinery was included (the notion of a {\it complete} subgraph), then Ara and Pardo, together with Pardo's colleague Mariangeles Moreno-Fr\'{i}as, had all the ingredients in hand to make the  following definition, and prove the subsequent theorem, in \cite{AMP07}.   (We state the definition and theorem here only for finite graphs; these results were established for more general graphs in \cite{AMP07}, and the general version will be discussed below.)  




\bigskip

\begin{definition}\label{def:graphalgebra}(\cite[p. 161]{AMP07})
Let $E$ be a finite graph, and let $K$ be a field.  We define the {\it graph $K$-algebra}  $L_K(E)$ associated with $E$ as the $K$-algebra generated by a set $\{p_v \  | \ v\in E^0\}$ together with a set $\{x_e, y_e \ | \ e\in E^1\}$, which satisfy the following relations:

(1)  $p_vp_{v'} = \delta_{v,v'}p_v$ for all $v,v' \in E^0.$

(2)  $p_{s(e)}x_e = x_ep_{r(e)} = x_e$ for all $e\in E^1.$

(3)  $p_{r(e)}y_e = y_ep_{s(e)} = y_e$ for all $e\in E^1.$

(4)  $y_ex_{e'} = \delta_{e,e'}p_{r(e)}$ for all $e,e' \in E^1$.

(5)  $p_v = \sum_{\{e\in E^1|s(e)=v\}}x_ey_e$ for every $v\in E^0$ that emits edges.
\end{definition}

We note that both the terminology used in this definition (``graph algebra"), as well as  the notation, is quite similar to the terminology and notation which was already being employed in the context of graph C$^*$-algebras.

\begin{theorem}\label{AMPfinite}  {\bf (The Ara / Moreno / Pardo Realization Theorem)} {\rm (\cite[Theorem 3.5]{AMP07})}
Let $E$ be a finite graph and $K$ any field.  Then there is a natural monoid isomorphism $\V(L_K(E)) \cong M_E$.  
\end{theorem}

By examining the proof of \cite[Theorem 3.5]{AMP07}, and using Bergman's Theorem, we can in fact restate this fundamental result as follows.  

\bigskip
\noindent
{\bf Theorem} {\bf \ref{AMPfinite}$^\prime$}  {\bf (The Ara / Moreno / Pardo Realization Theorem, restated)}   
{\it Let $E$ be a finite graph and $K$ any field.   Let $M_E$ be the monoid given by the specific set of generators and relations presented in Definition \ref{def:graphmonoid}.  Let $d$ denote the element $\sum_{v\in E^0}a_v$ of $M_E$.     Then $L_K(E) \cong B(M_E,d)$.    Consequently, $\V(L_K(E)) \cong M_E$.  Moreover,  $L_K(E)$ is hereditary.   }

\bigskip

In the same groundbreaking article \cite{AMP07}, Ara, Moreno, and Pardo were also  able to establish a connection between the 
 $\V$-monoids of $L_K(E)$ and $C^*(E)$.

 \begin{theorem}\label{IsoVmonoids}  {\bf (The Ara / Moreno / Pardo Monoid Isomorphism Theorem)} {\rm (\cite[Theorem 7.1]{AMP07})}
Let $E$ be a finite graph.  Then there is a natural monoid isomorphism $\V(L_\C(E)) \cong \V(C^*(E))$.  
\end{theorem}

We conclude our discussion of Historical Plot Line \#1  in the development of Leavitt path algebras by again quoting Enrique Pardo:   

\medskip

\begin{tabular}{@{\hspace{3ex}}p{29em}}
{\it For us the motivation was to give an algebraic framework to all these families of (purely infinite simple) $C^*$-algebras associated to combinatorial objects, say Cuntz-Krieger algebras and graph $C^*$-algebras. For this reason we always looked at properties that were known in $C^*$ case and were related to combinatorial information: we wanted to know which part of these results relies in algebraic
information, and which ones in analytic information. So, we looked at
K-Theory, stable rank, exchange property (in $C^*$-algebras this is real
rank zero property), prime and primitive ideals, the classification
problem and Kirchberg-Phillips Theorem...}
\end{tabular}

\medskip

We will visit each of these topics later in the article.

\bigskip

 {\bf
\subsubsection{ \ Historical Plot Line \#2: Leavitt path algebras as quotients of quiver  algebras}
}


\bigskip

The second historical plot line begins with the author's interest in Leavitt's algebras, specifically the algebras $L_K(1,n)$.   For instance, these algebras were used in \cite{Abra97}  to produce non-IBN rings having unexpected isomorphisms between their matrix rings; were used again in \cite{Abra01} to solve a question (posed in \cite{NTV}) about strongly graded rings;  and were subsequently investigated yet again in \cite{AbAn02}, in joint work with P.N. \'{A}nh of the R\'{e}nyi Institute of Mathematics (Hungarian Academy of Sciences, Budapest).

During a Spring 2001 visit to the University of Iowa, \'{A}nh met the analyst Paul Muhly.\footnote{The \'{A}nh / Muhly meeting was quite fortuitous.  \'{A}nh was a visiting research guest of  Kent Fuller at the University of Iowa during Spring Semester 2001. Fuller regularly went to lunch at various Iowa City restaurants with a group of his departmental colleagues, an excursion in which Paul Muhly was a frequent participant;   Fuller  of course invited   \'{A}nh to join in.}  
Subsequently,  \'{A}nh invited Muhly to give a talk at the R\'{e}nyi Institute  (during a 2003 trip that Muhly and his wife were making to Budapest anyway, to visit their son); it was during this visit that the two mathematicians began to consider the potential for connections between various topics.   
   Muhly was one of the organizers of the May/June 2004  NSF -  CBMS conference\footnote{National Science Foundation - Conference Board in Mathematical Sciences.  The NSF-CBMS Regional Research Conferences in the Mathematical Sciences are a series of five-day conferences, each of which features a distinguished lecturer delivering ten lectures on a topic of important current research in one sharply focused area of the mathematical sciences.} ``Graph Algebras: Operator Algebras We Can See", delivered by Iain Raeburn,    held at the University of  Iowa.   Muhly consequently extended invitations to attend that conference to the author, to  \'{A}nh, and to a handful of other ring theorists.\footnote{V. Camillo, L. M\'{a}rki, and E. Ortega also attended.}  During conference coffee break discussions, the algebraists began to realize that when one considered the ``pre-completion" version  of the   graph C$^*$-algebras, the remaining algebraic structure  looked quite familiar, specifically, as some sort of modification of the well-known notion of a {\it quiver algebra} or {\it path algebra}.
\bigskip

\begin{definition}\label{def:pathalgebra}
Let $F = (F^0,F^1,s,r)$ be a graph and  $K$  any field.   The {\it path} $K${\it -algebra of} $F$ (also known as the {\it quiver} $K${\it -algebra} of $F$), denoted $KF$,  is the $K$-vector space having basis $F^0 \sqcup F^1$, with multiplication given by the $K$-linear extension of 
$$ p \cdot q  =\begin{cases}
 \ pq & \text{if }r(p)=s(q)\\
 \ 0 & {\rm otherwise.}\end{cases}
$$  (If $v\in F^0$ we denote $s(v) = v = r(v)$.)    
\end{definition}

Gonzalo Aranda Pino visited the author's home institution for the period July 2004 through December 2004.\footnote{As if the  \'{A}nh / Muhly meeting (and consequent  attendance of the ring theorists at the 2004 Iowa CBMS conference) was not fortuitous enough, it turned out that, many months prior to that conference, Mercedes Siles Molina had contacted the author regarding the possibility of having the author host one of her Ph.D. students for a six month visit at the University of Colorado, to commence July 2004.  That having been arranged, Gonzalo Aranda Pino arrived in Colorado Springs at precisely the time that this new idea was blossoming.}
   Early in Aranda Pino's visit, the author shared with him some of the ideas which had been discussed in Iowa City during the previous month.       A few weeks of collaborative effort subsequently led to the following.

\begin{definition}\label{def:extendedgraph}
Given a directed graph $E = (E^0, E^1, s,r)$  we define the  {\it extended graph of} $E$ as the  graph
$\widehat{E}=(E^0,E^1\sqcup (E^1)^*,s',r'),$ where $(E^1)^*=\{e_i^*:e_i\in  E^1\}$, and the functions $r'$ and $s'$ are
defined by setting  $r'|_{E^1}=r,\ s'|_{E^1}=s,\ r'(e_i^*)=s(e_i)$, and    $s'(e_i^*)=r(e_i).$
\end{definition}
\begin{definition}\label{def:Lpafinite}   
Let $E$ be a finite graph and $K$ any field. The {\it Leavitt path $K$-algebra $L_K(E)$} is defined as the path $K$-algebra $K\widehat{E}$, 
modulo the relations:

\smallskip
 \ \ \  (CK1) \  $e_i^*e_j=\delta_{ij}r(e_j)$ for every $e_j\in E^1$ and  $e_i^*\in (E^1)^*$. 
 
 \ \ \ \  (CK2) \ 
$v_i=\sum_{\{e_j\in E^1|s(e_j)=v_i\}}e_je_j^*$ for every  $v_i\in E^0$ which is not a sink.


\end{definition}

Some of the notation which was developed in the C$^*$-algebra context is also used in  the Leavitt path algebra world, e.g., the use of the ``CK" labels to denote the two key relations. (Cf. Definition \ref{graphCstardefinition}).  

With both Leavitt's Theorem (part 2 of Theorem \ref{LeavittsTheorem}) and The Simplicity Theorem for graph C$^*$-algebras (Theorem \ref{Cstarsimplicity}) in mind, the author and Aranda Pino focused their initial investigation on an internal, multiplicative question about  the algebras $L_K(E)$: for which graphs $E$ and fields $K$ is $L_K(E)$ simple?  Using techniques completely unlike those utilized to achieve Theorem  \ref{Cstarsimplicity}, the following result was established.  (See Notations \ref{graphnotation} and \ref{notation:heredandsat} below for definitions of appropriate terms.)

\begin{theorem}\label{thm:AAPsimple}{\bf (The Abrams / Aranda Pino Simplicity Theorem)} \cite[Theorem 3.11]{AAP05})     Let $E$ be a finite graph and $K$ any field.   Then $L_K(E)$ is simple if and only if the only hereditary saturated subsets of $E$ are trivial, and every cycle in $E$ has an exit.     
\end{theorem}

{\bf 
\subsubsection{ \ The confluence of the two Historical Plot Lines }
} 
By making the obvious correspondences $v \leftrightarrow p_v$, $e \leftrightarrow x_e$, and $e^* \leftrightarrow y_e$, we see immediately: 

\bigskip

\begin{center}
For a finite graph $E$ and field $K$, \\ the graph $K$-algebra   of Definition \ref{def:graphalgebra}  is the same algebra \\ as the Leavitt path $K$-algebra  of Definition \ref{def:Lpafinite}. 
\end{center}

\bigskip

It is of historical interest to note that  the work on  \cite{AAP05} was started in July 2004.  Subsequently, \cite{AAP05} was submitted for publication in September 2004, accepted for publication in June 2005,  appeared online in  September 2005, and appeared in print in November 2005. 
On the other hand, the work on  \cite{AMP07} was started in early 2004.  Subsequently,  \cite{AMP07} was submitted for publication in late 2004 (and posted on ArXiV at that time), and accepted for publication in early 2005, but  did not appear in print until April 2007.   
 So even though \cite{AAP05} appeared in print eighteen months prior to the appearance in print of \cite{AMP07}, in fact most the mathematical work done to produce the latter preceded that of the former.

 \medskip
 
 Both \cite{AAP05} and \cite{AMP07} should be viewed as the foundational articles on the subject.
\bigskip

%


%

\section{Leavitt path algebras of row-finite graphs: general properties and examples}\label{rowfinitesection}

Section \ref{intro}  of this article was meant to give the reader an overall view of the motivating ideas which led naturally to the construction of Leavitt path algebras.  Over the next three sections  we describe some of the key ideas and results for Leavitt path algebras arising from {\it row-finite} graphs. Subsequently,   in Section \ref{arbitrarygraphssection} we relax this hypothesis on the graphs.  (For those results which do not extend verbatim to the unrestricted case, we will indicate in the statement that the graph must be row-finite (or finite); otherwise, we will make no such stipulation in the statement.)     


\begin{notation}\label{rowfinitenotation} A vertex $v$ in a graph $E = (E^0, E^1, s,r)$ is called {\it regular} in case  $0 < |s^{-1}(v)| < \infty$; otherwise, $v$ is called {\it singular}.   Specifically, if $s^{-1}(v) = \emptyset$ then $v$ is called a {\it sink}, while $v$ is called an {\it infinite emitter} in case  $|s^{-1}(v)|$ is infinite.    $E$ is called {\it row-finite} in case $E$ contains no infinite emitters.   
\end{notation}

Here is the formal definition of a Leavitt path algebra arising from a row-finite graph.  

\begin{definition}\label{def:Lpageneral}
Let $E = (E^0, E^1, s, r)$ be a row-finite graph and $K$ any field.  Let $\widehat{E}$ denote the extended graph of $E$.   The {\it Leavitt path $K$-algebra $L_K(E)$} is defined as the path $K$-algebra $K\widehat{E}$, 
modulo the relations:

\smallskip
(CK1) \ $e_i^*e_j=\delta_{ij}r(e_j)$ for every $e_j\in E^1$ and  $e_i^*\in (E^1)^*$. 

 (CK2) \ $v_i=\sum_{\{e_j\in E^1|s(e_j)=v_i\}}e_je_j^*$ for every  non-sink $v_i\in E^0$. 
 
 \smallskip
 \noindent
Equivalently, we may define $L_K(E)$ as the free associative $K$-algebra on generators $E^0 \sqcup E^1 \sqcup (E^1)^*$, modulo the relations 

(1)  $v{v'} = \delta_{v,v'}v$ for all $v,v' \in E^0.$

(2)  $s(e)e = er(e) = e$ for all $e\in E^1.$

(3)  $r(e)e^* = e^*s(e) = e^*$ for all $e\in E^1.$

(4)  $e^*e' = \delta_{e,e'}r(e)$ for all $e,e' \in E^1$.

(5)  $v = \sum_{\{e\in E^1|s(e)=v\}}ee^*$ for every non-sink $v\in E^0$. 
\end{definition}

It is established in \cite{TomfordeUniqueness} that the expected map from $L_\C(E)$ to $C^*(E)$ is in fact injective.  With this 
and the construction of the graph $C^*$-algebra $C^*(E)$, we get

\begin{proposition}\label{densesubalgebra}
For any graph $E$,   $L_\C(E)$ is isomorphic to a dense $\ast$-subalgebra of $C^*(E)$. 
\end{proposition}

The interplay between graphs and algebras will play a major role in the theory.   It is important to note at the outset that in general, if $F$ is a subgraph of $E$, then $L_K(F)$ need {\it not} correspond to a subalgebra of $L_K(E)$, because the (CK2) relation imposed at a vertex $v$ in $L_K(F)$   need not be the same  as the relation imposed at $v$ in $L_K(E)$.      For a row-finite graph $E$, a subgraph $F$ is said to be {\it complete} in case, whenever $v\in F^0$, then either $s^{-1}_F(v) = \emptyset$, or $s^{-1}_F(v) = s^{-1}_E(v)$.  (In other words, if $v\in F^0$, then either $v$ emits no edges in $F$, or emits the same edges in $F$ as it does in $E$.)   Perhaps not surprisingly, when $F$ is a complete subgraph of $E$, then there is an injection of algebras $L_K(F) \hookrightarrow L_K(E)$.   Moreover, 

\begin{proposition}\label{rowfinitedirectlimit}\cite[Lemma 3.2]{AMP07}   
The assignment $E \mapsto L_K(E)$ can be extended to a functor $L_K$ from the category of row-finite graphs and complete graph inclusions to the category of $K$-algebras and (not necessarily unital) algebra homomorphisms.   The functor $L_K$ commutes with direct limits.  It follows that every $L_K(E)$ for a row-finite graph $E$ is the direct limit of graph algebras corresponding to finite graphs. 
\end{proposition}  
 
Because of Proposition \ref{rowfinitedirectlimit}, it is often the case that a result which holds for the Leavitt path algebras of  finite graphs can be extended to the row-finite case.   

\begin{definition}\label{def:Efamily} Let $E$ be any graph and $A$ any $K$-algebra.   A {\it Leavitt $E$-family in}    $A$ is a subset $\mathcal{S} =     \{ a_v \ | \ v\in E^0\} \cup \{ b_e \ | \ e\in E^1\} \cup \{ c_e \ | \ e\in E^1\}$ of $A$ for which

(1)  $a_va_{v'} = \delta_{v,v'}a_v$ for all $v,v' \in E^0.$

(2)  $a_{s(e)}b_e = b_ea_{r(e)} = b_e$ for all $e\in E^1.$

(3)  $a_{r(e)}c_e = c_ea_{s(e)} = c_e$ for all $e\in E^1.$

(4)  $c_eb_{e'} = \delta_{e,e'}a_{r(e)}$ for all $e,e' \in E^1$.

(5)  $a_v = \sum_{\{e\in E^1|s(e)=v\}}b_ec_e$ for every non-sink $v\in E^0$. 

\end{definition}

By the description of $L_K(E)$ as a quotient of a free associative $K$-algebra modulo the germane relations given in Definition \ref{def:Lpageneral}, we immediately get the following result, which often proves to be quite useful in the subject. 

\begin{proposition}\label{prop:homstoEfamilies}{\bf (Universal Homomorphism Property of Leavitt path algebras)}  
Let $E$ be a graph, and suppose $\mathcal{S}$ is a Leavitt $E$-family in the  $K$-algebra $A$.  Then there exists a unique $K$-algebra homomorphism $\varphi: L_K(E) \rightarrow A$ for which $\varphi(v) = a_v, \varphi(e) = b_e,$ and $\varphi(e^*) = c_e$ for all $v\in E^0$ and $e\in E^1$.   
\end{proposition}

\begin{notation}\label{graphnotation}


A sequence of edges $\alpha = e_1, e_2, ... , e_n$ in a graph $E$ for which $r(e_i) = s(e_{i+1})$ for all $1\leq i \leq n-1$ is called a {\it path of length} $n$.  We typically denote such $\alpha$ more simply by $ e_1e_2 \cdots  e_n$.    Each vertex $v$ of $E$ is viewed as a path of length $0$.   The set of paths of length $n$ in $E$ is denoted by $E^n$; the set of all paths in $E$ is denoted ${\rm Path}(E)$.  So we have ${\rm Path}(E) = \sqcup_{n\in \Z^+}E^n$.  

For $\alpha = e_1e_2\cdots e_n \in {\rm Path}(E)$, $s(\alpha)$ denotes $s(e_1)$, $r(\alpha)$ denotes $r(e_n)$, and ${\rm Vert}(\alpha)$ denotes the set $\{s(e_1), s(e_2), \dots, s(e_n), r(e_n) \}$.   The path $ e_1 e_2 \cdots e_n$ is {\it closed} if $s(e_1) = r(e_n)$.  A closed path $ c = e_1 e_2 \cdots e_n$ is {\it simple} in case $s(e_i) \neq s(e_1)$ for all $2\leq i \leq n$. Such a simple closed path $c$ is said to be {\it based at} $v=s(e_1)$.      
 A simple closed path $ c = e_1 e_2 \cdots e_n$ is a  {\it cycle} in case there are no repeats in the  list of vertices $s(e_1), s(e_2), ..., s(e_n)$.   $E$ is called {\it acyclic} in case there are no cycles in $E$.  
 
An {\it exit}  for a path $ e_1 e_2  \cdots e_n$  is an edge $f \in E^1$ for which $s(f) = s(e_i)$ and $f\neq e_i$ for some $1\leq i \leq n$.

The graph $E$ satisfies {\it Condition (L)} in case every cycle in $E$ has an exit.

The graph $E$ satisfies {\it Condition (K)} in case no vertex in $E$ is the base of exactly one simple closed path in $E$.

\end{notation}

If $\alpha = e_1 e_2 \cdots e_n$ is a path in $E$, then we may view $\alpha$ as an element of the path algebra $KE$, and as an element of the Leavitt path algebra $L_K(E)$ as well.  (In this sense, concatenation in the graph $E$ is interpreted as multiplication in $KE$ or $L_K(E)$.)   We denote by $\alpha^*$ the element $e_n^* \cdots e_2^* e_1^*$ of $L_K(E)$.   We often refer to a path $\alpha = e_1 e_2 \cdots e_n$ of $E$ (viewed as an element of $L_K(E)$)  as a {\it real path}, while an element of $L_K(E)$ of the form $\alpha^* = e_n^* \cdots e_2^* e_1^*$  is called a {\it ghost path}.     Here are some easily verified basic properties of Leavitt path algebras. 

\begin{proposition}\label{properties}    Let $E$ be any graph and $K$ any field  

(1)   Every nonzero element $r$ of $L_K(E)$ may be written ({\it not necessarily uniquely}) as 
$$r = \sum_{i=1}^n k_i \alpha_i \beta_i^*,$$
where $k_i \in K^\times$, and $\alpha_i, \beta_i \in {\rm Path}(E)$ with $r(\alpha_i) = r(\beta_i)$ for $1\leq i \leq n$.    

\smallskip

(2)     For each $\alpha \in {\rm Path}(E)$, $\alpha^* \alpha = r(\alpha)$.

\smallskip

(3)      The natural $K$-algebra map $KE \rightarrow L_K(E)$ is a one-to-one homomorphism.  

\smallskip

(4)      $L_K(E)$ is unital (with multiplicative identity $\sum_{v\in E^0}v$) if and only if $E^0$ is finite.  In general, $L_K(E)$ has a set of enough idempotents, consisting of finite sums of distinct vertices.

\smallskip

(5)     The map $\ast: L_K(E) \rightarrow L_K(E)$ induces an isomorphism  $L_K(E) \cong L_K(E)^{op}$.  In particular, for Leavitt path algebras, the categories of left $L_K(E)$-modules and right $L_K(E)$-modules are isomorphic.   
\end{proposition}

\subsection{Examples of familiar / ``known"  algebras which arise as Leavitt path algebras}

 We saw in Section \ref{intro}  how specific algebras arise from Bergman's Theorem, starting with a specified monoid.    We re-examine those here, and present additional examples  as well.

\begin{example}\label{example:Mn(K)}{\bf \ Full matrix $K$-algebras.} \ \ 
Let $A_n$ denote the graph
$$ \xymatrix{ \bullet^{v_1} \ar[r]^{e_1} & \bullet^{v_2} \ar[r]^{ \hspace{-.4in} e_2}   & \ \ \  \cdots  \ \ \ \bullet^{v_{n-1}} \ar[r]^{\hspace{.2in} e_{n-1}} &\bullet^{v_{n}}} $$
\noindent
Then $L_K(A_n) \cong {\rm M}_n(K)$.   This is not hard to see.  We present two different approaches, in order to play up the germane ideas.   

 The first approach:   consider the standard matrix units $\{ E_{i,j} \ | \ 1\leq i,j \leq n\}$ in  ${\rm M}_n(K)$.  Since each vertex (other than $v_n$) emits a single edge, the (CK2) relation at these vertices becomes $e_ie_i^* = v_i$.   Using this, it is straightforward to  
 verify that the set 
 \smallskip
 
\hspace{-.15in}  $\mathcal{S} = \{ E_{i,i} \ | \ 1\leq i \leq n\} \  \cup \ \{  E_{i,i+1} \ | \ 1\leq i \leq n-1\} \ \cup \ \{ E_{i+1,i} \ | \ 1\leq i \leq n-1\}$

\smallskip
 \noindent
  is an $A_n$-family in ${\rm M}_n(K)$.  So  the Universal Homomorphism Property ensures the existence of a $K$-algebra homomorphism $\varphi$ for which $\varphi(v_i) = E_{i,i}$, $\varphi(e_i) = E_{i,i+1}$, and $\varphi(e_i^*) = E_{i+1,i}$.  That $\varphi$ is an isomorphism is easily checked (for instance, by constructing the expected  function  $\psi: {\rm M}_n(K) \rightarrow L_K(A_n)$, and verifying that $\psi = \varphi^{-1}$). 
  
 The second approach:   we analyze the monoid $M_{A_n}$, define $d = \sum_{i=1}^n a_{v_i}$, and see easily that $(M_{A_n},d) \cong (\Z^+,n)$.  With the relations describing $M_{A_n}$, it is clear that $B(M_{A_n},d) \cong {\rm M}_n(K).$  Now Theorem \ref{AMPfinite}$^\prime$ applies.  

\end{example}

Full matrix rings over $K$ arise as the Leavitt path algebras of graphs other than the $A_n$ graphs.  In Theorem \ref{finiteacyclic}  below we will justify the isomorphisms asserted in the next two examples.  These two examples play up the fact that non-isomorphic graphs may have isomorphic Leavitt path algebras.  (This observation lies at the heart of much of the current research activity in Leavitt path algebras.)  

\begin{example}{\bf \ Full matrix $K$-algebras, revisited.} \ \ 
For $n\in \N$ let $B_n$ denote the graph 
$$\xymatrix{ \bullet^{w_1} \ar[dr]^{e_1} & \bullet^{w_2} \ar[d]^{e_2} & \bullet^{w_3} \ar[dl]_{e_3} \\
 & \bullet^v &  \ar@{.>}[l] \\
\bullet^{w_{n-1}} \ar[ur]^{e_{n-1}}  & \ar@{.>}[u] & \ar@{.>}[ul]}$$
Then $L_K(B_n) \cong {\rm M}_{n}(K).$  
\end{example}

\begin{example}\label{example:Dn}{\bf \ Full matrix $K$-algebras, again revisited.} \ \ 
For $n\in \N$ let $D_n$ denote the graph
\smallskip

$$
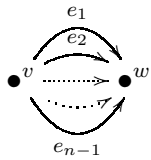
\xymatrix{ \bullet^v \ar@{.>}[r] \ar@/^{10pt}/ [r]^{e_2}    \ar@/^{20pt}/ [r]^{e_1}  \ar@{.>} @/^{-10pt}/ [r]   \ar@/^{-20pt}/ [r]_{e_{n-1}}    & \bullet^w}$$

\noindent
Then $L_K(D_n) \cong {\rm M}_{n}(K)$. 

\end{example}

Proceeding in a manner similar to that utilized in Example \ref{example:Mn(K)}, one can easily establish the following two claims.  (See Example \ref{FundamentalExamples}.)

\begin{example}\label{example:Laurent}{\bf \ The Laurent polynomial $K$-algebra.} \ \ 
Let $R_1$ denote the graph
$$ 
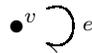
\xymatrix{ \bullet^v \ar@(ur,dr)^e} $$
Then $L_K(R_1) \cong K[x,x^{-1}]$, the Laurent polynomial algebra.  The isomorphism is clear:   $v \mapsto 1$, $e \mapsto x$, and $e^* \mapsto x^{-1}$.
\end{example}

Here is the Fundamental Example of Leavitt path algebras. 

\begin{example}\label{example:Leavitt}{\bf \  Leavitt $K$-algebras.} \ \ 
For $n\geq 2$, let $R_n$ denote the graph
$$ 
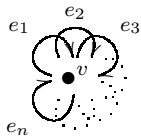
\xymatrix{ \bullet^v       \ar@(l,u)^{e_1}      \ar@(ul,ur)^{e_2}   \ar@(u,r)^{e_3}   \ar@{.} @(ur,dr)  \ar@{.} @(dl,dr)  \ar@{.} @(r,d)^{} \ar@(d,l)^{e_n} }$$
Then $L_K(R_n) \cong L_K(1,n)$, the Leavitt algebra of order $n$.  The isomorphism is clear:   using the description of the generators and relations for $L_K(1,n)$ given in $(\dagger)$ above, $v \mapsto 1$, $e_i \mapsto x_i$, and $e_i^* \mapsto y_i$.
\end{example}

\begin{example}\label{example:Toeplitz}{\bf \ The Toeplitz $K$-algebra.} \ \ 
For any field $K$, the {\it Jacobson algebra}, described in \cite{Jacobson1950}, is the $K$-algebra
$$A = K\langle x,y \ | \ xy = 1\rangle.$$
This algebra was the first example appearing in the literature of an algebra which is not {\it directly finite}, that is, in which there are elements $x,y$ for which $xy=1$ but $yx\neq 1$.  
Let $\mathcal{T}$ denote the ``Toeplitz graph"
$$\xymatrix{ \bullet^v \ar@(ul,dl)_e \ar[r]^f & \bullet^w}$$

\noindent
Then $L_K(\mathcal{T}) \cong A$.   The isomorphism is not hard to write down explicitly.     First, the set 

\hspace{.75in} $\mathcal{S} = \{yx, 1-yx\} \cup \{y^2x, y-y^2x\} \cup \{yx^2, x - yx^2\}$

\smallskip
\noindent
 is easily shown to be a $\mathcal{T}$-family in $A$, so by the Universal Homomorphism Property of Leavitt path algebras there exists a $K$-algebra homomorphism $\varphi: L_K(\mathcal{T}) \rightarrow A$ for which $\varphi(v) = yx$, $\varphi(w) = 1-yx$, $\varphi(e) = y^2x$, $\varphi(f) = y - y^2x$, $\varphi(e^*) = yx^2$, and $\varphi(f^*) = x - yx^2$.   On the other hand, we define $X = e^* + f^*, Y = e+f$ in $A$.  Using (CK1) and (CK2) we get easily that $XY = 1$.  This gives a $K$-algebra homomorphism $\psi: A \rightarrow L_K(\mathcal{T})$, the algebra extension of  $x \mapsto X$ and $y \mapsto Y$.   It is easy to check that $\varphi$ and $\psi$ are inverses.   
\end{example}

\begin{example}\label{example:matricesoverL(E)}{\bf \ Full matrix $K$-algebras over $L_K(E)$.} \ \ 
Let $E$ be any graph, $K$ any field,  and $n\in \N$.  The graph $E(n)$ is defined as follows.  For each $v\in E^0$, one  adds to $E$ the following vertices and edges

\hspace{1in} $\xymatrix{\bullet^{v^{n}} \ar[r]^{f_v^{n-1}} & \bullet^{v^{n-1}} \ar@{.>}[r] & \cdots \ar[r]^{f_v^2} & \bullet^{v^{2}} \ar[r]^{f_v^1} & }, $

\noindent
where $r(f_v^1) = v.$  Then $L_K(E(n)) \cong {\rm M}_n(L_K(E))$.  (See \cite[Proposition 9.3]{AT}.)
\end{example}

\begin{example}\label{example:Minfty}{\bf \ Infinite matrix $K$-algebras.} \ \ 
Let $I$ be any set.  We denote by ${\rm M}_{I}(K)$ the set of those $I \times I$ matrices $M$, having  entries in $K$, for which $M_{i,j} = 0$ for at most finitely many pairs $(i,j)$.   Then ${\rm M}_{I}(K)$ is a $K$-algebra, which is unital if and only if $I$ is finite (and in this case  ${\rm M}_{I}(K)$ consists of {\it all} $I \times I$ matrices  having  entries in $K$).  When $I$ is infinite, then ${\rm M}_{I}(K)$ has a set of enough idempotents, consisting of finite sums of distinct matrix units of the form $E_{i,i}$.   

If $A_\N$ denotes the graph 
$$ \xymatrix{ \bullet^{v_1} \ar[r]^{e_1} & \bullet^{v_2} \ar[r]^{e_2}   &  \bullet^{v_3}  \ar@{.>}[r] 
&} $$
then $L_K(A_\N) \cong M_\N(K)$.   

More generally, for any infinite set $I$, let $B_I$ denote the graph having  vertices $\{v\} \cup \{w_i \ | \ i\in I\}$, and edges $\{ e_i \ | \ i\in I\}$, with $s(e_i) = w_i$ and $r(e_i) = v$ for all $i\in I$.   Then $L_K(B_I) \cong {\rm M}_I(K)$.

\end{example}




\section{Internal / multiplicative properties of Leavitt path algebras}
\label{section:InternalStructure}
Not surprisingly, a number of the key results in the subject focus on passing structural information from the directed graph $E$ to the Leavitt path algebra $L_K(E)$, and vice versa; i.e.,  results of the form
 \begin{equation} \tag{$\dagger \dagger$}
E \mbox{ has graph property } \mathcal{P} \ \Longleftrightarrow \ L_K(E)  \mbox{ has ring property } \mathcal{Q}.\end{equation}

\smallskip
\noindent
The Simplicity Theorem (Theorem \ref{thm:AAPsimple})  is the quintessential result of this type.   We will describe a number of additional such results in this section and the next.  In the author's opinion,  these results are quite interesting, some even {\it remarkable}, in their own right.      Just as compellingly, some of these results have been utilized to produce heretofore unrecognized classes of algebras having interesting ring-theoretic properties.

Looking ahead:  in contrast, in the next section, we will engage in a   discussion of the equally important ``external / module-theoretic"  properties of Leavitt path algebras.  As described in Section \ref{intro}, the ``internal / multiplicative" and ``external / module-theoretic" properties form the historical foundations of the subject.   We will see in the final section that these also drive much of the current investigative energy.

\subsection{Finite dimensional Leavitt path algebras}\label{findimsubsection}

We start by analyzing the Leavitt path algebras of finite acyclic graphs.  From a ring-theoretic point of view, these  turn out to be the most basic (least interesting?) of all the Leavitt path algebras.  

\begin{theorem}\label{finiteacyclic}{\bf Structure Theorem of Leavitt path algebras for finite acyclic graphs.}
Let $E$ be a finite acyclic graph and $K$ any field.   Let $w_1, \dots, w_t$ denote the sinks of $E$.  (At least one sink must exist in any finite acyclic graph.)   For each $w_i$, let $N_i$ denote the number of elements of ${\rm Path}(E)$ having range vertex $w_i$.   (This includes $w_i$ itself, as a path of length $0$.)   Then $$L_K(E) \cong \bigoplus_{i=1}^t {\rm M}_{N_i}(K).$$
\end{theorem}
{\it Sketch of Proof.}    For each sink $w_i$ consider the ideal $I(w_i)$ of $L_K(E)$. If $\alpha, \beta \in {\rm Path}(E)$ have $r(\alpha) = r(\beta) = w_i$, then $\alpha \beta^* = \alpha w_i \beta^* \in I(w_i)$.   Using the (CK1) relation with the fact that $w_i$ is a sink, one shows easily that the set of $N_i^2$ elements $\{ \alpha \beta^* \ | \ \alpha, \beta \in {\rm Path}(E), \ r(\alpha)=r(\beta) = w_i\}$ is a set of matrix units, which yields that $I(w_i) \cong {\rm M}_{N_i}(K).$   That the sum $\sum_{i=1}^t I(w_i)$ is direct follows by again using  the hypothesis that the $w_i$ are sinks.  Now let $\gamma \delta^*$ be any monic monomial in $L_K(E)$.  If $v = r(\gamma)$ is a sink, then $\gamma \delta^* \in \sum_{i=1}^t I(w_i)$.  Otherwise, the (CK2) relation may be invoked at $v$, and we may write 
$$\gamma \delta^* = \gamma v \delta^* = \gamma (\sum_{e\in s^{-1}(v)}ee^*) \delta^* = \sum_{e\in s^{-1}(v)}(\gamma e) (\delta e)^*.$$
If $r(e)$ is a sink, then the expression $(\gamma e) (\delta e)^*$ is in $\sum_{i=1}^t I(w_i)$; if not, then in the same manner one can use the (CK2) relation at $r(e)$ to rewrite  $(\gamma e) (\delta e)^*$.  Since $E$ is finite and acyclic, the process must terminate with expressions of the desired form.  \hfill $\Box$

\medskip

So for finite acyclic graphs, the resulting Leavitt path algebras are, among other things: unital semisimple; left artinian; and finite dimensional.  Indeed, any of these three ring/algebra-theoretic properties characterizes the Leavitt path algebras of finite acyclic graphs, thus yielding three examples of results of type $(\dagger \dagger)$.    Perhaps more importantly, Theorem \ref{finiteacyclic} yields a result of the following type:  among a certain class of graphs (specifically, finite acyclic), we can  determine, using easy-to-compute graph-theoretic properties, which of those graphs yield isomorphic Leavitt path algebras (specifically, those for which the number of sinks, and the corresponding $N_i$, are equal).    So Theorem  \ref{finiteacyclic} may be viewed as a very basic type of Classification Theorem.     

\begin{notation}\label{notation:heredandsat}
Let $E$ be any graph, and let $v,w\in E^0$.  We write $v\geq w$ in case  there exists  $p \in {\rm Path}(E)$ for which $s(p)=v$ and $r(p)=w$.    

Let $X$ be a subset of $E^0$.     $X$ is called {\it hereditary} in case, whenever $v\in X$ and $w\in E^0$ and $v\geq w$, then $w\in X$.   $X$ is called {\it saturated} in case, whenever $v\in E^0$ is regular and $r(s^{-1}(v))\subseteq X$, then $v\in X$.   (Less formally:   $X$ is saturated in case whenever $v$ is a non-sink in $E$ which emits finitely many edges, and the range vertices of all of those edges are in $X$, then $v$ is in $X$ as well.)  

Clearly both $E^0$ and $\emptyset$ are hereditary  saturated subsets of $E^0$, and clearly the intersection of any collection of hereditary saturated subsets of $E^0$ is again hereditary saturated.   If $S$ is any subset of $E^0$, then  $\overline{S}$ denotes the smallest hereditary saturated subset of $E^0$ which contains $S$; $\overline{S}$ is called the {\it hereditary saturated closure of} $S$.  (Such exists by the previous observation.)  
\end{notation}

The interplay between vertices $E^0$ of $E$ on the one hand (viewed as idempotent elements of $L_K(E)$), and ideals of $L_K(E)$ on the other, plays a central role in the ideal structure of $L_K(E)$.   This connection clearly brings to light  the roles of the two (CK) relations in this context.

\begin{proposition}  Let $E$ be any graph and $K$ any field.  
Let $I$ be an ideal of $L_K(E)$.  Then $X = E^0 \cap I$ is a hereditary saturated subset of $E^0$.  
\end{proposition}

{\it Proof.}
Let $v \in X$.  If $w\in E^0$ for which there exists $\alpha \in {\rm Path}(E)$ with $s(\alpha) = v$ and $r(\alpha)=w$, then by Proposition \ref{properties}(2)  $$w = \alpha^*\alpha = \alpha^* \cdot v \cdot \alpha$$ in $L_K(E)$, so that $w\in I$, and thus in $X$.  So $X$ is hereditary.  On the other hand, suppose $v \in E^0$ has the property that $|s^{-1}(v)|$ is finite, and that $r(e) \in X$ for each $e\in s^{-1}(v)$.   But by (CK2) $$v = \sum_{e \in s^{-1}(v)} ee^* = \sum_{e \in s^{-1}(v)} e \cdot r(e) \cdot e^*$$ in $L_K(E)$, so that $v\in I$, and thus in $X$.   So $X$ is saturated.  \hfill $\Box$

\subsection{The $\Z$-grading, and graded ideals}  

A $K$-algebra $R$ is {\it $\Z$-graded} in case $R = \oplus_{i\in \Z}R_i$ as $K$-vector spaces, in such a way that $R_i \cdot R_j \subseteq R_{i+j}$ for all $i,j \in \Z$.  The subspaces $R_i$ are called the {\it homogeneous components} of $R$.   The Leavitt path algebras admit a $\Z$-grading, as follows.  Any path $K$-algebra $K\widehat{E}$ of an extended graph $\widehat{E}$ is $\Z$-graded, by setting ${\rm deg}(v) = 0$ for $v\in E^0$,  and ${\rm deg}(e) = 1$, ${\rm deg}(e^*) = -1$ for $e\in E^1$, and extending additively and multiplicatively.  Since the two sets of relations (CK1) and (CK2) consist of homogeneous elements of degree $0$ with respect to this grading on $K\widehat{E}$, the grading passes to the quotient algebra $L_K(E)$.     In particular, for $m\in \Z$, the homogeneous component $L_K(E)_m$  of degree $m$ consists of $K$-linear combinations of elements of the form $\alpha \beta^*$, where $r(\alpha) = r(\beta)$, and $\ell(\alpha) - \ell(\beta) = m$.   

A two-sided ideal $I$ in a $\Z$-graded ring $R$ is called a {\it graded ideal} in case, whenever $s \in I$ and $s = \sum_{j\in \Z} s_j$ is the decomposition of $s$ into homogeneous components, then $s_j \in I$ for each $j\in \Z$.   It is easy to show that if a two-sided ideal $I$ in a $\Z$-graded ring is generated by homogeneous elements of degree $0$, then $I$ is a graded ideal.   In particular, for any set of vertices $X \subset E^0$, the ideal $I(X)$ of $L_K(E)$ is graded.    In contrast, not all ideals of a Leavitt path algebra are necessarily graded; for instance, the ideal $I(1+x) \subset K[x,x^{-1}] \cong L_K(R_1)$ is not graded, as neither $1$ nor $x$ is in $I(1+x)$.

So on the one hand any ideal $I$ of $L_K(E)$ gives rise to the hereditary saturated subset $I \cap E^0$ of $E^0$, while on the other, any subset $X$ of $E^0$ gives rise to the graded ideal $I(X)$ of $L_K(E)$.   The perhaps-expected connection is the following.

\begin{proposition} \label{prop:latticeiso} \cite[Theorem 5.3]{AMP07} 
Let $E$ be a row-finite graph.   Then there is a lattice isomorphism between the lattice of graded ideals $\{I\}$  of $L_K(E)$ and the lattice of hereditary saturated subsets $\{X\}$  of $E^0$:
$$I \mapsto I\cap E^0    \ ,  \ \ \ \ \ \ X \mapsto I(X)  \ .$$
\noindent
In particular, every graded ideal of $L_K(E)$ is generated by vertices.
\end{proposition}

{\it Sketch of Proof.}   It is not hard to show that $I(X) = I(\overline{X})$.  On the other hand, if $v\in I(X)$, 
then by using an explicit, iterative description of the hereditary saturated closure of a set, one can show that $v\in \overline{X}$.      \hfill $\Box$ 

\medskip

The connection between these two lattices does not hold verbatim in case $E$ contains infinite emitters, as we will see in Section \ref{arbitrarygraphssection}.   

It was shown by Bergman \cite{BergUnpub} that if $R$ is a $\Z$-graded (unital) ring, then the Jacobson radical $J(R)$ is necessarily a graded ideal.   (See also \cite[Theorem 2.5.40]{RowenBook}.)  Using that $J(R)$ contains no nonzero idempotents in any ring $R$, Proposition \ref{prop:latticeiso} yields the following nice ``internal" result about Leavitt path algebras.

\begin{corollary}\label{cor:zeroradical}
Let $E$ be any graph and $K$ any field.  Then $L_K(E)$ has zero Jacobson radical.   
\end{corollary}


\subsection{Ideals in Leavitt path algebras}\label{generalidealssubsection}

In general, loosely speaking, the two key players in the graph $E$ which drive the ideal structure of  $L_K(E)$ are  the vertices, and  the cycles without exits.  While the hereditary saturated subsets  will dictate the graded structure of $L_K(E)$,  the cycles without exits (when $E$ contains such)  provide additional structural nuances.   The following result provides some motivation as to why this should be the case.   For an element $p(x) = \sum_{i=m}^N k_i x^i  \in K[x,x^{-1}]$ (with $k_N \in K^\times$), and a cycle $c$ in the graph $E$, we denote by $p(c)$ the element $\sum_{i=m}^N k_i c^i $ of $L_K(E)$, where $c^i = (c^*)^{-i}$ whenever $i < 0$, and $c^0 = s(c)$.  

\begin{theorem}\label{reductiontheorem}{\bf The Reduction Theorem.}  \cite[Proposition 3.1]{Malaga4}
Let $E$ be any graph and $K$ any field.  Let $0 \neq x \in L_K(E)$.  Then there exist $\alpha, \beta \in {\rm Path}(E)$ for which either:

(1) $ \alpha^* x \beta = kv$ for some $v\in E^0$ and $k\in K^\times$,  \ \ or

(2)  $ \alpha^* x \beta = p(c)$, where $0 \neq p(x) \in K[x,x^{-1}]$ and $c$ is a cycle without exits.  

\smallskip
\noindent
In other words,  we can transform (via multiplication by real paths and/or ghost paths) any element of $L_K(E)$ to either a nonzero multiple of a vertex, or to a nonzero polynomial in a cycle without exits.   
\end{theorem}

{\it Sketch of Proof.}    The proof uses an idea similar to the one Leavitt used in his proof of the Simplicity Theorem for $L_K(1,n)$ (\cite[Theorem 2]{Leav64}).   Essentially, starting with $x$, one shows that there is a path $\gamma$ in $E$ for which $0\neq x\gamma \in KE$.  This is done by finding $v\in E^0$ for which $xv \neq 0$, then writing $xv \in L_K(E)$ in a form which minimizes the length of the ghost terms from among all possible representations of $xv$, and then applying an induction argument.  With this in hand, one then modifies $x \gamma$ via left multiplication by terms of the form $\delta^*$ to ``reduce" $x\gamma$ to one of the two indicated forms.     \hfill $\Box$

 \medskip


\begin{definition}  For a hereditary saturated subset $H$ of $E$, let $C_H$ denote the set of cycles $c$ in $E$ for which ${\rm Vert}(c) \cap H = \emptyset$, and for which $r(e)\in H$ for every exit $e$ of $c$.      

For any subset $C$ of $C_H$, consider a set $P = P(C) = \{ p_c(x) \ | \ c\in C\}$ of noninvertible, nonzero elements of $K[x]$.    Let $P_C$ denote the subset  $\{ p_c(c) \ | \ c\in C\}$ of $L_K(E)$. 
\end{definition}


For instance, for the Toeplitz graph $\mathcal{T}$ described in Example \ref{example:Toeplitz},  let $H$ be the (only nontrivial) hereditary saturated subset $\{w\}$.  Then the cycle $e$ is in $C_H$.  For any polynomial $p(x) \in K[x]$, we may form the ideal $I(w, p(e))$ of $L_K(\mathcal{T})$ generated by the two elements $w$ and $p(e)$. 


In a similar way, for general graphs, using the data provided by a hereditary saturated subset $H$ of $E^0$, a set $C= C(H)$ of cycles which miss $H$ but all of whose exits land in $H$, and nontrivial polynomials $P = P(C) = P(C(H))$ in $K[x]$ (one for each element of $C$), we can build an ideal in $L_K(E)$, namely, the ideal generated by $H$ together with elements of $L_K(E)$ of the form $p_c(c)$.  Rephrased, starting with such $(H,C(H),P(C(H))$, we can build the ideal $I(H,\{p_c(c)\}_{c\in C})$.      Indeed, this process gives {\it all} the ideals of $L_K(E)$.

\begin{theorem}\label{structureofidealsrowfinite}{\bf Structure Theorem of Ideals}\cite[Theorem 2.8.10]{TheBook} 
Let $E$ be a row-finite graph.   Then every ideal of $L_K(E)$ is of the form $I(H,\{p_c(c)\}_{c\in C})$ as described above.   
\end{theorem}

Indeed, with  not-hard-to-anticipate order relations defined on triples of the form $(H, C, P)$, there is a stronger form of Theorem \ref{structureofidealsrowfinite}, one which gives a lattice isomorphism between the set of appropriate triples and the lattice of two-sided ideals of $L_K(E)$.

There are some immediate consequences of Theorem \ref{structureofidealsrowfinite}.  The most noteworthy of these is the   Simplicity Theorem (Theorem \ref{thm:AAPsimple}):  that $L_K(E)$ is simple if and only if the only hereditary saturated subsets of $E^0$ are $\emptyset$ and $E^0$, and every cycle in $E$ has an exit.   (Of course the chronology here is reversed: the historically-significant Simplicity Theorem precedes the establishment of Theorem   \ref{structureofidealsrowfinite} by almost a decade.)    This is seen quite readily.   By Theorem \ref{structureofidealsrowfinite}, any ideal of $L_K(E)$ looks like $I(H,\{p_c(c)\}_{c\in C})$.  By hypothesis there are only two possibilities for $H$.  When $H = E^0$ then $C(H)$, and therefore $P(C(H))$, is empty, so that the only ideal of this form is $I(E^0) = L_K(E)$.  On the other hand,  when $H = \emptyset$, then, as by hypothesis every cycle in $E$ has an exit, we get that $C(H)$, and therefore $P(C(H))$, is empty here as well.   So the only ideal of this second form is $I(\emptyset) = \{0\}$, and the Simplicity Theorem follows.

Returning yet again to the Toeplitz graph $\mathcal{T}$ of Example \ref{example:Toeplitz}, we see as a consequence of Theorem \ref{structureofidealsrowfinite} that the complete set of ideals of $L_K(\mathcal{T})$ consists of the three graded ideals $I(\emptyset) = \{0\}$, $I(w)$, and $I(E^0) = L_K(\mathcal{T})$, together with the nongraded ideals of the form $I( w, p(e))$, where $p(x) \in K[x]$ is a polynomial of degree at least $1$ for which $p(0) \neq 0$.

Considering the stronger (admittedly unstated) form of Theorem \ref{structureofidealsrowfinite}, a second consequence (also a statement of type $(\dagger \dagger)$)  is the following description of the Leavitt path algebras satisfying the chain conditions on two-sided ideals.

\begin{proposition}\label{chainconditions}
Let $E$ be a row-finite graph and $K$ any field.

(1)    $L_K(E)$ has the descending chain condition on two-sided ideals if and only if $E$ satisfies Condition (K), and the descending chain condition holds in the lattice of hereditary saturated subsets of $E^0$  (\cite[Theorem 3.9]{AbramsBellColakRanga}).

(2)     $L_K(E)$ has the ascending chain condition on two-sided ideals if and only if $L_K(E)$ has the ascending chain condition on graded two-sided ideals, if and only if the ascending chain condition holds in the lattice of hereditary saturated subsets of $E^0$.  In particular, the Leavitt path algebra $L_K(E)$ for every finite graph $E$ has the a.c.c. on two-sided ideals  (\cite[Theorem 3.6]{AbramsBellColakRanga}).

\end{proposition}

\bigskip




  \begin{center}
   Discussion:    {\bf The Rosetta Stone.} \label{RosettaStone}
     \end{center}
     
Of great interest   in the study of Leavitt path algebras is the observation that many of the results in the subject seem to (quite  mysteriously) mimic  corresponding results for graph C$^*$-algebras.   For example, comparing the Simplicity Theorem for Leavitt path algebras (Theorem \ref{thm:AAPsimple}) with the Simplicity Theorem for graph C$^*$-algebras (Theorem \ref{Cstarsimplicity}), we see that the conditions on $E$ which yield simplicity of the associated graph algebra are identical in both cases.   Suffice it to say that the proofs of the two Simplicity Theorems utilize significantly different tools one from the other.  More to the point, even with the close relationship between $L_\C(E)$ and $C^*(E)$ in mind (cf. Proposition \ref{densesubalgebra}), it is currently not understood  as to whether either one of the Simplicity Theorems  should ``directly" imply the other.   

We provide in Appendix 1 a list of additional situations  in which an algebraic property of $L_{\mathbb{C}}(E)$ is analogous to a topological property of $C^*(E)$, and for which the necessary and sufficient graph-theoretic property of $E$ is identical in each case. 
A systematic reason which would explain the existence of  so many such examples is usually referred to as the ``Rosetta Stone of Graph Algebras".     A good reference which contains in one place a discussion of both  Leavitt path algebra and graph C$^*$-algebra properties is \cite{Malaga2006}.  We note that even the seemingly most basic of questions, ``if $L_\C(E) \cong L_\C(F)$ as rings, is $C^*(E)\cong C^*(F)$ as C$^*$-algebras?" (and its converse), has only been answered (in the affirmative) for restricted classes of graphs; the question in general remains open (see \cite{AT}).     The search for the Rosetta Stone comprises one of the many current lines of research in the field.    

\subsection{Matrix rings over the Leavitt algebras}

There are too many additional ``internal / multiplicative" properties of Leavitt path algebras to include them all in this article.     For a number of reasons (its connection to the Rosetta Stone and its important consequences outside of Leavitt path algebras, to name two), we spend some space here describing  the Isomorphism Question for Matrix Rings over Leavitt algebras.    

\medskip

We reconsider the Leavitt algebras $L_K(1,n)$ for $n\geq 2$, the motivating examples of Leavitt path algebras.   Fix $n$ and $K$, and let $R$ denote $L_K(1,n)$.   By construction we have ${}_RR \cong {}_RR^n$ as left $R$-modules; so by taking endomorphism rings and using the standard representation of these endomorphism rings as matrix rings, we get $R \cong {\rm M}_n(R)$ as $K$-algebras.     Indeed, since ${}_RR \cong {}_RR^{1+j(n-1)}$ for all $j\in \N$,  we similarly get $R \cong {\rm M}_{1+j(n-1)}(R)$ as $K$-algebras for all $j\in \N$.   Now starting from a different point of view:  once we have established a ring isomorphism $S \cong {\rm M}_\ell(S)$ for some ring $S$ and some $\ell \in \N$, by taking $\ell \times \ell$ matrix rings of both sides $t$ times, we get    $S \cong {\rm M}_{\ell^t}(S)$ for any $t\in \N$.  In particular,  we have $R \cong {\rm M}_{n^t}(R)$ for all $t\in \N$; indeed, using the previous observation, we have more generally  that $R \cong {\rm M}_{(1+j(n-1))^t}(R)$ as $K$-algebras for all $j,t\in \N$.      

The question arises:   if $R = L_K(1,n)$ is isomorphic as $K$-algebras to some $p \times p$ matrix ring over itself, must $p$ be an integer of the form   $(1+j(n-1))^t $?   It is not hard to give an example where the answer is negative:  one can show (by explicitly writing down matrices which multiply correctly) that $R = L_K(1,4)$ has $R \cong {\rm M}_2(R)$, and $2$ is clearly not of the indicated form when $n=4$.          But an analysis of this particular case leads easily to the observation that if $d \ | \ n^t$ for some $t\in \N$, then $R \cong  {\rm M}_d(R)$ (by an explicitly described isomorphism). 

The upshot of the previous observations is the natural question:  $$\mbox{Given } n\in \N, \mbox{ for which } d\in \N \mbox{ is }L_K(1,n) \cong {\rm M}_d(L_K(1,n)) \mbox{ as $K$-algebras?}$$
\noindent
The analogous question was posed for matrix rings over the Cuntz algebras $\mathcal{O}_n$ in \cite{PaschkeSalinas}:  given $ n\in \N$, for which $ d\in \N $ is  $ \mathcal{O}_n \cong {\rm M}_d( \mathcal{O}_n)$ as C$^*$-algebras?   The resolution of this analogous question required many years of effort.  In the end, the solution may be obtained as a consequence of the Kirchberg Phillips Theorem:
$ \mathcal{O}_n \cong {\rm M}_d( \mathcal{O}_n)$ if and only if ${\rm g.c.d.}(d,n-1) = 1$.   So while the C$^*$-algebra question was resolved for matrices over the Cuntz algebras, the  solution did not shed any light on the analogous Leavitt algebra question, both because  the C$^*$-solution  required analytic tools, and because it did not produce an explicit isomorphism between the germane algebras. 

An easy consequence of  \cite[Theorem 5]{Leav62} is that, when  ${\rm g.c.d.}(d,n-1) > 1$,  then $L_K(1,n) \not\cong {\rm M}_d(L_K(1,n)) $.  With this and the Cuntz algebra result in hand, it made sense to conjecture that $L_K(1,n) \cong {\rm M}_d(L_K(1,n)) $ if and only if ${\rm g.c.d.}(d,n-1) = 1$.   Clearly if $d \ | \ n^t$ for some $t\in \N$ then ${\rm g.c.d.}(d,n-1) = 1$, so that the conjecture is validated in this situation.  The key idea was to  explicitly produce an isomorphism in situations more general than this.    The method of attack was clear:  one reaches the desired conclusion by finding a subset of ${\rm M}_d(L_K(1,n))$ of size $2n$ which both behaves as in $(\dagger)$, and generates ${\rm M}_d(L_K(1,n))$ as a $K$-algebra.  

  The smallest pair $d,n$ for which ${\rm g.c.d.}(d,n-1) = 1$ but $d  \not\vert \ n^t$ for any $t\in \N$ is the case $d = 3, n=5$.   Finding a subset of ${\rm M}_3(L_K(1,5))$ of size $2\cdot 5 = 10$ which behaves as in $(\dagger)$ is not hard; for instance, by (somewhat) mimicking the process used in the $d \ | \ n^t$ case, one is led to consider these five matrices in ${\rm M}_3(L_K(1,5))$  
$$\begin{pmatrix}x_{1}&0&0\\
    x_2&0&0\\
    x_3&0&0\end{pmatrix}, \hspace{.05in}
 \begin{pmatrix}x_4&0&0\\
    x_5&0&0\\
    0&1&0\end{pmatrix}, \hspace{.05in}
 \begin{pmatrix}0&0&{x_1}^2\\
    0&0&x_2x_1\\
    0&0&x_3x_1\end{pmatrix}, \hspace{.05in}
  \begin{pmatrix}0&0&x_4x_1\\
    0&0&x_5x_1\\
    0&0&x_2\end{pmatrix}, \hspace{.05in}
 \begin{pmatrix}0&0&x_3\\
    0&0&x_4\\
    0&0&x_5\end{pmatrix}, \hspace{.05in}  $$
\normalsize
together with their ``duals" 
$$\begin{pmatrix}y_{1}&y_2&y_3\\
    0&0&0\\
    0&0&0\end{pmatrix}, \hspace{.05in}
 \begin{pmatrix}y_4&y_5&0\\
    0&0&1\\
    0&0&0\end{pmatrix}, \hspace{.05in}
 \begin{pmatrix}0&0&0\\
    0&0&0\\
    {y_1}^2&y_1y_2&y_1y_3\end{pmatrix}, \hspace{.05in}
  \begin{pmatrix}0&0&0\\
    0&0&0\\
    y_1y_4&y_1y_5&y_2\end{pmatrix}, \hspace{.05in}
 \begin{pmatrix}0&0&0\\
    0&0&0\\
    y_3&y_4&y_5\end{pmatrix}. \hspace{.05in} $$
\noindent
Although these ten matrices satisfy $(\dagger)$, they do not generate all of ${\rm M}_3(L_K(1,5))$ (in retrospect, one can show that these ten matrices do not generate the matrix unit $e_{1,3}$, for example).     

The breakthrough came from a process which 
involves viewing matrices over Leavitt algebras as Leavitt path algebras for various graphs, and then manipulating the underlying graphs appropriately.    This process led to the consideration of the following (very similar, yet) different set of five matrices in ${\rm M}_3(L_K(1,5))$

$$\begin{pmatrix}x_{1}&0&0\\
    x_2&0&0\\
    x_3&0&0\end{pmatrix}, \hspace{.05in}
 \begin{pmatrix}x_4&0&0\\
    x_5&0&0\\
    0&1&0\end{pmatrix}, \hspace{.05in}
 \begin{pmatrix}0&0&{x_1}^2\\
    0&0&x_2x_1\\
    0&0&x_3x_1\end{pmatrix}, \hspace{.05in}
  \begin{pmatrix}0&0&x_4x_1\\
    0&0&x_5x_1\\
    0&0&x_2\end{pmatrix}, \hspace{.05in}
 \begin{pmatrix}0&0&x_4\\
    0&0&x_3\\
    0&0&x_5\end{pmatrix}, \hspace{.05in} $$
\normalsize
together with their duals
$$\begin{pmatrix}y_{1}&y_2&y_3\\
    0&0&0\\
    0&0&0\end{pmatrix}, \hspace{.05in}
 \begin{pmatrix}y_4&y_5&0\\
    0&0&1\\
    0&0&0\end{pmatrix}, \hspace{.05in}
 \begin{pmatrix}0&0&0\\
    0&0&0\\
    {y_1}^2&y_1y_2&y_1y_3\end{pmatrix}, \hspace{.05in}
  \begin{pmatrix}0&0&0\\
    0&0&0\\
    y_1y_4&y_1y_5&y_2\end{pmatrix}, \hspace{.05in}
 \begin{pmatrix}0&0&0\\
    0&0&0\\
    y_4&y_3&y_5\end{pmatrix}. \hspace{.05in} $$
\noindent
The only differences between the two sets of ten matrices lie in the fifth and tenth matrices, where two of the entries have been interchanged.  It is now not hard to show that this second set of ten matrices satisfies $(\dagger)$, {\it and} generates ${\rm M}_3(L_K(1,5))$ as a $K$-algebra.    The underlying idea which prompted the interchange of entries is purely number-theoretic, and is fully described in Appendix 2.   In short, the integer $3$ is used to partition the set $\{1,2,3,4,5\}$ into the subsets $\{1,4\} \sqcup \{2,3,5\}$;  then, in order to build the first five matrices of this second set, one  inserts monomials having left-most factor $x_t$ into row $i$ in such a way that $i$ and $t$ are in the same subset with respect to this partition.  So putting the term $x_3$ in row 1 and $x_4$ in row 2 (as is done in the fifth matrix of the first displayed set)  will not work; on the other hand, putting $x_4$ in row 1 and $x_3$ in row 2 is consistent with this partition, and leads to a collection with the desired properties.    Once this observation was made, the generalization to arbitrary $d,n$ was not overly difficult.

\begin{theorem}\label{thm:matricesoverLeavittalgebras}\cite[Theorems 4.14 and 5.12]{AbAnP}
Let $n\in \N$ and let $K$ be any field.   Then 
$$L_K(1,n) \cong {\rm M}_d(L_K(1,n)) 
\ \Longleftrightarrow \ \ {\rm g.c.d.}(d,n-1)=1.$$
More generally, $${\rm M}_d(L_K(1,n)) \cong {\rm M}_{d'}(L_K(1,n))  
 \ \Longleftrightarrow \ \ {\rm g.c.d.}(d,n-1)={\rm g.c.d.}(d',n-1).$$

\noindent
Moreover, when ${\rm g.c.d.}(d,n-1)={\rm g.c.d.}(d',n-1)$,  an isomorphism  ${\rm M}_d(L_K(1,n)) \rightarrow {\rm M}_{d'}(L_K(1,n))  $ can be explicitly described.  
\end{theorem}

There are two historically important consequences of the explicit construction of the isomorphisms which yield Theorem \ref{thm:matricesoverLeavittalgebras}.  
First, this context is one of the few places where a result from one side of the graph algebra universe yields a result in the other.  Specifically, when $K = \C$ and ${\rm g.c.d.}(d,n-1)=1$, the explicit nature of an isomorphism $L_\C(1,n) \cong {\rm M}_d(L_\C(1,n))$ constructed in the proof of Theorem \ref{thm:matricesoverLeavittalgebras} allows (by a straightforward completion process) for the {\it explicit} construction of an isomorphism  $\mathcal{O}_n \cong {\rm M}_d(\mathcal{O}_n)$.   (The description of such an explicit isomorphism came as more than a bit of a surprise to some researchers in the C$^*$-community.)
Second, the explicit construction led to the resolution of a longstanding question in group theory.   In the mid 1970's, G. Higman produced, for each pair  $r,n \in \N$ with $n\geq 2$, an infinite, finitely presented simple group, denoted $G^+_{n,r}$.  (The groups  $G^+_{n,r}$ are called the {\it Higman-Thompson} groups.) Higman was able to establish some sufficient conditions regarding isomorphisms between these groups, but did not have a complete classification. 
However, in 2011, Enrique Pardo showed how the construction given in the proof of Theorem \ref{thm:matricesoverLeavittalgebras} could be brought to bear in this regard. 
\begin{theorem}\label{Higmansimplegroups}   \cite[Theorem 3.6]{Pardo2011}  $$G_{n,r}^+ \cong G_{m,s}^+ \ \  \Longleftrightarrow \ \ m=n \ \mbox{ and }  \ {\rm g.c.d.}(r,n-1) = {\rm g.c.d.}(s,n-1).$$ 
\end{theorem}
{\it Sketch of Proof.} The ($\Longrightarrow$) direction was already known by Higman.  Conversely,  one first shows that $G_{n,\ell}^+$ can be realized as an appropriate subgroup of the invertible elements of $M_\ell(L_\mathbb{C}(1,n))$ for any $\ell \in \N$.    Then  one verifies that the explicit isomorphism  from $M_r(L_\mathbb{C}(1,n))$ to $M_s(L_\mathbb{C}(1,n))$ provided in the proof of  Theorem \ref{thm:matricesoverLeavittalgebras} takes $G_{n,r}^+$ onto $G_{n,s}^+$.   \hfill $\Box$

\medskip

For any three positive integers $t,n,r$ (with $n\geq 2$), Brin \cite{Brin}  constructed a group (denoted $tV_{n,r}$)  which can be viewed as a $t$-dimensional analog of the Higman-Thompson group,  in that $1V_{n,r} \cong G_{n,r}^+$. (The groups  $tV_{n,r}$ are called the {\it Brin-Higman-Thompson} groups.)  On the other hand, for $t$ a positive integer and $n\geq 2$, one may consider the $t$-fold tensor product algebra $L_K(1,n)^{\otimes t}$ of $L_K(1,n)$ with itself $t$ times.  (We will more fully consider such tensor products in the following subsection.)     In \cite{DMP}, Dicks and Mart\'{i}nez-P\'{e}rez beautifully generalize Pardo's Theorem \ref{Higmansimplegroups} by showing  that  $tV_{n,r}$ is isomorphic to an appropriate subgroup of the invertible elements of ${\rm M}_r(L_K(1,n)^{\otimes t})$ (specifically, the positive unitaries), and subsequently use this isomorphism to establish that   $tV_{n,r} \cong t'V_{n',r'}$ if and only if $t=t'$, $n=n'$, and ${\rm g.c.d.}(r, n-1) =  {\rm g.c.d.}(r', n'-1).$  Along the way, Dicks and Mart\'{i}nez-P\'{e}rez present a streamlined, somewhat more intuitive proof of  Theorem \ref{thm:matricesoverLeavittalgebras}.

\subsection{Tensor products of Leavitt path algebras}\label{tensorsubsection}

Of fundamental importance in the theory of graph C$^*$-algebras is the fact that $\mathcal{O}_2 \otimes \mathcal{O}_2 \cong \mathcal{O}_2$ (homeomorphically).   This isomorphism is not explicitly described; rather, it follows (originally) from some deep work done by Elliott (and streamlined in  \cite{Rordam}).  The isomorphism  $\mathcal{O}_2 \otimes \mathcal{O}_2 \cong \mathcal{O}_2$  is utilized in the proof of the Kirchberg Phillips Theorem.   (The C$^*$-algebra $\mathcal{O}_2$ is nuclear, so that there is no ambiguity in forming this tensor product.)

In the context of the previous paragraph, together with the Rosetta Stone discussion,  it is then natural to ask:   is $L_K(1,2) \otimes_K L_K(1,2) \cong L_K(1,2)$?       This question had been posed as early as 2006, and was the focus of sustained investigative effort for a number of years.   
The resolution of this question {\it in the negative} came in early 2011, in the form of three  different approaches by three different investigative teams.      

The first proof (unpublished), offered by Warren Dicks, utilized a classical result of Cartan and Eilenberg \cite[Theorem X1.3.1]{CartanEilenberg}, which yields that  the flat dimension of a tensor product is at least the sum of the flat dimensions of the two algebras.   By Theorem \ref{AMPfinite}$^\prime$, the global dimension of $L_K(1,2)$ (indeed, of any Leavitt path algebra) is at most $1$. ({\it Global dimension at most }$1$ is equivalent to {\it hereditary}.) Consequently, the flat dimension of a Leavitt path algebra $L_K(E)$ equals $1$ precisely when $L_K(E)$ is not von Neumann regular (i.e., when there are $L_K(E)$-modules which are not flat).  But it had been shown in \cite{AbRanga} that if $E$ is a graph containing at least one cycle, then $L_K(E)$ is not von Neumann regular, so, in particular, $L_K(1,2)\cong L_K(R_2)$ is not von Neumann regular.   So the flat dimension, and therefore  also the global dimension, of $L_K(1,2) \otimes_K L_K(1,2)$ is at least $2$, so that $L_K(1,2)\otimes L_K(1,2)$ cannot be a Leavitt path algebra (again using Theorem \ref{AMPfinite}$^\prime$), and so can't be isomorphic to $L_K(1,2)$.

A second proof (unpublished) was offered by Jason Bell and George Bergman.  Effectively, Bell and Bergman explicitly constructed a left $L_K(1,2)$-module $M$ (involving functions on $[0,1) \subseteq \mathbb{R}$ having finite support in $\mathbb{Q}$ of the form $n / 2^j$), and showed that the left $L_K(1,2) \otimes L_K(1,2)$-module $M \otimes M$ has projective dimension $2$, so that $L_K(1,2) \otimes L_K(1,2)$ has global dimension at least $2$, and thus (arguing as did Dicks)  cannot be isomorphic to any Leavitt path algebra.

The third approach to verifying that $L_K(1,2) \otimes L_K(1,2) \not\cong L_K(1,2)$ is the most general of the three.   Utilizing Hochschild homology, Ara and Corti\~{n}as in \cite{AraCortinas} showed (among many other things) the following, from which the result of interest follows immediately.  

\begin{theorem}\label{AraCortinas}
Suppose  $\{E_i\}_{i=1}^m$ and $\{F_j\}_{j=1}^n$ are finite graphs, each containing at least one cycle, and let $K$ be any field.  If $\otimes_{i=1}^m L_K(E_i)$ is Morita equivalent to $ \otimes_{j=1}^n L_K( F_j)$, then $m=n$. 
\end{theorem}

 Two currently unresolved questions about the tensor products of Leavitt path algebras will be given in Section \ref{section:CurrentLines}.

\subsection{Some additional internal / multiplicative properties of Leavitt path algebras}\label{additionalinternalsubsection}

We conclude the section by presenting five additional multiplicative properties of Leavitt path algebras:  primeness; the center; Gelfand Kirillov dimension; wreath products; and the simplicity of the corresponding bracket Lie algebra.  

\smallskip

A ring $R$ is {\it prime} in case for any two-sided ideals $I,J$ of $R$, if $IJ = \{0\}$ then $I=\{0\}$ or $J = \{0\}$.    
A graph $E$ is called {\it downward directed} if, for any two vertices $v,w\in E^0$, there exists a vertex $u\in E^0$ for which $v\geq u$ and $w\geq u$.    

\begin{theorem}\label{primeness}\cite{APPSMIndiana} 
Let $E$ be any graph and $K$ any field.  Then $L_K(E)$ is prime if and only if $E$ is downward directed. 
\end{theorem}

{\it Sketch of Proof.}  ($\Rightarrow$) If $R$ denotes $L_K(E)$, and  $v,w \in E^0$, then the ideals $RvR$ and $RwR$ are each nonzero, so that $RvRRwR \neq \{0\}$, so that $vRw \neq \{0\}$, which yields a nonzero element of the form $\alpha \beta^*$ with $s(\alpha) = v$, and $r(\beta^*) = s(\beta)=w$, so that $u = r(\alpha)$ has the desired property.   

($\Leftarrow$)   The converse can be proved `elementwise', but it is easier to invoke \cite[Proposition 5.2.6(1)]{NvO}, which implies that for a $\Z$-graded ring, primeness is equivalent to graded primeness.  So we need only check that if $I,J$ are nonzero {\it graded} ideals, then $IJ \neq \{0\}$.   But by Proposition \ref{prop:latticeiso} (or its generalization Theorem \ref{GradIdgen} given below), any nonzero graded ideal contains a vertex; so if $v\in E^0 \cap I$ and $w\in E^0 \cap J$, and $u\in E^0$ with $v\geq u$ and $w\geq u$, then $0 \neq u = u^2 \in IJ$.   \hfill $\Box$

\bigskip

For a ring $R$, the {\it center} $Z(R) = \{ r \in R \ | \ rx = xr \mbox{ for all } x \in R\}$.   It is well-known that $Z({\rm M}_n(K)) = K\cdot I_n$ (where $I_n$ denotes the identity matrix in ${\rm M}_n(K)$).   Additionally, this easily yields  that  the center of ${\rm M}_\N(K)$ is  $\{0\}$.    The following result includes these observations as specific cases.

\begin{theorem}\label{center}(\cite{ArandaPinoCrow})
Let $E$ be a row-finite graph.     Suppose $L_K(E)$ is simple (see Theorem \ref{thm:AAPsimple}).    If $E^0$ is finite, then $Z(L_K(E)) = K\cdot 1_{L_K(E)}$.   If  $E^0$ is infinite, then $Z(L_K(E)) = \{0\}$.
\end{theorem}

For a $K$-algebra $A$, the Gelfand-Kirillov Dimension $GK_K(A)$ is an algebraic invariant of $A$ which, loosely speaking,  measures how far $A$ is from being finite dimensional.  (Finite dimensional algebras have GK dimension $0$.   On the other hand, the free associative $K$-algebra on two generators has GK dimension $\infty$. Such an algebra is said to have {\it exponential growth}; otherwise, the algebra has {\it polynomially bounded growth}.  See e.g. \cite{KrauseLenagan} for a full description.)  
If $C$ and $C'$ are two disjoint cycles (i.e., ${\rm Vert}(C) \cap {\rm Vert}(C') = \emptyset$), the symbol $C \Rightarrow C'$ indicates that there is a path which starts in ${\rm Vert}(C)$ and ends in ${\rm Vert}(C')$.      A sequence of disjoint cycles $C_1, ..., C_k$ is a {\it chain of length $k$} in case $C_1 \Rightarrow \cdots \Rightarrow C_k$.  Let $d_1$ denote the maximal length of a chain of cycles in $E$, and let $d_2$ denote the maximal length of a chain of cycles each of which has an exit.  
\begin{theorem}\label{GKdimension}\cite[Theorem 5]{AAJZ}
Let $E$ be a finite graph and $K$ any field.

(1)  $L_K(E)$ has exponential growth if and only if there exist a pair of distinct cycles in $E$ which are not disjoint.

(2)  In case $L_K(E)$ has polynomially bounded growth, then the GK dimension of $L_K(E)$ is ${\rm max}(2d_1 - 1, 2d_2)$.
\end{theorem}

Further results regarding Leavitt path algebras of polynomially bounded growth, and of the automorphism groups of some specific such algebras, are presented in \cite{AAJZ-PNAS}.

\bigskip

For a countable dimensional $K$-algebra $C$ and  ring-theoretic property $\mathcal{P}$,  an {\it affinization of $A$ with respect to $\mathcal{P}$}  is an embedding of $C$ in a finitely generated (i.e., affine) $K$-algebra $D$, for which, if $C$ has $\mathcal{P}$, then so does $D$.         

Let $E$ be a row-finite graph and $A$ any associative $K$-algebra.   In \cite{AA2014}, the authors present the construction of the {\it wreath product}, denoted  $ A \ wr \ L_K(E)$.    In case $W$ is a hereditary saturated subset of $E^0$, then the wreath product construction allows for the realization of $L_K(E)$ as the wreath product of two Leavitt path algebras, namely, as $L_K(W) \ wr \ L_K(E/W)$.     Furthermore, let $\mathcal{T}$ be the Toeplitz graph of Example \ref{example:Toeplitz}.  Then the  wreath product  $A \ wr \ L(\mathcal{T})$  is isomorphic to a $K$-algebra  of the form $K[x , x^{-1}] + {\rm M}_\N(A)$ (with multiplication explicitly described).   This algebra  can then be embedded in an algebra of the form $K[x , x^{-1}] + {\rm RCFM}_\N(A)$, where ${\rm RCFM}_\N(A)$ is the (unital) ring of those $\N \times \N$ matrices with entries in $A$, for which each row and each column contains at most finitely many nonzero entries.   One may then build in a natural way an affine  $K$-algebra $B$, generated by four elements, for which $K[x , x^{-1}] + {\rm M}_\N(A) \subset B \subset K[x , x^{-1}] + {\rm RCFM}_\N(A)$.   

\begin{theorem}\label{affinization}\cite{AA2014}   For an associative $K$-algebra $A$ let $B$ be the affine $K$-algebra described above.

(1)  There exists a unital algebra $A$ for which $B$ is an affinization of $K[x , x^{-1}] + {\rm M}_\N(A)$ with respect to the property {\rm non-nil Jacobson radical}.

(2)  There exists a unital  algebra $A$ for which  $B$ is an affinization of $K[x , x^{-1}] + {\rm M}_\N(A)$ with respect to the property {\rm non-nilpotent locally nilpotent radical}. 
\end{theorem}
\vspace{-.05in}
\ \ \ \ \ \ Both of the constructs mentioned in Theorem \ref{affinization} give a systematic approach to what had been previously longstanding ring-theoretic questions. 

\bigskip

For a $K$-algebra $R$, the corresponding {\it bracket Lie algebra} $[R,R]$ consists of $K$-linear combinations of elements of the form $xy - yx$ with $x,y\in R$.    $[R,R]$ is a Lie algebra, with the usual bracket operation.   A Lie algebra $L$ is called {\it simple} in case $[L,L] \neq \{0\}$, and the only Lie ideals of $L$ are $\{0\}$ and $L$.   
 Let $E$ be a finite graph, and write $E^0 = \{v_i \mid 1 \leq i \leq m\}$.  If  $v_i$ is a not a sink, for each $1\leq j \leq m$  let $a_{ij}$ denote the number of edges $e \in E^1$ such that $s(e) = v_i$ and $r(e) = v_j$. In this situation,  define $B_i = (a_{ij}) - \epsilon_i \in \Z^{m}$ (where $\epsilon_i$ is the element of $\Z^{m}$ which is $1$ in the $i$-th coordinate, and zero elsewhere).     On the other hand, if $v_i$ is a sink, let $B_i = (0) \in \Z^{m}.$ 
\begin{theorem} \label{simplebracket} \cite[Theorem 23]{AbramsMesyan} 
Let $K$ be a field, and let $E$ be a finite  graph having at least two vertices for which $L_K(E)$ is simple. Write $E^0 = \{v_1, \dots, v_m\}$, and for each $\, 1\leq i \leq m$ let $B_i$ be as above. Then the Lie $K$-algebra $\, [L_K(E),L_K(E)]$ is simple if and only if $\, (1,\dots, 1) \not\in \mathrm{span}_K\{B_1, \dots, B_m \}$.
\end{theorem}

As it turns out, the condition given in Theorem \ref{simplebracket} for the simplicity of $\, [L_K(E),L_K(E)]$ depends not only on the structure of $E$ but also on the characteristic of $K$ (see  \cite[Examples 28 and 29]{AbramsMesyan}).   The $K$-dependence of a result  about Leavitt path algebras is very  much the exception.  But for one intriguing additional example, see  Theorem \ref{nofreequotients} and the subsequent discussion.   

By introducing and utilizing the notion of a {\it balloon} over a subset of $E^0$, Alahmedi and Alsulami are able to extend Theorem \ref{simplebracket} to all row-finite graphs (specifically, the simplicity of $L_K(E)$ is not required); see \cite[Theorem 2]{AA-Lie}.   For instance, it is shown in \cite{AA-Lie}  that the graph $E$ given here

$$
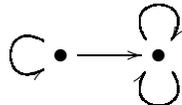
\xymatrix{ \bullet \ar@(ul,dl) \ar[r] & \bullet \ar@(ul,ur) \ar@(dr,dl)}$$

\bigskip
\noindent
has the property that the Lie algebra $[L_K(E), L_K(E)]$ is simple, even though the Leavitt path algebra $L_K(E)$ is not simple.    

In related work \cite{AA-Skew}, the same two authors analyze the simplicity of the Lie algebra of $\ast$-skew-symmetric elements of a Leavitt path algebra.

%
%

\section{Module-theoretic properties of Leavitt path algebras}\label{ModuleSection}

The module theory of Leavitt path algebras has for the most part been focused on the structure of the finitely generated projective $L_K(E)$-modules, owing to the Ara / Moreno / Pardo Realization Theorem (Theorem \ref{AMPfinite}) describing $\mathcal{V}(L_K(E))$.   In this section we take a closer look at the structure of these projectives, specifically, the purely infinite modules.  Of central interest here   is the question of whether or not the analog of the Kirchberg Phillips Theorem (Theorem \ref{KPforgraphalgebras}) holds for Leavitt path algebras; we present in Theorem \ref{restrictedKP} the Restricted Algebraic KP Theorem.  We next look at the structure of some simple (non-projective) $L_K(E)$-modules.
We conclude by considering some monoid-theoretic properties of $\mathcal{V}(L_K(E))$.      

\smallskip

%
%

\subsection{Purely infinite simplicity}

We have seen that the cycle structure of the graph $E$, and the existence of exits for those cycles,  is a significant factor driving the algebraic structure of the Leavitt path algebra $L_K(E)$.  
We have also seen behavior in the Leavitt algebras $L_K(1,n)$  that at first glance seems somewhat exotic:   ${}_RR \cong {}_RR^n$ as left $R$-modules.    Specifically, the module ${}_RR$ has the property that ${}_RR \cong {}_RR \oplus P$ where $P \neq \{0\}$; i.e., that ${}_RR$ has a nontrivial direct summand which is isomorphic to itself.  ({\it Nontrivial} here means that the complement of the direct summand is nonzero.)   This same sort of behavior is manifest in $L_K(E)$ when $E$ has cycles with exits.  

\begin{remark} Let $\alpha\in {\rm Path}(E)$, and let $r(\alpha) = w$.  Then $L_K(E)\alpha\alpha^* \cong L_K(E)w$ as left $L_K(E)$-modules, since,  if  $\varphi = \rho_{\alpha}: L_K(E)\alpha\alpha^* \rightarrow L_K(E)w$ denotes right multiplication by $\alpha$, then it is easy to show that $\varphi^{-1} = \rho_{\alpha^*}$.
\end{remark}

\begin{proposition}\label{cycleswithexitsgivepis}
Suppose  $c$  is a cycle in a graph $E$ based at a vertex $v$,   and suppose  $e$  is an exit for $c$  with  $s(e) = v$.  Then the left $L_K(E)$-module $L_K(E)v$ has a nontrivial direct summand isomorphic to itself.  
\end{proposition}

{\it Proof.} 
Clearly  $L_K(E)v = L_K(E)cc^* + L_K(E)(v - cc^*)$.  But the sum is direct:  if $xcc^* = y(v-cc^*)$ for $x,y\in L_K(E)$, then multiplying  both sides on the right by $cc^*$ yields $xcc^* = y(cc^* - cc^*) = 0.$            That    $L_K(E)v \cong L_K(E)cc^*$ as left $L_K(E)$-modules follows from the previous Remark.     Since $e$ is an exit for $c$ we have $e^*c = 0$ by (CK1).   Now to show that  the complement $L_K(E)(v - cc^*)$ is nonzero,   assume to the contrary that $v-cc^* = 0$.  Multiplying both sides on the left by $ee^*$ gives   $ee^* - 0 = 0$, thus giving $ee^* = 0$,  which is impossible.   \hfill $\Box$  

\medskip

\medskip

A left $R$-module $M$ is called {\it  infinite} in case $M \cong M \oplus N$ with $N \neq \{0\}$.   An idempotent $x \in R$ is called {\it  infinite} in case $Rx$ is  infinite.  The ring  $R$ is called {\it purely infinite simple} in case $R$ is simple, and each nonzero left ideal of $R$ contains an  infinite idempotent.  Purely infinite simple rings were first introduced in \cite{AGPpis};  the idea was born in the context of C$^*$-algebras.
Clearly a purely infinite module can satisfy neither of the two chain conditions, nor can it have finite uniform dimension.

With the Simplicity Theorem in hand, and with Proposition \ref{cycleswithexitsgivepis} as guidance,  some medium-level effort yields the following.

\begin{theorem}\label{pistheorem}{\bf The Purely Infinite Simplicity Theorem} \cite{AAPpureinf}     Let $E$ be a  row-finite graph and $K$ any field.   Then $L_K(E)$ is purely infinite simple if and only if $L_K(E)$ is simple, and $E$ contains at least one cycle.   
\end{theorem}

Theorems \ref{finiteacyclic} and \ref{pistheorem} together yield what is typically called the {\it Dichotomy for simple Leavitt path algebras}:   for $L_K(E)$ simple, either $L_K(E)$ is purely infinite simple, or $L_K(E) \cong {\rm M}_n(K)$ for some $n\in \N$.  

\medskip

In the context of Leavitt path algebras, the purely infinite simple algebras play an especially intriguing role.   For any ring $S$, the Grothendieck group $K_0(S)$ is the universal group corresponding to the abelian monoid $\mathcal{V}(S)$.  (Here {\it universal} means that any homomorphism from $\mathcal{V}(S)$ to an abelian group $G$ necessarily factors through $K_0(S)$.)  When $\mathcal{V}(S) \cong \Z^+$ (as is often the case in general, e.g.,  when $S$ is a field or $S = \Z$), then one gets $K_0(S) \cong \Z$, by ``adding in the negatives".   As it turns out, however, if $S$ is purely infinite simple, then $\mathcal{V}(S) \setminus \{[0]\}$ is a group, precisely $ K_0(S)$.    This is perhaps counterintuitive at first glance:  although  $\mathcal{V}(S)$ has an identity element (namely, $[0]$), there still remains an identity element once $[0]$ is eliminated.    For instance, if $R = L_K(1,n)$, then $R \oplus R^{n-1} = R^n \cong R$.   Using this, it's trivial to conclude that $[R^{n-1}]$ is an identity element for  $\mathcal{V}(R) \setminus \{[0]\} = \{ [R], [R^2], \dots, [R^{n-1}]\}$.  The group $\mathcal{V}(R) \setminus \{[0]\} $ is clearly isomorphic to $ \Z / (n-1)\Z$.  

Although the converse is not true for arbitrary rings, when one restricts to the class of Leavitt path algebras, then the converse is true as well \cite{PardoMadVet}:   $L_K(E)$ is purely infinite simple if and only if $\mathcal{V}(L_K(E)) \setminus \{[0]\}$ is a group (necessarily $K_0(L_K(E))$).  Moreover, this group is easy to describe in this situation.   
As is standard, for a finite graph $E$ with $E^0 = \{ v_1, v_2, \dots, v_n\}$, the {\it incidence matrix} of $E$ is the $| E^0 | \times | E^0|$ $\Z^+$-valued matrix $A_E$, where $A_E(i,j)$ equals the number of edges $e$ for which $s(e) = v_i$ and $r(e) = v_j$.   
By interpreting the (CK2) relation as it plays out in $\mathcal{V}(L_K(E))$, one gets
\begin{proposition}\label{prop:Kzeroisoquotient}
Let $E$ be a finite graph, with $|E^0| = n$.  Suppose $L_K(E)$ is purely infinite simple.  Then 
$$K_0(L_K(E)) \cong  \Z^n /(I_n-A_E)\Z^n,$$
where $I_n$ denotes the $n\times n$ identity matrix.   
\end{proposition}
In other words, when $L_K(E)$ is purely infinite simple, then $K_0(L_K(E))$ is the cokernel of the linear transformation $I_n-A_E: \Z^n \rightarrow \Z^n$ induced by matrix multiplication. 

As an easy example of how this plays out in an already-familiar situation, suppose $E = R_m$, the graph with one vertex and $m$ loops.  
Then $A_E = (m)$, so $I - A_E$ is the $1 \times 1$ matrix  $(1-m)$, and $K_0(L_K(E)) \cong  \Z^1 /(1-m)\Z^1 = \Z / (m-1)\Z$, as we've seen previously.  

\bigskip

%
%

\subsection{Towards a Classification Theorem for purely infinite simple Leavitt path algebras}



In many endeavors in  which  an object from one class is associated to an object in another,  a fundamental question is to identify the stalks of the process; that is, determine which objects from the first class correspond to the same object in the second.     Asked in the current context:  if two graphs $E, F$ produce the ``same" Leavitt path algebra (up to isomorphism, or up to Morita equivalence, or up to some other ring-theoretic invariant), can anything be said about the relationship between $E$ and $F$?   As seen in Theorem \ref{finiteacyclic}, if $E$ and $F$ are finite acyclic graphs for which $L_K(E) \cong L_K(F)$,  then $E$ and $F$ have the  same number of sinks, and the same number of directed paths ending at those sinks.    (An additional easy consequence of Theorem \ref{finiteacyclic} is that if $L_K(E)$ is Morita equivalent to $L_K(F)$, then $E$ and $F$ have the same number of sinks.)  

We spend some time here investigating this question in the context of purely infinite simple Leavitt path algebras.    The reason is twofold:   this investigation plays up an important relationship between Leavitt path algebras and symbolic dynamics, and also provides the foundation for much of the current research focus in Leavitt path algebras.    The discussion here will be quite broad and intuitive;  for details, the standard reference is \cite{LindMarcus}.

For a finite directed graph $E$, one defines the notion of a ``flow" (essentially, ``flow of information") through the graph.   Two graphs $E$ and $F$ are ``flow equivalent" in case the collection of flows through $E$ match up appropriately with the collection of flows through $F$.  Two matrices with entries in $\Z^+$ are called {\it flow equivalent} in case the directed graphs corresponding to the two matrices are flow equivalent.     The directed graph $E$ (or the corresponding incidence matrix $A_E$) is called 
 
(1) \emph{irreducible} if  for any pair $v,w\in E^0$ 
there is a path from $v$ to $w$;

(2) \emph{essential} if there are neither sources nor sinks in $E$; and 

(3) \emph{trivial} if $E$ consists of a single cycle with no other
vertices or edges.

\medskip
\noindent
A deep, fundamental result in flow dynamics is 

{\bf Franks' Theorem.}  \cite{Franks} 
Suppose that $A$ and $B$ are non-negative irreducible essential  nontrivial square integer
matrices. 
Then $A$ and $B$ are flow equivalent if and only if

\smallskip

$
\Z^n /(I_n-A)\Z^n \cong \Z^m /(I_m-B)\Z^m  \ \mbox{ and } \
  \det(I_n-A) = \det(I_m-B).
$
\smallskip

\medskip

There are a number of ways to systematically modify a directed graph.   As an intuitive example,   {\it expansion at} $v$ modifies the graph $E$ to the graph $E_v$ as indicated here.  

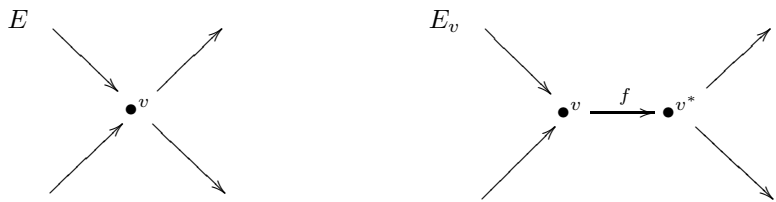
\footnotesize
 $$  
E \ \xymatrix{ {} \ar[rd]  &   {}  &  {} \\
{} & \bullet^v      {} \ar[ru] \ar[rd] & {}  \\
{} \ar[ru] & {}  & {}  }
 \hspace{1in}
E_v \  \xymatrix{ {} \ar[rd]  &  {} & {}  &  {} \\
{} & \bullet^v   \ar[r]^{f} & \bullet^{v^*}   {} \ar[ru] \ar[rd] & {}  \\
{} \ar[ru] & {}  & {} & {} }
$$
\normalsize
It can easily be shown that the graphs $E$ and $E_v$ are flow equivalent.     In a similar manner, one may describe five more systematic modifications of a graph (each having the property that the original graph is flow equivalent to the modified graph):  {\it contraction} (the inverse of expansion);  {\it out-split},  as well as its inverse {\it out-amalgamation};  and {\it in-split},  as well as its inverse {\it in-amalgamation}.   The specific descriptions of these ``graph moves" are given in Appendix 3. 

\smallskip
 The second deep, fundamental theorem germane to the current discussion is 

\smallskip

{\bf The Parry / Sullivan Theorem.}    Two finite directed graphs are flow equivalent if and only if one can be gotten from the other by a sequence of transformations involving these six graph moves.   







\smallskip

Combining Franks' Theorem with the Parry / Sullivan Theorem, we get 

\begin{theorem}\label{FranksandPStheorem}
Suppose $E$ and $F$ are irreducible essential nontrivial graphs.  Then 
$\Z^n /(I_n-A_E)\Z^n \cong \Z^m /(I_m-A_F)\Z^m
\ \mbox{and}  \ {\rm det}(I - A_E) = {\rm det}(I - A_F)$ 
 if and only if  $E$ can be obtained from $F$ by some sequence of graph moves, with each move one of the six types described above.  
\end{theorem}

We are now in position to present the (miraculous?) bridge between the ideas from flow dynamics and those of Leavitt path algebras.   First, using the Purely Infinite Simplicity Theorem (Theorem \ref{pistheorem}) and some straightforward graph theory,  it is not hard to show that  $E$ is irreducible, essential, and nontrivial if and only if $E$ has no sources and $L_K(E)$ is purely infinite simple.     Next,

\begin{proposition}\label{graphmovespreserveequivalence}  Suppose $E$ is a graph for which $L_K(E)$ is purely infinite simple.   Suppose $F$ is gotten from $E$ by doing one of the six aforementioned graph moves.    Then $L_K(E)$ and $L_K(F)$ are Morita equivalent.  In particular, $L_K(F)$ is purely infinite simple.    In addition, if $v$ is a source in $E$, and $F$ is gotten from $E$ by eliminating $v$ and all edges $e\in E^1$ having $s(e) = v$, then  $L_K(E)$ and $L_K(F)$ are Morita equivalent.   
\end{proposition}

{\it Sketch of Proof.}  It is not hard to show that an isomorphic copy of $L_K(E)$ can be viewed as a (necessarily full, by simplicity) corner of  $L_K(G)$ (or vice-versa), where $E$ and $G$ are related by one of the graph moves.      \hfill $\Box$ 

\medskip

The previous discussion yields the first of two desired results.

\begin{theorem}\label{MoritaKP}    Let $E$ and $F$ be finite graphs and $K$ any field.  Suppose $L_K(E)$ and $L_K(F)$ are purely infinite simple.      If 
$$K_0(L_K(E))  \cong K_0(L_K(F))  \ \ \  \mbox{and} \ \ \ {\rm det}(I - A_E) = {\rm det}(I - A_F),$$ then $L_K(E)$ and $L_K(F)$ are Morita equivalent.   
\end{theorem}

{\it Sketch of Proof.}   Suppose $E$ and/or $F$ have sources; then using   Proposition  \ref{graphmovespreserveequivalence} we may construct graphs $E'$ and $F'$ for which  $L_K(E')$ and $L_K(F')$ are purely infinite simple,  $L_K(E) $ is Morita equivalent to $L_K(E')$, and $L_K(F)$ is Morita equivalent to $ L_K(F')$, where $E'$ and $F'$ have no sources.   But since Morita equivalent rings have isomorphic $K_0$ groups, and because  (it's straightforward to show that) ${\rm det}(I - A_E) = {\rm det}(I - A_{E'})$ and ${\rm det}(I - A_F) = {\rm det}(I - A_{F'})$, we have that the  hypotheses of Theorem \ref{FranksandPStheorem} are satisfied for $E'$ and $F'$.  Thus   $F'$ can be gotten from $E'$ by a sequence of appropriate graph moves.  But again invoking Proposition  \ref{graphmovespreserveequivalence}, each of these moves preserves Morita equivalence.  So $L_K(E')$ is Morita equivalent to $L_K(F')$,   and the result follows.   \hfill $\Box$

\medskip

The third deep, fundamental result of interest here is 

\smallskip
{\bf Huang's Theorem.}  Suppose $L_K(E)$ is Morita equivalent to $L_K(F)$.   Further, suppose there is {\it some} isomorphism $\varphi: K_0(L_K(E)) \rightarrow K_0(L_K(F))$ for which $\varphi([L_K(E)]) = [L_K(F)]$.     
Then there is some Morita equivalence $\Phi: L_K(E){\rm -Mod} \rightarrow L_K(F){\rm -Mod}$ for which $\Phi |_{K_0(L_K(E))} = \varphi.$

\smallskip

Consequently:

\begin{theorem}\label{restrictedKP}  {\bf The Restricted Algebraic Kirchberg Phillips Theorem.}  \cite[Corollary 2.7]{Flow}    Let $E$ and $F$ be finite graphs and $K$ any field.  Suppose $L_K(E)$ and $L_K(F)$ are purely infinite simple.     If 
$$K_0(L_K(E))  \cong K_0(L_K(F)) 
 \mbox{ via an isomorphism for which } [L_K(E)] \mapsto   [L_K(F)],$$
 \vspace{-.25in}
$$ \ \ \  \mbox{and} \ \ \ {\rm det}(I - A_E) = {\rm det}(I - A_F),$$ then $L_K(E) \cong L_K(F)$. 
\end{theorem}  

{\it Sketch of Proof.}    For any Morita equivalence $\Phi: R-Mod \rightarrow S-Mod$, if $\Phi({}_RR) = {}_SS$, then $R \cong {\rm End}_R({}_RR) \cong {\rm End}_S(\Phi({}_RR)) \cong {\rm End}_S({}_SS) \cong S$ as rings.  Now apply Theorem \ref{MoritaKP} together with Huang's Theorem.   \hfill $\Box$

\medskip

As an example of how the Restricted Algebraic KP Theorem can be implemented, let $E$ be the graph

 \vspace{-.2in} 
\footnotesize
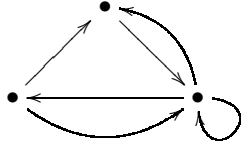
$$ \xymatrix{ {} & \bullet \ar[rd] & & {} \\
\bullet  \ar[ru]
\ar@/^{-15pt}/ [rr]&  & \bullet
\ar@(r,d)
 \ar[ll] \ar@/^{-10pt}/
[lu] & }$$
\normalsize

\vspace{.15in}
\noindent
Then using the description provided in Proposition \ref{prop:Kzeroisoquotient}, we get $K_0(L_K(E)) \cong \mathbb{Z} / 3\mathbb{Z} $;  moreover, under this isomorphism, $[L_K(E)] \mapsto 1$.  Easily we get 
     ${\rm det}(I - A_{E}) = -3 <0$.    But the Leavitt path algebra $L_K(R_4) \cong L_K(1,4)$ has precisely the same data associated with it, so we conclude that 
$L_K(E) \cong L_K(1,4).$

\medskip

In Section \ref{section:CurrentLines} we describe how the Restricted Algebraic Kirchberg Phillips Theorem has been acting as a springboard for much of the current research energy in the subject.

%
%

\subsection{Simple $L_K(E)$-modules}   

We now move our focus on $L_K(E)$-modules from  projectives to simples.   

\smallskip

Let $p$ be an {\it infinite path in} $E$; that is, $p$ is a sequence $ e_1e_2e_3\cdots$, where $e_i \in E^1$ for all $i\in \N$, and for which $s(e_{i+1}) = r(e_i)$ for all $i\in \N$.   (N.b.:  an infinite path in $E$ is not an element of ${\rm Path}(E)$,  nor  of the Leavitt path algebra $L_K(E)$.)   The set of infinite paths in $E$ is denoted by $E^\infty$.  
For $p = e_1e_2e_3\cdots \in E^\infty$  and $n\in \N$,  
$p_{>n}$ denotes the infinite path $e_{n+1}e_{n+2}\cdots$.

Let $c$ be a closed path in $E$.  Then  $c c c \cdots$ is an infinite path in $E$,  denoted by $c^\infty$, and called a    {\it cyclic infinite} path.   A closed path $d$ is  {\it irreducible} in case $d$ cannot be written as $e^j$ for any  closed path $e$ and $j > 1$.   For any closed path $c$ there exists an irreducible $d$  for which $c = d^n$; then $c^\infty = d^\infty$ as elements of $E^\infty$.

For $p,q \in E^\infty$, $p$ and $q$ are {\it tail equivalent} (written $p \sim q$) in case there exist integers $m,n$ for which $p_{>m}= q_{>n}$ (i.e., in case $p$ and $q$ eventually become the same infinite path).   For $p\in E^\infty$,  $[p]$ denotes the $\sim$ equivalence class of $p$.   
 An element $p$ of $E^\infty $ is {\it rational} in 
case $p \sim c^\infty$ for some irreducible closed path $c$; otherwise   $p$ is {\it irrational}.    For instance, in   $$R_2 \ = \ \xymatrix{\bullet^v \ar@(ul,dl)_e \ar@(ur,dr)^f} \ ,$$
\noindent
$q = efeffefffeffffe\cdots$ is an irrational infinite path.  In any graph $E$ for which there exists a vertex having two distinct irreducible closed paths based at that vertex, it is not hard to show that there are uncountably many irrational infinite paths in $E^\infty$.  Additionally,  there are infinitely many irreducible paths in such a situation (and thus infinitely many tail-inequivalent infinite rational paths);  for instance, any path of the form $ef^i$ for $i\in \Z^+$ is irreducible in  $R_2$.


\begin{definition} \label{Chendef}
{\rm 
Let $p$ be an infinite path in the graph $E$, and let $K$ be any field.  Let $V_{[p]}$ denote the $K$-vector space having basis $[p]$,  consisting of the distinct elements of $E^\infty$ which are tail-equivalent to $p$.    For  $v\in E^0$, $e\in E^1$, and $q = f_1f_2f_3\cdots  \in [p]$, define 
$$
 v \cdot q =
 \begin{cases}
q  &\text{if }  v = s(f_1) \\
0 &\text{otherwise,} 
\end{cases}
    \ \ \  \ \ \ \ 
 e\cdot q=
  \begin{cases}
eq  &\text{if }  r(e) = s(f_1) \\
0 &\text{otherwise,} 
\end{cases}
$$
$$\ \ \ \ \mbox{and} \qquad 
e^* \cdot q = 
 \begin{cases}
\tau_{>1}(q)  &\text{if }  e = f_1 \\
0 &\text{otherwise.} 
\end{cases}
$$
\noindent
Then  the $K$-linear extension of this action to all of $V_{[p]}$ gives a left $L_K(E)$-module structure on  $V_{[p]}$.   
}
\end{definition}

\begin{theorem} \label{Chentheoremforsimples}  (\cite[Theorem 3.3]{Chen}).  Let $E$ be any graph and $K$ any field.  Let $p\in E^\infty$.   Then the left $L_K(E)$-module $V_{[p]}$ described in  Definition \ref{Chendef}
 is simple.   Moreover, if $p,q \in E^\infty$, then $V_{[p]} \cong V_{[q]}$ as left $L_K(E)$-modules if and only if $p \sim q$, which happens precisely when $V_{[p]} = V_{[q]}$.
\end{theorem}

A module of the form $V_{[p]}$ as  in Theorem \ref{Chentheoremforsimples} is called a {\it Chen simple $L_K(E)$-module}.   In  \cite{AraRanga2},  Ara and Rangaswamy  describe those Leavitt path algebras $L_K(E)$  which admit at most countably many  simple left modules (Chen simples or otherwise) up to isomorphism.  Building on an observation made prior to Definition \ref{Chendef}, one sees that the structure of $K$ plays a role in this result, in that when $L_K(E)$ has this property, and $E$ contains cycles, then necessarily $K$ must be countable.

   It is possible to explicitly describe projective resolutions for the Chen simple modules.
 Let $\beta = e_1e_2 \cdots e_n \in {\rm Path}(E)$ or $\beta =  e_1e_2 \cdots \in E^\infty$.  
 For each $i\geq 0$ (and $i\leq n-1$ if $\beta = e_1 \cdots e_n \in {\rm Path}(E)$),  let 
$X_{i}(\beta) = \{f \in E^1 \ | \ s(f) = s(e_{i+1}),  \ \mbox{and} \ f \neq e_{i+1}\},$   and let  $J_i(\beta)$ be the left ideal $\sum_{ f \in X_i(\beta)} L_K(E) f^* \beta_i^*$ of $L_K(E)$.  
The following explicit description of projective resolutions of Chen simple modules follows from an elementwise analysis of the kernel of the appropriate right-multiplication map.  
(For an element $m$ in a left $L_K(E)$-module $M$, and any left ideal $I$ of $L_K(E)$, $\rho_m : I \rightarrow M$ denotes right multiplication by $m$.)  

\begin{theorem} \label{projresofVcinftyfromL(E)v}  \cite{AbManTo}    Let $E$ be any graph and $K$ any field.  

(1)  Let $c$ be an irreducible  closed path in $E$, with $v=s(c)$.  Then  $V_{[c^\infty]}$ is finitely presented (in fact, singly presented);  a projective resolution of $V_{[c^\infty]}$ is given by
$$ \xymatrix{ 0 \ar[r] &  L_K(E)v \ar[r]^{\rho_{c-v}} & L_K(E)v \ar[r]^{  \ \ \rho_{c^\infty}} &  V_{[c^\infty]} \ar[r] & 0}.$$

(2)  Let $p \in E^\infty$ be an  irrational infinite path in $E$ for which no element of  ${\rm Vert}(p)$ is an infinite emitter.  Then  $$ \xymatrix{ 0 \ar[r] & \oplus_{i=0}^\infty J_i(p) \ar[r] & L_K(E)v \ar[r]^{ \ \ \ \rho_p} &  V_{[p]} \ar[r] & 0}$$
 is a projective resolution of $V_{[p]}$.   In particular, $V_{[p]}$ is finitely presented if and only  if $X_i(p)$ is nonempty for at most 
finitely many $i\in \Z^+$.  
\end{theorem}

 Theorem \ref{projresofVcinftyfromL(E)v}  sharpens and clarifies some of the results of \cite{AraRanga1}.  
The explicit description of projective resolutions given in Theorem \ref{projresofVcinftyfromL(E)v}   can be used to (easily) show that $V_{[c^\infty]}$ is never projective, and that $V_{[p]}$ (for $p$ irrational) is not projective when  $V_{[p]}$ is not finitely presented (e.g., whenever $E$ is a finite graph).      Consequently, these two types of modules admit nontrivial extensions, some of which are captured in the following result.

\begin{theorem}\label{Ext(ST)} \cite{AbManTo}   Let $E$ be a finite graph and $K$ any field.  Let $T$ be a Chen simple module.   Denote by $U(T)$ the set
$ \{v\in E^0 \ | \ vT \neq \{0\} \}.$  For $p\in E^\infty$, denote by $r(X_i(p))$ the set $\{r(e_i) \ | \ e_i \in X_i(p)\}$.  

(1)  Let $d$ be an irreducible  closed path in $E$ with $v = s(d)$.
Then 
  
  \noindent
   ${\rm Ext}^1_{L_K(E)}(V_{[d^\infty]}, T)\neq \{0\}$ if and only if 
   $vT \neq \{0\}.$

(2)     Let $p$ be an irrational infinite path in $E$.   Then ${\rm Ext}^1_{L_K(E)}(V_{[p]},T) \neq \{0\} $ if and only if 
 $r(X_i(p)) \cap U(T) \neq \emptyset$ for infinitely many $i\geq 0$. 
\end{theorem}

As a consequence of Theorem \ref{Ext(ST)},  whenever $E$ is a graph containing at least one cycle, then (non-projective) indecomposable $L_K(E)$-modules of any desired finite length  can be constructed.

We close this subsection on simple $L_K(E)$-modules by noting that Rangaswamy \cite{Rangasimples} has given a construction of such modules arising from the  infinite emitters $v$ of $E^0$.

\bigskip

%
%

\subsection{Additional module-theoretic properties of $L_K(E)$}
\label{Additionalmodulestuff}

The previous discussion in this section first focused on projective modules, then on non-projective simple modules, over Leavitt path algebras.   We conclude the section by mentioning some monoid-theoretic properties of $M = \mathcal{V}(L_K(E))$.
 As the $\mathcal{V}$-monoid of a ring, $M$ is of course conical, and contains a distinguished element (as described prior to Theorem  \ref{BergmansTheorem}).     But there are two important additional properties of $\mathcal{V}(L_K(E))$, both of which yield information about the decomposition of projective $L_K(E)$-modules.

Suppose that $M$ is a left $R$-module which admits two direct sum decompositions $M = A_1\oplus A_2 = B_1 \oplus B_2$.   We ask whether there is necessarily some relationship between the two decompositions, indeed, whether there is some compatible ``refinement" of  these which allows for the systematic  formation  of each of the summands.   More formally, suppose $A_1\oplus A_2 = B_1\oplus B_2$ as left $R$-modules.  Then a {\it refinement} of this pair of direct sums consists of   left $R$-modules $M_{11}, M_{21}, M_{12},$ and $M_{22}$, for which: 
$$A_1 = M_{11} \oplus M_{12}, \  \ A_2 = M_{21} \oplus M_{22},$$
$$ \ B_1 = M_{11} \oplus M_{21}, \ \  B_2 = M_{12} \oplus M_{22}.$$

A second type of decomposition of modules relates to cancellation of direct summands.  Clearly in general an isomorphism $A \oplus C \cong B \oplus C$ of left $R$-modules need not imply $A \cong B$.  A  germane example here is this:   if $R = L_K(1,n)$, and $A = \{0\}$, $B = {}_RR^{n-1}$, and $C =  {}_RR$, then we have $A \oplus C \cong B \oplus C$ (since ${}_RR \cong {}_RR^{n}$), but obviously $A \not\cong B$.     In various situations it is natural to   require a stronger relationship  between such isomorphic direct sums, prior to trying to cancel $C$. One possible approach is as follows.  
 A ring $R$ is called  {\em separative} in case it satisfies the following
property: If $A,B,C\in \mathcal{V}(R)$ satisfy $A\oplus C\cong B\oplus C$, 
and $C$ is isomorphic to direct summands of both $A^n$ and $B^n$ for
some $n\in \mathbb N$, then $A\cong B$.   (Note that this additional condition obviously renders moot the previous example.)

\begin{theorem}
\label{refinementandseparativity} 
Let $E$ be a row-finite graph and $K$ any field.

1)  \cite[Proposition 4.4]{AMP07}  The monoid $M_E$  is a refinement monoid.   Consequently, $\mathcal{V}(L_K(E))$ is a refinement monoid.  

2)  \cite[Theorem 6.3]{AMP07} The monoid $M_E$ is separative. Consequently, the monoid $\mathcal{V}(L_K(E))$, and thus the ring $L_K(E)$,  is separative. 
\end{theorem}

{\it Sketch of Proof.}   (1)   is established by a careful analysis of the generators and relations which produce the graph monoid $M_E$.    On the other hand, (2) 
follows in part from results of Brookfield
\cite{Brookfield} on {\it primely generated} refinement monoids.  \hfill $\Box$ 

\medskip

In fact, the class of primely generated refinement
monoids satisfies many other nice cancellation properties, e.g. {\it unperforation}.   We will revisit refinement monoids at the end of  Section \ref{section:CurrentLines}.

%
%

\section{Classes of algebras related to, or motivated by,  Leavitt path algebras of row-finite graphs}\label{arbitrarygraphssection}

\smallskip

Historically, Leavitt path algebras were first defined only in the context of row-finite graphs.   Subsequently, the more general definition of  Leavitt path algebras for countable graphs (\cite{AAPHouston}), and then truly arbitrary graphs (\cite{GoodearlLimits}), appeared in the literature.   The original notion of a Leavitt path algebra for row-finite graphs has been generalized in other ways as well, including:  the construction of Leavitt path algebras for separated graphs;  Cohn path algebras;  Kumjian-Pask algebras of higher ranks graphs; Leavitt path rings; and more.    In this section we give an overview of some of these Leavitt-path-algebra-inspired structures.


\subsection{Leavitt path algebras for arbitrary graphs.}   Suppose $E$ is a graph which contains an infinite emitter $v$; that is, the set $s^{-1}(v) = \{e \in E \ | \ s(e) = v\}$ is infinite.   Then in a purely ring-theoretic context, the symbol $\sum_{e \in s^{-1}(v)}ee^*$, which would be the natural generalization of the (CK2) relation imposed at $v$,  is not defined.   Even in the analytic context of graph C$^*$-algebras, where convergence properties might allow for some sort of appropriate interpretation of an infinite sum, an expression of the form  $\sum_{e \in s^{-1}(v)}s_es_e^*$ proves to be problematic, in part owing to the fact that $\{s_es_e^* \ | \ e \in s^{-1}(v) \}$ is an infinite set of orthogonal projections.   

So, somewhat cavalierly, we simply choose not to invoke any (CK2)-like relation at infinite emitters.    We recall that a vertex $v\in E^0$ is  {\it regular} in case $0 < |s^{-1}(v)| < \infty$.

\begin{definition}\label{def:Lpaarbitrary}
Let $E = (E^0, E^1, s, r)$ be any  graph, and $K$ any field.  Let $\widehat{E}$ denote the extended graph of $E$.   The {\it Leavitt path $K$-algebra $L_K(E)$} is defined as the path $K$-algebra $K\widehat{E}$, 
modulo the relations:

\smallskip
(CK1) \ $e^*e'=\delta_{e,e'}r(e)$ for all $e, e' \in E^1$. 

 (CK2) \ $v=\sum_{\{e\in E^1|s(e)=v\}}ee^*$ for every  regular vertex $v\in E^0$. 
 
 \smallskip

Equivalently, we may define $L_K(E)$ as the free associative $K$-algebra on generators $E^0 \sqcup E^1 \sqcup (E^1)^*$, modulo the relations 

(1)  $v{v'} = \delta_{v,v'}v$ for all $v,v' \in E^0.$

(2)  $s(e)e = er(e) = e$ for all $e\in E^1.$

(3)  $r(e)e^* = e^*s(e) = e^*$ for all $e\in E^1.$

(4)  $e^*e' = \delta_{e,e'}r(e)$ for all $e,e' \in E^1$.

(5)  $v = \sum_{\{e\in E^1|s(e)=v\}}ee^*$ for every regular $v\in E^0$. 

\end{definition}

So the definition of a Leavitt path algebra for arbitrary graphs is essentially word-for-word identical to that for row-finite graphs (since ``regular" and ``non-sink" are identical properties in the row-finite case);  there is simply no (CK2) relation imposed at any vertex which is  the source vertex of infinitely many edges.

The generalization from Leavitt path algebras of row-finite graphs to those of arbitrary graphs was achieved in two stages.  Owing to the hypotheses typically placed on the corresponding graph C$^*$-algebras (in order to ensure separability), the initial extension for Leavitt path algebras was to graphs having countably many vertices and edges.   It is shown in \cite{AAPHouston} that the Leavitt path algebra of any such countable graph is Morita equivalent to the Leavitt path algebra of a suitably defined row-finite graph, using the {\it desingularization} process.   Subsequently, the foundational results regarding   Leavitt path algebras for arbitrary graphs were  presented in \cite{GoodearlLimits}.   Among other things, Goodearl established a suitable definition and context for morphisms between graphs (so-called {\it CK-morphisms}).  He was then able to show that direct limits exist in the appropriately defined graph category (denoted {\bf CKGr}), and that the functor $L_K$ from {\bf CKGr} to the category of $K$-algebras preserves direct limits.

The generalization to Leavitt path algebras of arbitrary graphs (from those of row-finite graphs) indeed expands the Leavitt path algebra universe.   For instance, it was shown in \cite{AbramsRanga}  that $L_K(E)$ is Morita equivalent to $L_K(F)$ for some row-finite graph $F$ if and only if $E$ contains no uncountable emitters (i.e., in case the set $s^{-1}(v)$ is at most countable for each $v\in E^0$).    So, for instance, let $I$ be an uncountable set, 
and let $D_I$ denote the graph consisting of two vertices $v,w$, and edges $\{e_i \ | \ i\in I\}$, where $s(e_i) = v$ and $r(e_i) = w$.  
Then $L_K(D_I)$ is isomorphic to the (unital) $K$-algebra generated by ${\rm M}_I(K) \sqcup \{ {\rm Id} \}$, where ${\rm Id}$ is the $I \times I$ identity matrix.  So $L_K(D_I)$ is not Morita equivalent (let alone,  isomorphic) to the Leavitt path algebra of any row-finite graph.   Similarly, if $R_c$ denotes the ``rose with uncountably infinitely many petals" graph, then $L_K(R_c)$ is not Morita equivalent to $L_K(F)$ for any row-finite graph $F$.   


In this expanded universe of Leavitt path algebras for arbitrary graphs, many of the results established in the row-finite case generalize verbatim, but many do not.   One of the main differences is that in the general case, we may pick up many new idempotents inside $L_K(E)$ for which there are no counterparts in the row-finite case.   For instance, let $v\in E^0$, and let $e\in s^{-1}(v)$.   Then the element $x = v - ee^*$ of $L_K(E)$ is easily shown to be an idempotent.  If $v$ is a regular vertex, then $x = \sum_{f \in s^{-1}(v), f\neq e}ff^*$ by the (CK2) relation.  On the other hand, if $v$ is an infinite emitter, then $x$ has no such analogous representation.   

We recall the graph-theoretic ideas given in Notation \ref{notation:heredandsat}:   a subset  $X$ of $E^0$ is {\it hereditary} in case, whenever $v\in X$ and $w\in E^0$ and $v\geq w$, then $w\in X$;   $X$ is  {\it saturated} in case, whenever $v\in E^0$ is regular and $r(s^{-1}(v))\subseteq X$, then $v\in X$.

\begin{definition}\label{BreakingVertex}
{\rm
Let $E$ be any graph, and let $H$ be a hereditary subset  of $E^0$.  A vertex  $v \in E^0$ is a  {\it breaking vertex of} $H$ in case  $v$ is in the set
$$ B{_H}=\{v\in E^0\setminus H \ \vert \ |s^{-1}(v)| = \infty  \ \text{and} \ 0 < \vert s^{-1}(v) \cap r^{-1}(E^0\setminus H)\vert < \infty\}.$$

\noindent
In words, $B_H$ consists of those vertices which are infinite emitters, which do not belong to $H$, and for which the ranges of the edges they emit are all, except for a finite (but nonzero) number, inside $H$.
For $v\in B{_H}$, define  $$v^H \ = \  v-\sum_{e\in s^{-1}(v) \cap r^{-1}(E^0\setminus H) }ee^\ast,$$ and, for any subset $S\subseteq B{_H}$, define
  $S^H=\{v^H \ \vert \ v\in S\}$.
}
\end{definition}
 
Of course a row-finite graph contains no breaking vertices, so that this concept does not play a role in the study of Leavitt path algebras arising from such graphs.  Also, we note that both $B_{E^0}$ and $B_\emptyset$ are empty.   
To help clarify the concept of breaking vertex, we offer the following example.
\begin{example}\label{InfiniteClockExample}
Let  $C_{\N}$ be the  {\it infinite clock graph} 
  pictured here
 $$  
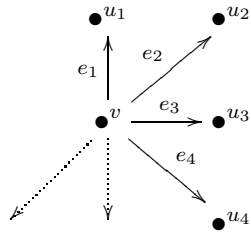
\xymatrix{ & {\bullet}^{u_1} & {\bullet}^{u_2} \\  & {\bullet}^v  
  \ar[u]^{e_1}  \ar[ur]^{e_2}  \ar[r]^{e_3}  \ar[dr]^{e_4}  \ar@{.>}[d] \ar@{.>}[dl] 
\ar@{}[dl]  & {\bullet}^{u_3} \\ &  & {\bullet}^{u_4}}$$
Let $U$ denote the set $\{u_i \ \vert \ i\in \N\} = C_\N^0 \setminus \{v\}$.  Any subset of $U$ is a hereditary subset of $C_\N^0$.  We note also that, since saturation applies only to regular vertices, any subset of $U$ is saturated as well.

If $H \subseteq U$ has $U\setminus H$ infinite, or if $H = U$, then $B_H = \emptyset$.   On the other hand, if $ U \setminus H$ is finite, then $B_H = \{v\}$, and in this situation, $v^H = v - \sum_{\{i  |  r(e_i) \in U \setminus H\}} e_ie_i^*$.    
\end{example}

It is clear that  for any  hereditary  saturated subset $H$ of a graph $E$,  and  for any $S\subseteq B_H$, the ideal
 $I(H\cup S^H)$ is a graded ideal, as it is generated by elements of $L_K(E)$ of degree zero.  It turns out that this process generates all the graded ideals of $L_K(E)$. 
 We denote by  ${\mathcal L}_{gr}(L_K(E))$   the collection of two-sided graded ideals of $L_K(E)$, and by $\mathcal{T}_E$ the collection of pairs $(H,S)$ where $H$ is a hereditary saturated subset of $E$,  and   $S\subseteq B_H$.   

\begin{theorem}  \label{GradIdgen}\cite[Theorem 5.7]{TomfordeUniqueness} 
Let $E$ be an arbitrary graph and $K$ any field. Then
there is a bijection 
$
\varphi:   \mathcal{L}_{gr}(L_K(E))  \to  \mathcal{T}_E, \   \mbox{given by } \ 
 I  \mapsto  (I \cap E^0, S)
 $
 where $ S=\{v\in B_H\ \vert \ v^H \in I\}$ for $H= I \cap E^0$. The inverse is given by 
$
\varphi^{-1}:  \mathcal{T}_E \to         \mathcal{L}_{gr}(L_K(E)),  \  \mbox{via}  \  
 (H, S) \mapsto  I(H\cup S^H).
$
\end{theorem}

There is an appropriate lattice structure which can be defined in $\mathcal{T}_E$ so that the map $\varphi$ is a lattice isomorphism.   In addition, there is a generalization of Theorem \ref{GradIdgen} to the lattice of {\it all} ideals of $L_K(E)$, see \cite[Theorem 2.8.10]{TheBook}.



We close the subsection by presenting a result which is of interest in its own right (it provided a systematic approach to answering a decades-old question of Kaplansky), and which will reappear later in the context of the Rosetta Stone.     An algebra $A$ is called {\it left primitive} in case $A$ admits a faithful simple left module.   It was shown in \cite{APPSMIndiana} that for row-finite graphs, $L_K(E)$ is primitive if and only if $E$ is downward directed and satisfies Condition (L).  
However, the extension of this result to arbitrary  graphs requires an extra condition.   The graph $E$ has the {\it Countable Separation Property} in case there exists a countable set $S\subseteq E^0$ with the property that for every $v\in E^0$ there exists $s\in S$ for which $v\geq s$.  

\begin{theorem}\label{primitiveLpa}\cite[Theorem  5.7]{AbBellRanga}
Let $E$ be an arbitrary graph and $K$ any field.  Then $L_K(E)$ is primitive if and only if $E$ is downward directed, $E$ satisfies Condition (L), and $E$ has the Countable Separation Property.
\end{theorem}

\subsection{Leavitt path algebras of separated graphs}

The (CK2) condition imposed at any regular vertex in the definition of a Leavitt path algebra may be modified in various ways.  Such is the motivation for the discussion in both this and the following subsection.   All of these ideas appear in \cite{AraGoodearl}. 

\smallskip

 In the (CK2) condition which appears in the definition of the Leavitt path algebra $L_K(E)$, the edges emanating from a given regular vertex $v$ are treated as a single entity, and the relation $v = \sum_{e \in s^{-1}(v)} ee^*$ is imposed.  More generally, one may partition the set  $s^{-1}(v)$ into disjoint nonempty subsets, and then impose a (CK2)-type relation corresponding exactly to those subsets.   More formally, 
a {\it separated graph} is a pair $(E,C)$, where $E$ is a graph, $C = \sqcup_{v\in E^0}C_v$, and, for each $v\in E^0$, 
 $C_v$ is a partition of $s^{-1}(v)$ (into pairwise disjoint nonempty subsets). (In
case $v$ is a sink, $C_v$ is taken to be the empty family of subsets of $s^{-1}(v)$.)

\begin{definition}\label{def:Lpaofseparated}  Let $E$ be any  graph and $K$ any field.  Let $\widehat{E}$ denote the extended graph of $E$, and $K\widehat{E}$ the path $K$-algebra of $\widehat{E}$.    The {\it Leavitt path algebra of the separated graph $(E,C)$ with coefficients in
the field $K$} is the quotient of $K\widehat{E}$ by the ideal generated by these two types of relations: 

(SCK1) for each $X\in C$, $e^*f = \delta_{e,f}r(e)$ for all $e,f\in X$,    and

(SCK2) for each non-sink $v\in E^0$, $v = \sum_{e\in X}ee^*$ for every finite  $X\in C_v$.   

\end{definition}

So the usual Leavitt path algebra $L_K(E)$ is exactly $L_K(E,C)$, where each $C_v$ is defined to be the subset $\{s^{-1}(v)\}$ if $v$ is not a sink, and $\emptyset$ otherwise.   
Leavitt path algebras of separated graphs include a much wider class of algebras than those which arise as Leavitt path algebras in the standard construction.    For instance, the algebras of the form $L_K(m,n)$ for $m\geq 2$ originally studied by Leavitt in \cite{Leav62}  do not arise as $L_K(E)$ for any graph $E$.  On the other hand,  as shown in \cite[Proposition 2.12]{AraGoodearl},  $L_K(m,n)$ ($m\geq 2$) appears as a full corner of the Leavitt path algebra of an explicitly described separated graph (having two vertices and $m+n$ edges).   In particular, $L_K(m,n)$ is Morita equivalent to the Leavitt path algebra of a separated graph.

Of significantly more importance is the following Bergman-like realization result, which shows that the collection of  Leavitt path algebras of separated graphs is extremely broad.

\begin{theorem}\label{thm:realizationseparated}\cite[Section 4]{AraGoodearl}   Let $M$ be any conical abelian monoid.   Then there exists a graph $E$, and partition $C = \sqcup_{v\in E^0}C_v$, for which $\mathcal{V}(L_K(E,C)) \cong M$.  
\end{theorem}

Consequently, $\mathcal{V}(L_K(E,C))$ need not share the separativity nor the refinement properties of the standard Leavitt path algebras $L_K(E)$.   Furthermore, the ideal structure of $L_K(E,C)$ is in general significantly more complex than that of $L_K(E)$, 
but a description of the idempotent-generated ideals can be achieved (solely in terms of graph-theoretic information).

\subsection{Cohn path algebras}
In the previous subsection we saw one way to modify the (CK2) relation, namely, by imposing it on subsets of $s^{-1}(v)$ for $v\in E^0$.  

\smallskip

A second way to modify the (CK2) relation is to simply eliminate it.


\begin{definition}\label{def:Cohnpathalgebra}   Let $E$ be any graph and $K$ any field.   The {\it Cohn path algebra} $C_K(E)$ is the path $K$-algebra $K\widehat{E}$ of the extended graph of $E$, modulo the relation

(CK1)   $e^*f = \delta_{e,f}r(e)$  for each $e,f\in E^1$.

\end{definition}
The terminology ``Cohn path algebra" postdates the Leavitt path algebra terminology, and owes to the fact that for each $n\geq 1$, the algebra $C_K(R_n)$ (for $R_n$ the rose with $n$ petals graph) is precisely the algebra $U_{1,n}$ described and investigated by Cohn in \cite{Cohn66}.    

   Indeed, even the case $n=1$ is of interest here:  
   $C_K(R_1)$ is the unital $K$-algebra $A$ generated by an element $e$ for which $e^*e = 1$ (and no other relation involving $e$).  Thus we get that $C_K(R_1)$ is exactly the Jacobson algebra described in Example \ref{example:Toeplitz}, so that (using the computation presented in that Example), we have $C_K(R_1) \cong L_K(\mathcal{T})$, the Leavitt path algebra of the Toeplitz graph.    Pictorially,  $$C_K( \ \ \ \ \ \xymatrix{ \bullet \ar@(ul,dl)}) \ \cong \ L_K(\ \ \ \ \ \xymatrix{ \bullet \ar@(ul,dl) \ar[r] & \bullet}).$$
\noindent
This isomorphism between a Cohn path algebra and a Leavitt path algebra is not a coincidence.

\begin{theorem}\label{CohnisLeavitt}\cite[Section 1.5]{TheBook} 
Let $E$ be any graph.   Then there exists a graph $F$ (which is explicitly constructed from $E$)  for which $C_K(E) \cong L_K(F)$.   That is, every Cohn path algebra is isomorphic to a Leavitt path algebra.
\end{theorem}

In particular, the explicit construction mentioned in Theorem \ref{CohnisLeavitt} of the graph $F$ from the graph $E$ in case $E = R_1$ yields that $F = \mathcal{T}$.    So although at first glance the Cohn path algebra construction seems less restrictive than the Leavitt path algebra construction, the collection of algebras which arise as $C_K(E)$ is (properly) contained in the collection of algebras which arise as $L_K(E)$.      (One way to see that the containment is proper is to note that the Cohn path algebra $C_K(E)$ has Invariant Basis Number for any finite graph $E$; see \cite{AbramsKanuni}.)   

One may view the Leavitt path algebras and Cohn path algebras as occupying the opposite ends of a spectrum:  in the former, we impose the (CK2) relation at all (regular)  vertices, while, in the latter, we do not impose it at any of the vertices.  The expected middle-ground construction may be formalized:  if $X$ is any subset of the regular vertices ${\rm Reg}(E)$ of $E$, then the {\it Cohn path $K$-algebra relative to $X$}, denoted $C_K^X(E)$, is the algebra $C_K(E)$, modulo the (CK2) relation imposed {\it only} at the vertices $v\in X$.    So $C_K(E) = C_K^\emptyset (E)$, while $L_K(E) = C_K^{{\rm Reg}(E)}(E)$.   Theorem \ref{CohnisLeavitt} generalizes appropriately from Cohn path algebras to relative Cohn path algebras.

\subsection{Additional constructions}

We close this section with a description of four additional Leavitt-path-algebra-inspired constructions. 

\medskip

{\it Cohn-Leavitt algebras.}  \ The following (not unexpected) mixing-and-matching of the Leavitt path algebras of separated graphs with the relative Cohn path algebras has been defined and studied in \cite{AraGoodearl}.  

\begin{definition}
Let $(E,C)$ be a separated graph.  Let $C_{{\rm fin}}$ denote the subset of $C$ consisting of those $X$ for which $|X|$ is finite.  Let $S$ be any subset of $C_{{\rm fin}}$.  Denote by $CL_K(E,C,S)$ the quotient of the path $K$-algebra $K\widehat{E}$, modulo the relations (SCK1) of Definition \ref{def:Lpaofseparated}, together with the relations (SCK2) for the sets $X \in S$.   $CL_K(E,C,S)$ is called the {\it Cohn-Leavitt algebra of the triple $(E,C,S)$}.
\end{definition}

\medskip

{\it Kumjian-Pask algebras.} \  Any directed graph $E = (E^0, E^1, s, r)$ may be viewed as a category $\Gamma_E$; the objects of $\Gamma_E$ are  the   vertices $E^0$, and, for each pair $v,w\in E^0$, the morphism set ${\rm Hom}_{\Gamma_E}(v,w)$ consists of those elements of ${\rm Path}(E)$ having source $v$ and range $w$.   Composition is concatenation.      As well, the set $\Z^+$ is a category with one object, and morphisms given by the elements of $\Z^+$, where composition is addition.   In this level of abstraction, the length map $\ell: {\rm Path}(E) \rightarrow \Z^+$ is a functor, which satisfies the following {\it factorization} property:   if $\lambda \in {\rm Path}(E)$ and $\ell(\lambda) = m+n$, then there are unique $\mu, \nu \in {\rm Path}(E)$ such that $\ell(\mu) = m, \ell(\nu) = n$, and $\lambda = \mu \nu$.     
Conversely, we may view a category as the morphisms of the category, where the objects are identified with the identity morphisms.  Then any category $\Lambda$ which admits a functor $d: \Lambda \rightarrow \Z^+$ having the factorization property can be viewed as a directed graph $E_\Lambda$ in the expected way.

With these observations as motivation, one defines a {\it higher rank} graph, as follows.

\begin{definition}\label{kgraphdef}
Let $k$ be a positive integer.  View the additive semigroup $(\Z^+)^k$ as a category with one object, and view a category as the morphisms of the category, where the objects are identified with the identity morphisms.   A {\it graph of rank $k$} (or simply a {\it $k$-graph}) is a countable category $\Lambda$, together with a functor $d: \Lambda \rightarrow (\Z^+)^k$, which satisfies the factorization property:  if $\lambda \in \Lambda$ and $d(\lambda) = \overline{m} + \overline{n}$ for some $\overline{m},  \overline{n} \in (\Z^+)^k$, then there exist unique $\mu, \nu \in \Lambda$ such that $d(\mu) = \overline{m}, d(\nu) = \overline{n}$, and $\lambda = \mu \nu$.    (So the usual notion of a  graph is a $1$-graph in this more general context.)

Given any  $k$-graph $(\Lambda, d)$ and field $K$, one may define   the {\it Kumjian-Pask}  $K$-algebra $KP_K(\Lambda,d)$.  (We omit the somewhat lengthy details of the construction; see \cite{AaHCR} for the complete description.)    In case $k=1$, $KP_K(\Lambda,d)$ is the Leavitt path algebra $L_K(E_\Lambda)$.    \end{definition}
\medskip

{\it The regular algebra of a graph.} \ The following construction should be viewed not as a method to generalize the notion of Leavitt path algebra, but rather to use the properties of Leavitt path algebras as a tool to answer what at first glance seems to be an unrelated question.   The ``Realization Problem for von Neumann Regular Rings" asks whether every countable conical refinement monoid can be realized as the monoid $\mathcal{V}(R)$ for some von Neumann regular ring $R$.   It was shown in \cite{AbRanga}  that the only von Neumann regular Leavitt path algebras are those associated to acyclic graphs, so it would initially seem that Leavitt path algebras would not be fertile ground in the context of the Realization Problem.   
Nonetheless, Ara and Brustenga developed an elegant construction which provides the key connection.  Using the {\it algebra of rational power series on $E$}, and appropriate localization techniques ({\it inversion}), they showed how to construct a $K$-algebra $Q_K(E)$ with the following properties.

\begin{theorem}\label{AraBrustvNr}\cite[Theorem 4.2]{AraBrustenga}
Let $E$ be a finite graph and $K$ any field.   Then there exists a  $K$-algebra $Q_K(E)$ for which:

(1)   there is an inclusion of algebras $L_K(E) \hookrightarrow Q_K(E)$,

(2)  $Q_K(E)$ is unital von Neumann regular, and

(3) $\mathcal{V}(L_K(E)) \cong \mathcal{V}(Q_K(E))$.
\end{theorem}
Consequently, using the Realization Theorem (Theorem \ref{AMPfinite}$^\prime$), Theorem \ref{AraBrustvNr}  yields that any monoid which arises as the graph monoid $M_E$ for a finite graph $E$ has a positive solution to the Realization Problem.  This result represented (at the time) a  significant broadening of the class of monoids for which the Realization Problem had a positive solution.     The result extends relatively easily to row-finite graphs (see \cite[Theorem 4.3]{AraBrustenga}), with the proviso that $Q_K(E)$ need not be unital in that generality.

 \bigskip

{\it Non-field coefficients.} \   While nearly all of the energy expended on understand  $L_K(E)$ has focused on the graph $E$, one may also relax the requirement that the coefficients be taken from a field $K$.  For a commutative unital ring $R$ and graph $E$ one may form the {\it path ring $RE$ of $E$ with coefficients in $R$} in the expected way;   it is then easy to see how to subsequently define the  {\it Leavitt path ring $L_R(E)$ of $E$ with coefficients in $R$}.   While some of the results given when $R$ is a field do not hold verbatim in the more general setting (e.g.,  the Simplicity Theorem), one can still understand much of the structure of $L_R(E)$ in terms of the properties of $E$ and $R$; see e.g. \cite{Tom}.

\bigskip

With  these many generalizations of Leavitt path algebras  having now been noted, a comment on the extremely robust interplay between algebras and C$^*$-algebras is in order.    In some situations,  the C$^*$-ideas preceded the algebra ideas; in other situations,  the opposite; and in still others, the ideas were introduced simultaneously.  

\smallskip

\footnotesize $\bullet$ \normalsize      Leavitt   \cite{Leav62}  built the Leavitt algebras $L_\C(1,n)$  (1962); subsequently, Cuntz \cite{Cuntz}   built their C$^*$-counterparts, the Cuntz algebras $\mathcal{O}_n$  (1977).   (Cuntz's results were achieved independently from the work of Leavitt.) 
 
\footnotesize $\bullet$ \normalsize    Graph C$^*$-algebras of row-finite graphs  were then introduced in  \cite{BPRS00} (2000); these  in turn motivated the definition of Leavitt path algebras  of row-finite graphs in \cite{AAP05} and \cite{AMP07} (2005).

 \footnotesize $\bullet$ \normalsize  Graph C$^*$-algebras of countable graphs which contain infinite emitters were introduced in \cite{FLR} (2000);  these motivated the definition of Leavitt path algebras of such graphs in \cite{AAPHouston} (2006).

  \footnotesize $\bullet$ \normalsize Leavitt path algebras for arbitrary graphs were first given complete consideration in \cite{GoodearlLimits} (2009).   The initial  study of  C$^*$-algebras corresponding to arbitrary graphs appears  in \cite{AT} (2013), where this notion was utilized to give the first  systematic construction of C$^*$-algebras which are prime but not primitive.

 \footnotesize $\bullet$ \normalsize  C$^*$-algebras of higher rank graphs were formalized in \cite{KumjianPask} (2000); the corresponding Kumjian-Pask algebras were introduced in \cite{AaHCR} (2014).   
 
 \footnotesize $\bullet$ \normalsize  In the context of separated graphs,  both Leavitt path algebras {\it and} graph C$^*$-algebras of these objects were  introduced essentially simultaneously in the articles \cite{AraGoodearlCstar} (2011), and 
 \cite{AraGoodearl} (2012).

\medskip

%
%

\section{Current lines of research in Leavitt path algebras}\label{section:CurrentLines}

In the previous five sections we have given an overview of the subject of Leavitt path algebras.  In this final section we consider some of the important current research problems in the field.   For additional information, see 
``The graph algebra problem page": 

\begin{center}
www.math.uh.edu/tomforde/GraphAlgebraProblems/ListOfProblems.html
\end{center}
\noindent
This website was built and is being maintained by Mark Tomforde of the University of Houston.  

\medskip

We have previously discussed the (currently unresolved)  Rosetta Stone Question  for graph algebras.  More information about the Rosetta Stone is presented in Appendix 1.  

\medskip

\subsection{The Classification Question for purely infinite simple Leavitt path algebras, a.k.a. {\it The Algebraic Kirchberg Phillips Question}} 

We start with what is generally agreed to be the most compelling unresolved  question in the subject of Leavitt path algebras, stated concisely as: 
\begin{center}
The Algebraic Kirchberg Phillips Question:  \\
{\it Can we drop the hypothesis on the determinants in Theorem \ref{restrictedKP}?}
\end{center}
More formally, the Algebraic KP Question is the following ``Classification Question".  
Let $E$ and $F$ be finite graphs, and $K$ any field.  Suppose $L_K(E)$ and $L_K(F)$ are purely infinite simple.     If $K_0(L_K(E))  \cong K_0(L_K(F))$  
via an isomorphism for which $[L_K(E)] \mapsto   [L_K(F)],$
 is it necessarily the case that $L_K(E) \cong L_K(F)$?   
 
 The  name given to the Question derives from the previously mentioned Kirchberg Phillips Theorem for C$^*$-algebras (see the discussion prior to Theorem \ref{KPforgraphalgebras}), which yields as a special case that if $E$ and $F$ are finite graphs, and if $C^*(E)$ and $C^*(F)$ are purely infinite simple graph C$^*$-algebras with 
$K_0(C^*(E))  \cong K_0(C^*(F))$  via an isomorphism for which  $[C^*(E)] \mapsto   [C^*(F)],$
then $C^*(E) \cong C^*(F)$ (homeomorphically).  In particular, the determinants of the appropriate matrices play no role.  
 
 Intuitively, the Question asks whether  or not the integer ${\rm det}(I - A_E)$ can be ``seen" or ``recovered"  inside $L_K(E)$ as an isomorphism invariant.   There is indeed a way to interpret ${\rm det}(I - A_E)$ in terms of the cycle structure of $E$, see e.g. \cite{ChrisArXiV};  but this interpretation has not (yet?) been useful in this context.  
 
 With the Restricted Algebraic Kirchberg Phillips Theorem having been established, there are three possible answers to the Algebraic Kirchberg Phillips Question:

\smallskip

  {\it No.}   \ That is,  if the two graphs $E$ and $F$ have ${\rm det}(I - A_E) \neq {\rm det}(I - A_F)$, then $L_K(E) \not\cong L_K(F)$ for any field $K$.  
 
 \smallskip
  {\it Yes.} \ That is,   the existence of an isomorphism of the indicated type between the $K_0$ groups  is sufficient to yield an isomorphism of the associated Leavitt path algebras, for any field $K$.
 
 \smallskip
  {\it Sometimes.} That is,  for some pairs of graphs $E$ and $F$, and/or for some fields $K$, the answer is No, and for other pairs the answer is Yes.

\medskip

 One of the elegant aspects of the Algebraic KP Question is that its answer will be interesting, regardless of which of the three possibilities turns out to be correct.     If the answer is {\it No}, then isomorphism classes of purely infinite simple Leavitt path algebras will match exactly the flow equivalences classes of the germane set of graphs, which would suggest that there is some deeper, as-of-yet-not-understood connection between the two subjects.  If the answer is {\it Yes}, this would yield further compelling evidence for the existence of a Rosetta Stone, since then the Leavitt path algebra and graph C$^*$-algebra results would be exactly analogous.   If the answer is {\it Sometimes}, then (in addition to providing quite a surprise to those of us working in the field) this would likely require the development and utilization of a completely new set of tools in the subject. (Indeed, the {\it Sometimes} answer might be the most interesting of the three.)


Using a standard tool (the Smith Normal Form of an integer-valued matrix), it is not hard to show that the cardinality of the group $K_0(L_K(E))$ is $|{\rm det}(I - A_E)|$ in case $K_0(L_K(E))$ is finite, and the cardinality is infinite precisely when ${\rm det}(I - A_E) = 0.$  
So the Algebraic KP Question admits a somewhat more concise version:  If the signs of ${\rm det}(I - A_E)$ and ${\rm det}(I - A_F)$ are different, is it the case that $L_K(E) \not\cong L_K(F)$?   

The analogous question about Morita equivalence asks  whether or not we can drop the determinant hypothesis from Theorem 
\ref{MoritaKP}.   But the two questions will have the same answer: if isomorphic $K_0$ groups yields Morita equivalence of the Leavitt path  algebras, then the Morita equivalence together with Huang's Theorem will yield isomorphism of the algebras.  

Suppose $E$ is a finite graph for which $L_K(E)$ is purely infinite simple.  There is a way to associate with $E$ a graph $E_{-}$, for which $L_K(E_{-})$ is purely infinite simple, for which $K_0(L_K(E)) \cong K_0(L_K(E_{-}))$, and for which ${\rm det}(I - A_E) = - {\rm det}(I - A_{E_{-}})$.  This is called the ``Cuntz splice" process, which appends to a vertex $v \in E^0$ two additional vertices and six additional edges, as shown here pictorially:  
$$\qquad
E_{-} \ \  =  ((( \ E \ ))) 
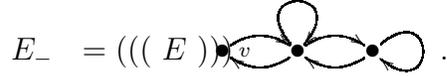
\hspace{-.1in}  \xy
(10,0)*{\bullet}="v";
(13,0)*{\scriptstyle v};
(20,0)*{\bullet}="v1";
(30,0)*{\bullet}="v2";
{\ar@/^/     "v" ;"v1" };
{\ar@/^/     "v1";"v"  };
{\ar@/^/     "v1";"v2" };
{\ar@/^/     "v2";"v1" };
{\ar@(ul,ur) "v1";"v1" };
{\ar@(ur,dr) "v2";"v2" };
\endxy.
$$
Although  the isomorphism between $K_0(L_K(E))$ and $K_0(L_K(E_{-}))$ need not in general send $[1_{L_K(E)}]$ to $[1_{L_K(E_{-})}]$,    the Cuntz splice process allows us an easy way to produce many specific examples of pairs of Leavitt path algebras to analyze in the context of the Algebraic KP Question.  The most ``basic" pair of such algebras arises from the following two graphs:

$$E_2 \  =  \ \ \ \
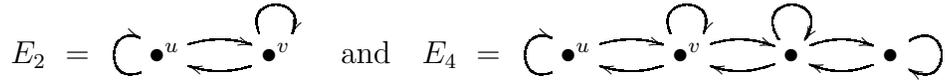
\xymatrix{
  \bullet^{u} \ar@(dl,ul)[] \ar@/^/[r]
& \bullet^{v} \ar@(ul,ur)[]  \ar@/^/[l]
}\ \ \
 \ \mbox{and}
\ \ \
E_4 \   = \ \ \ \
\xymatrix{
  \bullet^{u} \ar@(dl,ul)[] \ar@/^/[r]
& \bullet^{v} \ar@(ul,ur)[] \ar@/^/[r] \ar@/^/[l]
& \bullet^{} \ar@(ul,ur)[] \ar@/^/[r] \ar@/^/[l]
& \bullet^{} \ar@(ur,dr)[] \ar@/^/[l]
}
$$
\noindent
We note that $E_4 =  (E_2)_{-}.$  It is not hard to establish that
$$(K_0(L(E_2)),[1_{L(E_2)}])=(\{0\},0) = (K_0(L(E_4)),[1_{L(E_4)}]); $$
$$
\mbox{det}(I-A_{E_2})=-1;  \ \ \mbox{and}  \ \
\mbox{det}(I-A_{E_4})=1.$$

\noindent
 Is $L_K(E_2) \cong L_K(E_4)$?

\medskip

Here is an alternate approach to establishing the (analytic) Kirchberg Phillips Theorem (Theorem \ref{KPforgraphalgebras}) in the limited context of graph C$^*$-algebras.   Using the same symbolic-dynamics  techniques as those used to establish Theorem \ref{restrictedKP}, one can establish the C$^*$-version of the Restricted Algebraic Kirchberg Phillips Theorem (i.e., one which involves the determinants).   One then ``crosses the determinant gap" for a single pair of algebras, by showing that $C^*(E_2) \cong C^*(E_4)$; this is done using a powerful analytic tool (KK-theory).  Finally, again using analytic tools, one shows that this one particular   crossing of  the determinant gap allows for the crossing of the gap for all germane pairs of graph C$^*$-algebras.    
But neither KK-theory, nor the tools which yield the extension from one crossing to all crossings, seem to accommodate analogous algebraic techniques.

The pair $\{E_2,E_4\}$  can appropriately  be viewed  as the ``smallest" pair of graphs of interest in this context, as follows.  We say a graph has Condition (Sing) in case there are no parallel edges in the graph (i.e., that the incidence matrix $A_E$ consists only of $0$'s and $1$'s).     It can be shown that, up to graph isomorphism, there are 2 (resp., 34) graphs having two (resp., three) vertices, and having Condition (Sing), and for which the corresponding Leavitt path algebras are purely infinite simple.  (See \cite{Classification}.)   For each of these graphs $E$,  ${\rm det}(I - A_E) \leq 0$.  So finding an appropriate pair of graphs with (Sing)  and with unequal (sign of the) determinant requires at least one of the two graphs to contain at least four vertices.   




\medskip

To the author's knowledge, no Conjecture regarding what the answer to the Algebraic KP Question should be has appeared in the literature.  

\bigskip


 \subsection{The Classification Question for  graphs with finitely many vertices and infinitely many edges}

We consider now the collection $\mathcal{S}$ of those graphs $E$ having finitely many vertices, but (countably) infinitely many edges, and for which $L_K(E)$ is (necessarily unital) purely infinite simple.   The Purely Infinite Simplicity Theorem (Theorem \ref{pistheorem}) extends to this generality, so we can fairly easily determine whether or not a given graph $E$ is in  $\mathcal{S}$.   Unlike the case for finite graphs,  a description of $K_0(L_K(E))$ for $E \in \mathcal{S}$ cannot be given in terms of the cokernel of an integer-valued  matrix  transformation from $\Z^{|E^0|}$ to $\Z^{|E^0|}$.  Nonetheless, there is still a relatively easy way to determine   $K_0(L_K(E))$, so that this group remains   a useful player in this context.   

For a graph $E$ let ${\rm Sing}(E)$ denote the set of singular vertices of $E$, i.e., the set of vertices which are either sinks, or infinite emitters.    Ruiz and Tomforde in  \cite{RuizTomf} achieved the following.  

\begin{theorem}\label{RuizTomfthm}
Let $E,F \in \mathcal{S}$.    If $K_0(L_K(E)) \cong K_0(L_K(F))$ and $|{\rm Sing}(E)| = |{\rm Sing}(F)|$, then $L_K(E)$ is Morita equivalent to $L_K(F)$.   
\end{theorem}

So, while   ``the determinant of $I - A_E$" is clearly not defined  here in the usual sense (because at least one of the entries would be the symbol $\infty$), the isomorphism class of $K_0$ together with the number of singular vertices is enough information to determine Morita equivalence.    Although this is quite striking, it is not completely satisfying, in that it remains unclear whether or not $|{\rm Sing}(E)| $ is an algebraic property of $L_K(E)$.    

Continuing the search for a Classification Theorem which is cast  completely in terms of algebraic properties of the underlying algebras, the authors were able to show that for a certain type of field (those with {\it no free quotients}), there is such a result.   In a manner similar to the computation of $K_0(L_K(E))$ for $E\in \mathcal{S}$, there is a way to easily compute $K_1(L_K(E))$ as well.  

\begin{theorem}\label{nofreequotients}\cite[Theorem 7.1]{RuizTomf} 
Suppose $E,F \in \mathcal{S}$, and suppose that $K$ is a field with no free quotients.   Then $L_K(E)$ is Morita equivalent to  $L_K(F)$ if and only if $K_0(L_K(E)) \cong K_0(L_K(F))$ and $K_1(L_K(E)) \cong K_1(L_K(F))$.
\end{theorem}

The collection of fields having no free quotients includes algebraically closed fields, $\mathbb{R}$, finite fields, perfect fields of positive characteristic, and others.  However, the field $\mathbb{Q}$ is not included in this list.   Indeed, the authors in \cite[Example 10.2]{RuizTomf}  give an example of graphs $E,F \in \mathcal{S}$ for which  $K_0(L_\mathbb{Q}(E)) \cong K_0(L_\mathbb{Q}(F))$ and $K_1(L_\mathbb{Q}(E)) \cong K_1(L_\mathbb{Q}(F))$, but $L_\mathbb{Q}(E)$ is not Morita equivalent to $L_\mathbb{Q}(F)$.   
There are many open questions here.  For instance, might there be an integer $N$ for which, if $K_i(L_K(E)) \cong K_i(L_K(F))$ for all $0 \leq i \leq N$, then $L_K(E)$ and $L_K(F)$ are Morita equivalent for all fields $K$?   Of note in this context is that, unlike the situation for graph C$^*$-algebras (in which ``Bott periodicity" yields that $K_0$ and $K_1$ are the only distinct $K$-groups), there is no analogous result for the $K$-groups of Leavitt path algebras.  Further, although a long exact sequence for the $K$-groups of $L_K(E)$  has been computed in \cite[Theorem 7.6]{AraBrustCortMunster}, this sequence does not yield easily recognizable information about $K_i(L_K(E))$ for $i\geq 2$. 

Finally, a recent intriguing result presented in  \cite{GRTW-K6}  demonstrates that, if $K$ is a finite extension of $\mathbb{Q}$, then the pair consisting of ($K_0(L_K(E)), K_6(L_K(E))$) provides a complete invariant for the Morita equivalence classes of Leavitt path algebras arising from graphs in $\mathcal{S}$,  while none of the pairs $(K_0(L_K(E)), K_i(L_K(E)))$ for $1\leq i \leq 5$ provides such.

\subsection{Graded Grothendieck groups, and the corresponding Graded Classification Question}

The Algebraic Kirchberg Phillips Question, motivated by the corresponding C$^*$-algebra result, is not the only natural classification-type question to ask in the context of Leavitt path algebras.   Having in mind the importance that the $\Z$-grading on $L_K(E)$ has been shown to play in the  multiplicative structure, Hazrat in \cite{HazGradedGroth} has built the machinery which allows for the casting of an analogous question from the graded point of view.

There is a very well developed theory of graded modules over group-graded rings, see, e.g., \cite{NvO}. (The theory is built for all groups, and is particularly robust in case the group is $\Z$, the case of interest for  Leavitt path algebras.) If $A = \oplus_{t\in \Z}A_t$ is a $\Z$-graded ring and $M$ is a left $A$-module, then $M$ is {\it graded} in case $M = \oplus_{i \in \Z} M_i$, and $a_t m_i \in M_{t+i}$ whenever $a_t \in A_t$ and $m_i \in M_i$.     If $M$ is a $\Z$-graded $A$-module, and $j \in \Z$,  then the {\it suspension module} $M(j)$ is a graded $A$-module, for which $M(j) = M$ as $A$-modules, with $\Z$-grading given by  setting $M(j)_i = M_{j+i}$ for all $i,j \in \Z$. 

In a standard way, one can define the notion of a graded finitely generated projective module, and subsequently build the monoid $\mathcal{V}^{{\rm gr}}$ of isomorphism classes of such modules, with $\oplus$ as operation.    If $[M] \in \mathcal{V}^{{\rm gr}}$, then $[M(j)] \in \mathcal{V}^{{\rm gr}}$ for each $j\in \Z$, which yields a  $\Z$-action on  $\mathcal{V}^{{\rm gr}}$, and thus by extension gives $\mathcal{V}^{{\rm gr}}$ the structure of a $\Z[x,x^{-1}]$-module.   In a manner completely analogous to the non-graded case, one may define the graded Grothendieck groups $K_i^{{\rm gr}}$ for each $i\geq 0$; the suspension operation yields a $\Z[x,x^{-1}]$-module structure on these as well.

From this graded-module point of view, one can now ask about  structural information of the $\Z$-graded $K$-algebra $L_K(E)$  which might be gleaned from the  $K_i^{{\rm gr}}$ groups.  A reasonable initial question might be to see whether the graded version of the Kirchberg Phillips Theorem holds.  That is, suppose that $E$ and $F$ are finite graphs for which $L_K(E)$ and $L_K(F)$ are purely infinite simple, and suppose $K_0^{{\rm gr}}(L_K(E)) \cong K_0^{{\rm gr}}(L_K(F))$ as $\Z[x,x^{-1}]$-modules, via an isomorphism which takes $[L_K(E)]$ to $[L_K(F)]$.  Is it necessarily the case that $L_K(E) \cong L_K(F)$ as $\Z$-graded $K$-algebras?   

As it turns out, the purely infinite simple hypothesis is not the natural one to start with in the graded context.  In fact, Hazrat in \cite{HazGradedGroth} makes the following Conjecture, which at first glance might seem  somewhat audacious.   

\begin{conjecture}\label{Hazratconjecture}
Let $E$ and $F$ be any pair of row-finite graphs.  Then $L_K(E) \cong L_K(F)$ as $\Z$-graded $K$-algebras if and only if $K_0^{{\rm gr}}(L_K(E)) \cong K_0^{{\rm gr}}(L_K(F))$ as $\Z[x,x^{-1}]$-modules, via an order-preserving isomorphism which takes $[L_K(E)]$ to $[L_K(F)]$.
\end{conjecture}  

So Hazrat's conjecture, slightly rephrased, asserts that the graded $K_0$ (viewed with the $\Z[x,x^{-1}]$-module structure induced by the suspension operation), together with the natural order and position of the regular module, is a complete graded isomorphism invariant for the collection of {\it all} Leavitt path algebras over row-finite graphs.  (The order on $K_0(R)$ is induced by viewing the nonzero elements of $\mathcal{V}(R)$ as the positive elements.  The order on $K_0(R)$  plays no role in purely infinite simple rings, because every nonzero element of $\mathcal{V}(R)$ is positive in that case.) 

In \cite[Theorem 4.8]{HazGradedGroth}, Hazrat verifies Conjecture \ref{Hazratconjecture} in case the graphs $E$ and $F$ are {\it polycephalic} (essentially, mixtures of acyclic graphs, or graphs which can be described as ``multiheaded comets" or ``multiheaded roses" in which the cycles and/or  roses have no exits.)    

As mentioned in the Historical Plot Line \#1, 
in work that predates the introduction of the general definition of Leavitt path algebras  the four authors of 
\cite{AGGP04}   investigated the notion of a fractional skew monoid ring, which in particular situations is denoted $A[t_+,t_-, \alpha]$.  Recast in the language of Leavitt path algebras, the discussion in \cite[Example 2.5]{AGGP04} yields that,  when $E$ is an essential graph (i.e., has no sinks or sources), then $L_K(E) = L_K(E)_0[t_+,t_-, \alpha]$ for suitable elements $t_+,t_- \in L_K(E)$, and a corner-isomorphism  $ \alpha$   of the zero component $L_K(E)_0$.  
 
 When $E$ is a finite graph with no sinks, then $L_K(E)$ is strongly graded (\cite[Theorem 2]{Hazrat2013}), which yields (by a classical theorem of Dade) that the category of graded modules over $L_K(E)$ is equivalent to the category of (all) modules over the zero component $L_K(E)_0$.   Thus, when $E$ has no sinks, we have reason to expect that the zero component might play a role in the graded theory.     In a deep result (which relies heavily on ideas from symbolic dynamics),  Ara and Pardo \cite[Theorem 4.1]{AraPardo2014} prove the following modified version of Conjecture \ref{Hazratconjecture}.

\begin{theorem}\label{AraPardoHazratconj}
Let $E$ and $F$ be finite essential graphs.  Write $L_K(E) = $ 

\noindent
$ L_K(E)_0[t_+,t_-, \alpha]$ as described above.   Then the following are equivalent.

\smallskip
(1)  $K_0(L_K(E)_0) \cong K_0(L_K(F)_0) $ via an order-preserving $K[x,x^{-1}]$-module isomorphism 
which takes $[1_{L_K(E)_0}]$ to $[1_{L_K(F)_0}]$.

\smallskip

(2)  There exists a locally inner automorphism $g$ of $L_K(E)_0$ for which $$L_K(F) \cong   L_K(E)_0[t_+, t_-, g\circ \alpha]$$

\noindent
as $\Z$-graded $K$-algebras.

\end{theorem}

\bigskip

A complete resolution of Conjecture \ref{Hazratconjecture} currently remains elusive.   

\subsection{Connections to noncommutative algebraic geometry}


 One of the basic ideas of (standard) algebraic geometry is the correspondence between geometric spaces and commutative algebras.   Over the past few decades, significant research energy has been focused on appropriately extending this correspondence to the noncommutative case; the resulting theory is called noncommutative algebraic geometry.\footnote{Thanks to S. Paul Smith for providing much of the information contained in this subsection.}

  Suppose $A$ is a $\Z^+$-graded algebra (i.e., a $\Z$-graded algebra for which $A_n = \{0\}$ for all $n<0$).      Let ${\rm Gr}(A)$ denote the category of $\Z$-graded left $A$-modules (with graded homomorphisms), and let ${\rm Fdim}(A)$ denote the full subcategory of ${\rm Gr}(A)$ consisting of the graded $A$-modules which are the sum of their finite dimensional submodules.   Denote by ${\rm QGr}(A)$ the quotient category ${\rm Gr}(A) / {\rm Fdim}(A)$.     The category ${\rm QGr}(A)$ turns out to be one of the fundamental constructions in noncommutative algebraic geometry.  
  In particular, if $E$ is a directed graph, then the path algebra $KE$ is $\Z^+$-graded in the usual way (by setting ${\rm deg}(v) = 0$ for each vertex $v$, and ${\rm deg}(e) = 1$ for each edge $e$), and so one may construct the category ${\rm QGr}(KE)$.   
 
 Let $E^{\rm nss}$ denote the graph gotten  by repeatedly removing all sinks and sources (and their incident edges) from $E$.   
 
 \begin{theorem} \cite[Theorem 1.3]{Smith} \label{SmithTheorem}
 Let $E$ be a finite graph.  Then there is an equivalence of categories
 $${\rm QGr}(KE) \ \sim \ {\rm Gr}(L_K(E^{\rm nss})).$$
 Moreover, since $L_K(E^{\rm nss})$ is strongly graded, then these categories are also equivalent to the full category of modules over the zero-component $(L_K(E^{\rm nss}))_0$.  
  \end{theorem}
  
 So the Leavitt path algebra construction arises naturally in the context of noncommutative algebraic geometry.   (The appearance of Leavitt path algebras in this setting is clarified by the notion of a Universal Localization, see e.g. \cite{Sch}.)     
 
 In general, when the $\Z^+$-graded $K$-algebra $A$ arises as an appropriate graded deformation of the standard polynomial ring $K[x_0, ..., x_n]$, then ${\rm QGr}(A)$ shares many similarities with projective $n$-space $\mathbb{P}^n$; parallels between them have been studied extensively (see e.g. \cite{SVdB}). However, in general, an algebra of the form $KE$ does not arise in this way; and for these, as asserted in \cite{Smithnote},  ``it is much harder to see any geometry hiding in ${\rm QGr}(KE)$."     In specific situations there are some geometric perspectives available (see e.g. \cite{Smith2014}), but the general case is not well understood.

\subsection{Tensor products}

As described in Section  \ref{tensorsubsection}, the algebras $L_K(1,2) \otimes L_K(1,2)$ and $L_K(1,2)$ are not isomorphic.    However, the following related questions are still unresolved.

\smallskip

(1)    Does there exist a (necessarily injective) nonzero  ring homomorphism $\varphi: L_K(1,2)\otimes_K L_K(1,2) \rightarrow L_K(1,2)$?

(2)  Is $L_K(1,2) \otimes_K L_K(1,2)$ isomorphic to $L_K(1,2) \otimes_K L_K(1,3)$?

\subsection{The Realization Problem for von Neumann regular rings}  

Although significant progress has been made in resolving the Realization Problem  for von Neumann regular rings (see the discussion prior to Theorem \ref{AraBrustvNr}), there is as of yet not a complete answer.  An excellent survey of the main ideas relevant to this endeavor can be found in \cite{Arasurvey}.  

  Using direct limit arguments, one can show that the graph monoid $M_E$ corresponding to a countable graph $E$ can be realized as $\mathcal{V}(R)$ for a von Neumann regular algebra $R$. Indeed, $M_E$ is constructed as a direct limit of  monoids of the form $M_F$, where the graphs $F$ are finite; in particular, $M_E$ is a direct limit of finitely generated refinement monoids. Furthermore,   $R$ can be constructed as a direct limit of  (von Neumann regular) quotient algebras of the form $Q_K(F)$ for $F$ finite.   

More generally, one can divide the (countable refinement) monoids arising in the Realization Problem into two types: {\it tame} (those which can be constructed as direct limits of finitely generated refinement monoids), and the others (called {\it wild}).  Investigations  (by Ara and Goodearl, see \cite{AraGoodearlMonoids}) continue into whether or not every finitely generated refinement monoid is realizable; whether or not the realization passes to direct limits; and whether or not there are wild monoids which are not realizable.

\bigskip
\bigskip

\section{ \ Appendix 1:  \  Some properties of $L_K(E)$ and $C^*(E)$ which suggest the existence of a Rosetta Stone}\label{RosettaStoneAppendix}

It has become apparent that there is a strong, but  mysterious, relationship between the structure of the Leavitt path algebra $L_{\C}(E)$ and the corresponding graph C$^*$-algebra $C^*(E)$.  In this context it is helpful to keep in mind that while $L_{\C}(E)$ may always be viewed as a dense $*$-subalgebra of $C^*(E)$ (see Proposition \ref{densesubalgebra}), the two algebras are in general clearly different as rings:  indeed, they coincide only when $E$ is finite and acyclic.   


We focus in this Appendix on finite graphs, so that the corresponding Leavitt path algebra $L_{\mathbb{C}}(E)$ or graph C$^*$-algebra $C^*(E)$ is unital (and $C^*(E)$ is separable as well).  But many of the observations we make here hold more generally.

Any C$^*$-algebra $A$ wears two hats: not only is $A$ a ring, but $A$  comes equipped with a topology as well, so that one may view the ring-theoretic structure of $A$ from a topological/analytic  viewpoint.
The standard example is this:   one may define the (algebraic) simplicity of the C$^*$-algebra either as a ring (no nontrivial two-sided ideals), or the (topological) simplicity as a topological ring (no nontrivial closed two-sided ideals).
In general, the algebraic and topological properties of a given C$^*$-algebra $A$ need not coincide. 

The graph $E$ is called {\it cofinal} in case every vertex of $E$ connects to every cycle and every sink of $E$.  (This turns out to be equivalent to $E^0$ having the property that the only hereditary saturated subsets of $E^0$ are $\emptyset$ and $E^0$.)

As a reminder:  $E$ has Condition (L) if every cycle in $E$ has an exit; 
 $E$ has Condition (K) if there is no vertex $v$ of $E$ which has exactly one simple closed path based at $v$; and 
$E$ is downward directed if for each pair of vertices $v,w$ of $E$ there exists a vertex $y$ for which $v\geq y$ and $w \geq y$.

\bigskip

 {\bf Property 1:   \ Simplicity}
 
 \smallskip

  {\it Algebraic:}  No nontrivial two-sided ideals.  

  {\it Analytic:}    No nontrivial closed two-sided ideals.  

 By Theorem \ref{thm:AAPsimple}, $L_{\C}(E)$ is simple if and only if $E$ is cofinal and has Condition (L).

 By \cite[Proposition 5.1]{BPRS00} (for the case without sources), and \cite{RaeburnBook} (for the general case), $C^*(E)$ is (topologically) simple if and only if $E$ is cofinal and has Condition (L).

 By \cite[p. 215]{Cuntzsimple}, for any unital C$^*$-algebra $A$, $A$ is topologically simple if and only if $A$ is algebraically simple.

\smallskip

 {\it Result}:   These are equivalent for any finite graph $E$:
 \begin{enumerate}
 
 \vspace{-.05in} 
 
 \item $L_{\C}(E)$ is  simple.

 \item $C^*(E)$ is (topologically) simple.

 \item $C^*(E)$ is (algebraically) simple.

 \item  $E$ is cofinal, and satisfies Condition (L).

     \end{enumerate}

\bigskip

{\bf  Property 2:  \ The $\mathcal{V}$-monoid}

\smallskip

(Much of this discussion is taken directly from \cite[Sections 2 and 7]{AMP07}.)

{\it Algebraic}:  For a ring $R$, $\mathcal{V}(R)$ is the monoid of isomorphism classes of finitely generated projective  left $R$-modules, with operation $\oplus$.  By \cite[Chapter 3]{Blackadar}, $\mathcal{V}(R)$ can be viewed as the set of equivalence classes $V(e)$ of idempotents $e$ in the (nonunital) infinite matrix ring $M_\N(R)$, with operation
$$V(e) + V(f) = V(
\left(%
\begin{array}{cc}
  e & 0  \\
  0 & f \\
\end{array}
\right)
).$$

{\it Analytic}:   For an operator algebra $A$, $\mathcal{V}_{MvN}(A)$ is the monoid of Murray - von Neumann equivalence classes of projections in $M_\N(A)$.

By \cite[4.6.2 and 4.6.4]{Blackadar}, whenever $A$ is a C$^*$-algebra, then $\mathcal{V}(A)$ agrees with $\mathcal{V}_{MvN}(A)$.

By \cite[Theorem 7.1]{AMP07}, the natural inclusion $\psi: L_\C(E) \rightarrow
C^*(E)$ induces a monoid isomorphism $\mathcal{V}(\psi): \mathcal{V}(L_\C(E)) \rightarrow \mathcal{V}(C^*(E))$.

By \cite[Theorem 3.5]{AMP07}, the monoid $\mathcal{V}(L_K(E))$ is independent of the field $K$;  specifically, $\mathcal{V}(L_K(E))\cong M_E$, the graph monoid of $E$.

\smallskip

{\it Result}:  For any finite graph $E$ and any field $K$, the following semigroups are isomorphic.

 \vspace{-.05in} 

\begin{enumerate}

\item the graph monoid $M_E$

\item  $\mathcal{V}(L_K(E))$

\item  $\mathcal{V}(C^*(E))$

\item  $\mathcal{V}_{MvN}(C^*(E))$

\end{enumerate}

\bigskip

 {\bf Property 3: \ Purely infinite simplicity}
 
 \smallskip

{\it Algebraic}:  $R$ is purely infinite simple in case $R$ is simple and every nonzero right ideal of $R$ contains an infinite idempotent.  (Source:  \cite[Definitions 1.2]{AGPpis}.)

{\it Analytic}:   The simple C$^*$-algebra $A$ is called purely infinite (simple)  if for every positive $x\in A$, the subalgebra $\overline{xAx}$ contains an infinite projection.  (Source:  \cite[p. 186]{Cuntz2}.)

By \cite[Theorem 1.6]{AGPpis}, (algebraic) purely infinite simplicity for unital rings is equivalent to:  $R$ is not a division ring, and for all nonzero $x\in R$  there exist $\alpha, \beta \in R$ for which $\alpha x \beta = 1$.

By \cite[Proposition 6.11.5]{Blackadar}, (topological) purely infinite simplicity for unital C$^*$-algebras is equivalent to:   $A\neq \C$ and for every $x \neq 0$ in $A$ there exist $\alpha, \beta \in A$ for which $\alpha x \beta = 1$.   (Remark:  Blackadar {\it defines} purely infinite simplicity this way, and then shows this definition is equivalent to Cuntz' definition given in \cite{Cuntz2}.)  Easily,  for any graph $E$, $C^*(E)$ is a division ring if and only if $E$ is a single vertex, in which case $C^*(E) = \C$.

Thus we have, for graph C$^*$-algebras, $C^*(E)$ is (algebraically) purely infinite simple if and only if $C^*(E)$ is (topologically) purely infinite simple.

By \cite[Theorem 11]{AAPpureinf}, $L_{\C}(E)$ is purely infinite simple if and only if $L_{\C}(E)$ is simple, and $E$ has the property that every vertex connects to a cycle.

By \cite[Proposition 5.3]{BPRS00}, $C^*(E)$ is (topologically) purely infinite simple if and only if $C^*(E)$ is simple, and $E$ has the property that every vertex connects to a cycle.

\smallskip

{\it Result}:   These are equivalent for any finite graph $E$:

\begin{enumerate}
 \vspace{-.05in}

\item $L_{\C}(E)$ is purely infinite simple.

\item  $C^*(E)$ is (topologically) purely infinite simple.
\item  $C^*(E)$ is (algebraically) purely infinite simple.

\item  $E$ is cofinal, every cycle in $E$ has an exit, and every vertex in $E$ connects to a cycle.

\end{enumerate}

\bigskip

 {\bf  Property 4: \ Exchange}
 
 \smallskip

{\it Algebraic}:  $R$ is an exchange ring if for any $a\in R$ there exists an idempotent $e\in R$ for which $e\in Ra$ and $1-e \in R(1-a)$.  (Note:   The original definition of {\it exchange ring} was given by Warfield, in terms of a property on direct sum decomposition of modules; this property clarifies the genesis of the name  {\it exchange}.   The definition given here is equivalent to Warfield's; this equivalence was shown independently by Goodearl and Warfield in \cite[discussion on p. 167]{GW}, and by Nicholson in \cite[Theorem 2.1]{Nicholson}.)  

{\it Analytic}:   For every $x>0$ there exists a projection $p$ such that $p\in Ax$ and $1-p \in A(1-x)$.
(We call this condition ``topological exchange".  Note:  There does not seem to be an explicit definition of  ``topological exchange ring" in the literature.)

By \cite[Theorem 4.5]{APS}. $L_{\C}(E)$ is an exchange ring if and only if $E$ satisfies Condition (K).

By \cite[Theorem 4.1]{JP} $C^*(E)$ has real rank zero if and only if $E$ satisfies Condition (K).

By \cite[Theorem 7.2]{AGOP98}, for a unital C$^*$-algebra $A$, $A$ has real rank zero if and only if $A$ is a topological exchange ring if and only if $A$ is an exchange ring.

\smallskip

{\it Result}:   These are equivalent for a finite graph $E$:

\begin{enumerate}
\vspace{-.1in} 

\item $L_{\C}(E)$ is an exchange ring.

\item  $C^*(E)$ is a (topological) exchange ring.

\item  $C^*(E)$ is an (algebraic) exchange ring.

\item  $E$ satisfies Condition (K).

\end{enumerate}


\bigskip

 {\bf Property 5: \ Primitivity}
 
 \smallskip

{\it Algebraic}:   $R$ is (left) primitive if there exists a simple faithful left $R$-module.

{\it Analytic}:   $A$ is (topologically) primitive if there exists an irreducible faithful $*$-representation of $A$.   (That is, there is a
faithful irreducible representation $\pi : A \to B(\mathcal{H})$ for a Hilbert
space $\mathcal{H}$.)

It is shown in \cite[Theorem 4.6]{APPSMIndiana} that $L_{\C}(E)$ is left (and / or right) primitive if and only if $E$ is downward directed and satisfies Condition (L).

It is shown in \cite[Proposition 4.2]{BHRS} that $C^*(E)$ is (topologically) primitive if and only if $E$ is downward directed and satisfies Condition (L).

It is shown in \cite[Corollary to Theorem 2.9.5]{DixmierBook} that a C$^*$-algebra is algebraically primitive if and only if it is topologically primitive.

\smallskip

{\it Result}:   These are equivalent for finite graphs:

\begin{enumerate}
\vspace{-.05in}

\item $L_{\C}(E)$ is primitive.

\item $C^*(E)$ is (topologically) primitive.

\item  $C^*(E)$ is (algebraically) primitive.

\item   $E$ is downward directed and satisfies Condition (L).   

\end{enumerate}

\vspace{-.05in}
(We note that the first three properties have been shown to be equivalent for arbitrary graphs as well, with the fourth condition being replaced by:  $E$ satisfies Condition (L), is downward directed, and has the Countable Separation Property.  See Theorem \ref{primitiveLpa} and \cite{ATprimitive}.)

\smallskip

It is interesting to note that in the situations in which we have a result which suggests the existence of a Rosetta Stone, the algebraic and topological conditions on $C^*(E)$ are identical.  Perhaps there is something in this observation which will lead to a deeper understanding of why there seems to be such a strong relationship between the properties of $L_{\C}(E)$ and $C^*(E)$.

There are indeed situations  where the analogies between the Leavitt path algebras and graph C$^*$ algebras are not as tight as those presented above.  

\medskip

 {\bf   A property for which the algebraic and analytic results are not identical:  \   Primeness}
 
 \smallskip

{\it Algebraic}:   $R$ is a prime ring in case $\{0\}$ is a prime ideal of $R$; that is, in case for any two-sided ideals $I,J$ of $R$, $I \cdot J = \{0\}$ if and only if $I= \{0\}$ or $J= \{0\}$.

{\it Analytic}:  $A$ is a prime C$^*$-algebra in case $\{0\}$ is a prime ideal of $A$; that is, in case for any closed two-sided ideals $I,J$ of $R$, $I \cdot J = \{0\}$ if and only if $I= \{0\}$ or $J= \{0\}$.

In \cite[Corollary 3.10]{APPSMIndiana} it is shown that $L_{\C}(E)$ is prime if and only if $E$ is downward directed.

But by \cite[Corollaire 1]{DixmierPaper}, any separable C$^*$-algebra is (topologically) prime if and only if it is (topologically)  primitive. So (for finite $E$) $C^*(E)$ is prime if and only if $C^*(E)$ is primitive, which by the previous discussion is if and only if $E$ is downward directed {\it and} satisfies Condition (L).
(We note that since  $I \cdot J = \{0\}$ implies $\overline{I}  \cdot \overline{J} = \{0\}$, it is straightforward to show that $A$ is algebraically prime if and only if $A$ is  analytically prime.)

So for example if $E$ is the graph with one vertex and one loop, then $L_{\C}(E)$ is prime (it's an integral domain, in fact), but $C^*(E)$ is not prime.  (It's not hard to write down nonzero continuous functions on the circle which are orthogonal.)

\medskip


There are a few other situations where the properties of $L_\C(E)$ and $C^*(E)$ do not match up exactly.  For instance, the only possible values of the (algebraic) stable rank of $L_\C(E)$ are $1,2$, and $\infty$;  as well, the only possible values of the (topological) stable rank of $C^*(E)$ are $1,2$, and $\infty$.  But among individual graphs, the values may be different:  if $E = R_1$, then the stable rank of $L_\C(R_1)$ is $2$, while the stable rank of  $C^*(R_1)$  is $1$.     

In addition, we have seen in Section \ref{tensorsubsection} that $\mathcal{O}_2 \otimes \mathcal{O}_2 \cong \mathcal{O}_2$, but $L_\C(1,2) \otimes_\C L_\C(1,2) \not\cong L_\C(1,2)$.  

\bigskip

{\bf Summary of Appendix 1}:  A ``Rosetta Stone for graph algebras"  refers to an  overarching principle which would allow an understanding as to why there is such an  extremely tight (but not perfect) relationship between various properties of Leavitt path algebras and graph C$^*$-algebras, as suggested by the examples given in this section.   Does such a Rosetta Stone exist?


%
%

\section{ \ Appendix 2: \ A number-theoretic observation}

Let $L_K(1,n)$ denote the Leavitt algebra of order $n$;  so $R = L_K(1,n)$ has the property that ${}_RR \cong {}_RR^n$ as left $R$-modules.  By Theorem \ref{thm:matricesoverLeavittalgebras} we have  
$$L_K(1,n) \cong {\rm M}_{d}(L_K(1,n)) \Leftrightarrow {\rm g.c.d.}(d,n-1) = 1.$$

\noindent
Indeed, when the appropriate number-theoretic condition is  satisfied then the isomorphism may be explicitly constructed.  

\bigskip

 The key to constructing such an isomorphism lies in considering a partition of $\{1,2,...,d\}$ into two nonempty disjoint subsets $S_1 \sqcup S_2$, described as follows.

Suppose $n$ is an integer having ${\rm g.c.d.}(d,n-1)=1$.    Write $n = qd + r$ with $1 \leq r \leq d$.  As ${\rm g.c.d.}(d, n-1) = 1$ we get ${\rm g.c.d.}(d, r-1) = 1$. 

For the current discussion we  focus only on $r$.   Note we have $r\geq 1$.   Let $s = d - (r-1)$.  Since ${\rm g.c.d.}(d, r-1) = 1$ we easily see ${\rm g.c.d.}(d, s) = 1$.       Now consider the sequence $\Sigma^{d,r}$, given by 
$$\Sigma^{d,r} \ = \ 1, \ 1+s, \ 1+2s, \ ..., \  1+ (d-1)s$$
of $d$ integers, interpreted ${\rm mod} \hspace{.025in} d$.   (Here we interpret $0 \ {\rm mod } \hspace{.025in}d$ as $d \ {\rm mod} \hspace{.025in}  d$.) 
Since ${\rm g.c.d.}(d, s) = 1$, elementary number theory gives that, as a set, the elements of $\Sigma^{d,r}$ form a complete set of residues ${\rm mod}  \hspace{.025in} d$.  

\medskip

In particular, for some $i_r$ ($1\leq i_r \leq d$) we have $1 + (i_r-1)s \equiv r-1$ ${\rm mod} \hspace{.025in}  d$.  

\medskip   

Now consider these two sequences:
$$\Sigma^{d,r}_1 = 1, \ 1+s, ..., \ 1+ (i_r - 1) s \ \ \ \ \  \   \Sigma^{d,r}_2 = 1+ i_r s, \ 1 + (i_r + 1)s,\  ... , \ 1 + (d-1)s.$$
So $\Sigma^{d,r}_1$ is just the first $i_r$ elements of $\Sigma^{d,r}$, and $\Sigma^{d,r}_2$ is the remaining $d-i_r$ elements. 

\medskip

We can also consider  the partition $S^{d,r} = S_1^{d,r} \sqcup S_2^{d,r}$ of $\{1,2,...,d\}$ which corresponds to  the elements of these two sequences:
$$ S^{d,r}_1 = \{1, \ 1+s, ..., \ 1+ (i_r - 1) s\} \ \ \ \ \ \  S^{d,r}_2 = \{1+ i_r s, \ 1 + (i_r + 1)s, \ ... , \ 1 + (d-1)s\}$$

So in particular $|S^{d,r}_1| = i_r$ and $|S^{d,r}_2| = d - i_r$.  Clearly $S^{d,r}_1 \neq \emptyset$.   But $S^{d,r}_2 \neq \emptyset$ as well, since $d \in S^{d,r}_2$.  This is because the first element $1+ i_r s$ of $\Sigma^{d,r}_2$ is always $d$, as $1+ i_r s = (1 + (i_r - 1)s)  + s =  (r-1) + (d - (r-1)) = d$.

\medskip

For notational convenience, if $d,r$ are fixed then we drop the superscript $d,r$ in the sequences and subsets.

\begin{example}\label{d=3 r=2 Example}
{\rm   The case $d=3, r=2$.     \ ${\rm g.c.d.}(3,2-1)=1$.     
$ r-1=1,$  $s=d-(r-1)=3-1=2$.

The sequence $\Sigma$ starts at 1, and increases by $s = 2$ each step, and we interpret ${\rm mod} \hspace{.025in} 3$ ($1\leq i \leq 3$).     
So  we get  the sequence 
 $ \Sigma = 1,3,2. $
Since $r-1=1$,  we get 
$$\Sigma_1 = 1 \ \ \ \Sigma_2 = 3,2$$
and so
$$S_1= \{1\},  \ \ S_2 =  \{2,3\}.  $$
}
\end{example}

\bigskip

\begin{example}
{\rm   The case $d=3, r=3$.  
${\rm g.c.d.}(3,3-1)=1$.     
$ r-1=2,$  $s=d-(r-1)=3-2=1$. 

\medskip

The sequence $\Sigma$ starts at 1, and increases by $s = 1$ each step, and we interpret ${\rm mod} \hspace{.025in} 3$ ($1\leq i \leq 3$).     
So  we get  the sequence 
  $ \Sigma = 1,2,3. $
Since $r-1=2$,  we get 
$$\Sigma_1 = 1,2 \ \ \ \Sigma_2 = 3$$
and so
$$S_1= \{1,2 \}  \ \ S_2 =  \{3\}.$$
}
\end{example}

\bigskip

\begin{example}
{\rm   The case $d=13, r=9$.  
${\rm g.c.d.}(13,9-1)=1$.     
$ r-1=8,$  $s=d-(r-1)=13-8=5$. 

\medskip

The sequence $\Sigma$ starts at 1, and increases by $s = 5$ each step, and we interpret ${\rm mod} \hspace{.025in} 13$ ($1\leq i \leq 13$).     
So  we get  the sequence 
 $ \Sigma = 1,6,11,3,8,13,5,10,2,7,12,4,9. $
Since $r-1=8$,  we get 
$$\Sigma_1 = 1,6,11,3,8 \ \ \ \  \Sigma_2 = 13,5,10,2,7,12,4,9$$
and so
$$S_1= \{1,3,6,8,11 \}  \  \ \ \ S_2 =  \{2,4,5,7,9,10,12,13\}.$$
}
\end{example}

\bigskip

By solving the congruence $1 + (i_r-1)s \equiv r-1$ ${\rm mod} \hspace{.025in}  d$, we easily get 
\begin{lemma}
 $i_r \equiv (r-1)^{-1}$  ${\rm mod}  \hspace{.025in} d$. 
\end{lemma}

In particular, if we have $1 \leq r \neq r' \leq d$ for which ${\rm g.c.d.}(d, r-1) = 1 = {\rm g.c.d.}(d, r'-1)$, then the two partitions  $S_1^{d,r} \sqcup S_2^{d,r}$ and $S_1^{d,r'} \sqcup S_2^{d,r'}$ of $\{1,2,...,d\}$ are necessarily different (since $1$ is in $S_1$, and the sizes of  $S_1^{d,r}$ and $S_1^{d,r'}$ are unequal).

\medskip

Given $d$, there exist $\varphi(d)$ (Euler $\varphi$-function) remainders which are relatively prime to $d$.   So there exist $\varphi(d)$ distinct partitions of $\{1,2,...,d\}$ which arise as $S_1^{d,r} \sqcup S_2^{d,r}$ for some $r$ having ${\rm g.c.d.}(d,r-1) = 1$.

\smallskip

We note that for any $d$ we always get these two partitions arising in the form $S^{d,r}$:
$$ \{1\} \sqcup \{2,3,...,d\} = S^{d,2}, \ \ \mbox{and} \ \ \{1,2,...,d-1\} \sqcup \{d\} = S^{d,d}.$$

It is easy to see that there are $(2^d - 2) / 2 = 2^{d-1} - 1$  ways to partition the set $\{1,2,...,d\}$ into two nonempty subsets $S_1 \sqcup S_2$ for which $1 \in S_1$.     Since $\varphi(d) < 2^{d-1} - 1$ for $d \geq 3$, we see that not all such two-nonempty-set partitions of $\{1,2,...,d\}$ can arise as $S_1^{d,r} \sqcup S_2^{d,r}$ for some $r$ having ${\rm g.c.d.}(d,r-1) = 1$.   For example, when $d=3$, the partition $\{1,3\} \sqcup \{2\}$ of $\{1,2,3\}$ does not arise in this way.  

\medskip

 {\bf We are interested in two related questions regarding the sequences described in this Appendix.}  
 
 \smallskip
 
(1)    For fixed $d,r$ having ${\rm g.c.d.}(d,r-1)=1$, do the sequences $\Sigma_1^{d,r}$ and  $\Sigma_2^{d,r}$ arise in  contexts other than that of isomorphisms between matrix rings over Leavitt algebras?   
 
(2)    Do the $\varphi(d)$ partitions of $\{1,2,...,d\}$ of the form $\{1,2,...,d\} = S_1^{d,r} \sqcup S_2^{d,r}$ (for some $r$ having ${\rm g.c.d.}(d,r-1) = 1$)   play a special role in any sorts of number-theoretic investigations?

\medskip

\begin{remark}
{\rm Referring back to how these sequences and partitions arose in the context of Theorem \ref{thm:matricesoverLeavittalgebras}, in that setting we start with $d,n$ having ${\rm g.c.d.}(d,n-1)=1$, write $n = qd+r$ with $1\leq r \leq d$, and then consider the partition $S_1^{d,r} \sqcup S_2^{d,r}$ of $\{1,2,...,d\}$.   We then use this partition of $\{1,2,...,d\}$ to build a partition of $\{1,2,...,n\}$  by simply extending the partition  $S_1^{d,r} \sqcup S_2^{d,r}$, ${\rm mod} \hspace{.025in} d$.  So, for instance, if $n=5$, $d=3$ then we get $5 = 1\cdot 3 + 2$.  We then consider the partition $\{1,2,3\} = S_1^{3,2} \sqcup S_2^{3,2} = \{1\} \sqcup \{2,3\}$, as described in Example \ref{d=3 r=2 Example}.   This then yields the partition $\{1,4\} \sqcup \{2,3,5\}$ of $\{1,2,3,4,5\}$ by simply extending ${\rm mod} \hspace{.025in} 3$. 

Specifically, in the proof of Theorem \ref{thm:matricesoverLeavittalgebras},  the ordering properties of the sequences $\Sigma_1$ and $\Sigma_2$ are not utilized, rather, only the partition $S_1^{d,r} \sqcup S_2^{d,r}$ of $\{1,2,...,d\}$ as sets is used.}

\end{remark}

%
%
%
%

\section{\ Appendix 3: \ The graph moves}


We give in this Appendix the formal definitions of each of the six  ``graph moves" which arise in the symbolic dynamics analysis associated to the Restricted Algebraic Kirchberg Phillips Theorem.   We conclude by presenting the ``source elimination" process as well.  

\begin{definition}\label{def_outsplit}
{\rm Let $E = ( E^0 , E^1 , r , s )$ be a directed graph. For each
$v \in E^0$ with $s^{-1}(v) \ne \emptyset$, partition the set
$s^{-1} (v)$ into disjoint nonempty subsets $\mathcal{E}^1_v ,
\ldots , \mathcal{E}^{m(v)}_v$ where $m(v) \ge 1$. (If $v$ is a sink,
then we put $m(v)=0$.) Let $\mathcal{P}$ denote the resulting
partition of $E^1$. We form the {\em out-split graph} $E_s({\mathcal
P})$ from $E$ using the partition $\mathcal{P}$ as follows:
\begin{align*}
E_s ( \mathcal {P} )^0 &= \{ v^i \mid v \in E^0 , 1 \le i \le m(v)
\}
\cup \{ v \mid m(v)  = 0 \} ,  \\
E_s ( \mathcal{P} )^1  &= \{ e^j \mid e \in E^1, 1 \le j \le m ( r
(e) ) \} \cup \{ e \mid m (r(e)) = 0 \} ,
\end{align*}

\noindent and define $r_{E_s ( \mathcal{P} )} , s_{E_s (
\mathcal{P} )} : E_s ( \mathcal{P} )^1 \rightarrow E_s (
\mathcal{P} )^0$ for each $e \in \mathcal{E}^i_{s(e)}$ by
$$
s_{E_s ( \mathcal{P} )} ( e^j ) = s(e)^i  \text{ and } s_{E_s (
\mathcal{P} )} ( e
) = s(e)^i,  \ \ \ 
r_{E_s ( \mathcal{P} )} ( e^j ) = r(e)^j \text{ and } r_{E_s (
\mathcal{P} )} ( e ) = r(e).
$$
\noindent Conversely, if $E$ and $G$ are graphs, and there exists a
partition ${\mathcal P}$ of $E^1$ for which $E_s({\mathcal P})=G$,
then $E$ is called an {\em out-amalgamation} of $G$. }
\end{definition}

\medskip

\begin{definition}
{\rm Let $E = (E^0, E^1, r, s)$ be a directed graph, and let $v
\in E^0$. Let $v^*$ and $f$ be symbols not in $E^0 \cup E^1$.   We
form the {\em expansion graph} $E_v$ from $E$ at $v$ as follows:
\begin{align*}
E_v^0      &= E^0 \cup \{ v^\ast \} \\
E_v^1      &= E^1 \cup \{ f \} \\
s_{E_v}(e) &= \left\{ \begin{array}{ll}
              v      & \textrm{ if $e = f$} \\
              v^\ast & \textrm{ if $s_E(e) = v$} \\
              s_E(e) & \textrm{ otherwise}
              \end{array}\right. \\
r_{E_v}(e) &= \left\{ \begin{array}{ll}
              v^\ast & \textrm{ if $e = f$} \\
              r_E(e) & \textrm{ otherwise}
              \end{array}\right. \\
\end{align*}
\noindent Conversely, if $E$ and $G$ are graphs, and there exists a
vertex $v$ of $E$ for which $E_v=G$, then $E$ is called a {\em
contraction} of $G$.}
\end{definition}

\medskip

\begin{definition}\label{def_insplit}
{\rm
Let $E = ( E^0 , E^1 , r , s )$ be a directed graph. For each $v
\in E^0$ with $r^{-1}(v) \ne \emptyset$, partition the set $r^{-1}
(v)$ into disjoint nonempty subsets $\mathcal{E}^v_1 , \ldots ,
\mathcal{E}^v_{m(v)}$ where $m(v) \ge 1$. (If $v$ is a source then
we put $m(v)=0$.) Let $\mathcal{P}$ denote the resulting partition
of $E^1$. We form the {\em in-split graph} $E_r({\mathcal P})$
from $E$ using the partition $\mathcal{P}$ as follows:
\begin{align*}
E_r ( \mathcal {P} )^0 &= \{ v_i \mid v \in E^0 , 1 \le i \le m(v)
\}
\cup \{ v \mid m(v)  = 0 \} ,  \\
E_r ( \mathcal{P} )^1  &= \{ e_j \mid e \in E^1, 1 \le j \le m ( s
(e) ) \} \cup \{ e \mid m (s(e)) = 0 \} ,
\end{align*}

\noindent and define $r_{E_r ( \mathcal{P} )} , s_{E_r (
\mathcal{P} )} : E_r ( \mathcal{P} )^1 \rightarrow E_r (
\mathcal{P} )^0$ by
\begin{align*}
s_{E_r ( \mathcal{P} )} ( e_j ) &= s(e)_j  \text{ and } s_{E_r (
\mathcal{P} )} ( e
) = s(e) \\
r_{E_r ( \mathcal{P} )} ( e_j ) &= r(e)_i \text{ and } r_{E_r (
\mathcal{P} )} ( e ) = r(e)_i \text{ where } e \in
\mathcal{E}^{r(e)}_i .
\end{align*}
\noindent Conversely, if $E$ and $G$ are graphs, and there exists a
partition ${\mathcal P}$ of $E^1$ for which $E_r({\mathcal P})=G$,
then $E$ is called an {\em in-amalgamation} of $G$.}
\end{definition}

\medskip

\begin{definition}
{\rm Let $E = (E^0, E^1, r, s)$ be a directed graph with at least
two vertices, and let $v \in E^0$ be a source.  We form the {\em
source elimination graph} $E_{\backslash v}$ of $E$ as follows:
\begin{align*}
E_{\backslash v}^0   &= E^0 \backslash \{ v \} \\
E_{\backslash v}^1   &= E^1 \backslash s^{-1}(v) \\
s_{E_{\backslash v}} &= s |_{E_{\backslash v}^1} \\
r_{E_{\backslash v}} &= r |_{E_{\backslash v}^1}
\end{align*}
}
\end{definition}






\begin{thebibliography}{}



\bibitem{Abra97}  Abrams, G.:  Non-induced isomorphisms of matrix rings.  Israel J. Math. {\bf 99}, 343-347 (1997)

\bibitem{Abra01}  Abrams, G.:  Invariant basis number and types for strongly graded rings.  J. Algebra {\bf 237}, 32-37 (2001)

\bibitem{AbAn02}   Abrams, G.,   \'{A}nh, P.N.: Some ultramatricial algebras which arise as intersections of Leavitt algebras.   J. Alg. App. {\bf 1}(4), 357--363 (2002)


\bibitem{Classification}  Abrams, G., \'{A}nh, P. N.,  Louly, A., Pardo, E.:   The classification question for Leavitt path algebras. J. Algebra {\bf 320}(5), 1983--2026  (2008)


\bibitem{AbAnP}  Abrams, G.,  \'{A}nh, P. N.,  Pardo, E.:   Isomorphisms between Leavitt algebras and their matrix rings. J. Reine Angew. Math. {\bf 624}, 103--132  (2008)


  \bibitem{TheBook}   Abrams, G., Ara, P., Siles Molina, M.:  Leavitt path algebras.   Lecture Notes in Mathematics, Springer Verlag (to appear)
  

 \bibitem{AAP05} 
      Abrams, G.,  Aranda Pino, G.:  The Leavitt path algebra of a graph.   J. Algebra  {\bf 293}, 319--334 (2005)


\bibitem{AAPpureinf} Abrams, G.,  Aranda Pino, G.: Purely infinite simple Leavitt path algebras.  J. Pure Appl. Algebra \textbf{207},  553--563 (2006)


 \bibitem{AAPHouston}   Abrams, G.,  Aranda Pino, G.:   The Leavitt path algebras of arbitrary graphs. Houston J. Math.  {\bf 34}(2), 423--442 (2008)
  
  
  \bibitem{AbBellRanga}  
  Abrams, G., Bell, J.P., Rangaswamy, K.M.: On prime non-primitive von Neumann regular algebras.   Trans. Amer. Math. Soc. {\bf 366}(5), 2375--2392 (2014) 
  

\bibitem{AbramsBellColakRanga}  Abrams, G., Bell, J. P.,  Colak, P.,  Rangaswamy, K.M.:   Two-sided chain conditions in Leavitt path algebras over arbitrary graphs. J. Algebra Appl. {\bf 11}(3),  1250044 - 23 pp.  (2012)


  \bibitem{AbramsKanuni}
  Abrams, G., Kanuni, M.:  Cohn path algebras have Invariant Basis Number.   ArXiV:  1303.2122v2 (2013)



\bibitem{Flow}  Abrams, G., Louly, A.,  Pardo, E.,  Smith, C.:   Flow invariants in the classification of Leavitt path algebras. J. Algebra {\bf 333}, 202--231  (2011)


  \bibitem{AbManTo}   Abrams, G., Mantese, F., Tonolo, A.:  Extensions of simple modules over Leavitt path algebras.  ArXiV:  1408.3580 (2014)
  
  
  \bibitem{AbramsMesyan}  
  Abrams, G., Mesyan, Z.:  Simple Lie algebras arising from Leavitt path algebras. J. Pure Appl. Algebra {\bf 216}(10), 2302--2313  (2012) 


\bibitem{AbRanga}
 Abrams, G., Rangaswamy, K. M.:  
 Regularity conditions for arbitrary Leavitt path algebras.  
Algebr. Represent. Theory   {\bf 13}(3), 319--334 (2010)
    
  
  \bibitem{AbramsRanga}  Abrams, G., Rangaswamy, K.M.:  Row-finite equivalents exist only for row-countable graphs.   in:  New trends in noncommutative algebra.  Contemporary  Math Series {\bf 562}, 1--10.   American Mathematical Society, Providence, RI  (2012) 
  
    
\bibitem{AT}  Abrams, G., Tomforde, M.: Isomorphism and Morita equivalence of graph algebras. Trans. Amer. Math. Soc. {\bf 363}(7), 3733--3767 (2011)


\bibitem{ATprimitive}   Abrams, G., Tomforde, M.:  A class of C$^*$-algebras that are prime but not primitive.    M\"unster J. Math. (to appear)
    
    \bibitem{AA2014}   Alahmadi, A., Alsulami, H.:   Wreath products by a Leavitt path algebra and affinizations.   Int. J. Algebra Comput.  DOI: 10.1142/S0218196714500295  (2014). 
    
  \bibitem{AAJZ}   Alahmadi, A.,  Alsulami, H.,  Jain, S. K.,  Zelmanov, E.:   Leavitt path algebras of finite Gelfand-Kirillov dimension.   J. Algebra Appl. {\bf 11}(6), 1250225 - 6 pp. (2012)
  
  \bibitem{AA-Skew}   Alahmadi, A., Alsulami, H.:    Simplicity of the Lie algebra of skew symmetric elements of a Leavitt path algebra.    arXiv:1304.2385 (2014)
   
\bibitem{AA-Lie}  Alahmadi, A., Alsulami, H.:      On the simplicity of Lie algebra of Leavitt path algebra.    arXiv:1304.1922  (2013)  
 
  
 \bibitem{AAJZ-PNAS}   Alahmedi, A.,  Alsulami, H., Jain, S.K.,  Zelmanov, E.:   Structure of Leavitt path algebras of polynomial growth. Proc. Natl. Acad. Sci. USA {\bf 110}(38), 15222Ð15224 (2013)

\bibitem{Arasurvey}  Ara, P.:  The realization problem for von Neumann regular rings.  in:  Ring Theory 2007, World Scientific Publishers, Hackensack, NJ, 21-37 (2009)

\bibitem{AraBrustenga} Ara, P.,  Brustenga, M.: 
   The regular algebra of a quiver. 
  J. Algebra {\bf 309}, 207--235 (2007) 
     

\bibitem{AraBrustCortMunster}
   Ara, P., Brustenga, M., Corti{\~n}as, G.:  $K$-theory of Leavitt path algebras.  
M\"unster J. Math. 
 {\bf 2},  5--33 (2009)
  
  

\bibitem{AraCortinas}  
Ara, P., Corti\~{n}as, G.:   Tensor products of Leavitt path algebras. Proc. Amer. Math. Soc. {\bf 141}(8), 2629 -- 2639 (2013)


\bibitem{AGGP04}
Ara, P., Gonz{\'a}lez-Barroso, M. A., Goodearl, K. R., Pardo, E.: 
  Fractional skew monoid rings.
  J. Algebra {\bf 278}, 104--126 (2004)
 

\bibitem{AraGoodearlCstar}   Ara, P., Goodearl, K.R.:  C$^*$-algebras of separated graphs.  J. Funct. Anal. {\bf 261}(9), 2540-2568 (2011)

\bibitem{AraGoodearl}  Ara, P., Goodearl, K.R.:  Leavitt path algebras of separated graphs.  J. Reine Angew. Math. {\bf 669}, 165--224 (2012)


\bibitem{AraGoodearlMonoids}  Ara, P., Goodearl, K.R.:  Tame and wild refinement monoids.   arXiv: 1405.7582.   



\bibitem{AGOP98} 
Ara, P.,  Goodearl, K. R.,  O'Meara, K. C., Pardo, E.:
    Separative cancellation for projective modules over exchange
             rings. 
 Israel J. Math.  {\bf 105}, 105--137 (1998)

 \bibitem{AGOR}  Ara, P., Goodearl, K.R., O'Meara, K.C., Raphael, R.:  $K_1$ of separative exchange rings and C$^*$-algebras of real rank zero.  Pacific J. Math \textbf{195}(2),  261--275 (2000)



\bibitem{AGPpis}
Ara, P.,  Goodearl, K. R.,  Pardo, E.: 
 {$K_0$} of purely infinite simple regular rings.
  $K$-Theory  {\bf 26}, 69--100 (2002)
    
 

\bibitem{AMP07}
Ara, P.,  Moreno, M. A., Pardo, E.: 
  Nonstable {$K$}-theory for graph algebras.  Algebr. Represent. Theory  {\bf 10}(2), 157--178 (2007)
   

\bibitem{AraPardo2014}  Ara, P., Pardo, E.:  Towards a $K$-theoretic characterization of graded isomorphisms between Leavitt path algebras.  
 $K$-theory.  (to appear)  
 
 \bibitem{AraRanga1}  Ara, P., Rangaswamy, K.M.:  Finitely presented modules over Leavitt path algebras.  J. Algebra {\bf 417}, 333-352 (2014)
 
 \bibitem{AraRanga2}   Ara, P., Rangaswamy, K.M.:  Leavitt path algebras with at most countably many irreducible representations.   Rev. Mat. Iberoam. (to appear)  arXiv:  1309.7917v2 

  \bibitem{AaHCR}  
  Aranda Pino, G., Clark, J., an Huef, A.,  Raeburn, I.:  Kumjian-Pask algebras of higher-rank graphs. Trans. Amer. Math. Soc. {\bf 365}(7), 3613--3641 (2013)
  


\bibitem{ArandaPinoCrow}
 Aranda Pino, G., Crow, K.: The center of a Leavitt path algebra.  
Rev. Mat. Iberoam.  {\bf 27}(2), 621--644 (2011)
     
     
     \bibitem{Malaga4}  Aranda Pino, G., Mart\'{i}n Barquero, D., Mart\'{i}n Gonz\'{a}lez, C., Siles Molina, M.:  The socle of a Leavitt path algebra.  J. Pure App. Alg. {\bf 212}, 500--509 (2008)
     

\bibitem{APS} 
Aranda Pino, G., Pardo, E., Siles Molina, M.:   Exchange Leavitt path algebras and stable rank.   J. Algebra \textbf{305}, 912--936  (2006)


    \bibitem{APPSMIndiana}  
    Aranda Pino, G.,  Pardo, E., Siles Molina, M.:   Prime spectrum and primitive Leavitt path algebras.   Indiana J. Math. \textbf{58}, 869--890  (2009)

\bibitem{Malaga2006}
  Aranda Pino, G., Perera, F.,  Siles Molina, M. (editors): 
   Graph algebras: bridging the gap between Analysis and Algebra.   
    M\'alaga University Press,  M\'alaga 
     (2007)
     

      \bibitem{BHRS}  
       Bates, T.,   Hong, J.,  Raeburn, I.,  Szyma\'{n}ski, W.:  The ideal structure of the C$^*$-algebras of infinite graphs.  Illinois J. Math. \textbf{46}(4), 1159--1176  (2002)

           \bibitem{BPRS00}   Bates, T.,   Pask, D.,   Raeburn, I.,  Szyma\'{n}ski, W.:  The C$^*$-algebras of row-finite graphs.  New York J. Math. \textbf{6}, 307-324  (2000)


\bibitem{BergUnpub}  Bergman, G.M.:   On Jacobson radicals of graded rings.  Unpublished

\bibitem{Berg74}
Bergman, G. M.:  
  Coproducts and some universal ring constructions.  
 Trans. Amer. Math. Soc.  {\bf 200}, 33--88 (1974)
        

 \bibitem{Blackadar}  Blackadar, B.:  K-theory for Operator Algebras, 2nd Edition.   MSRI Publications no. 5, Cambridge University Press, Cambridge (1998)

\bibitem{Brin}  Brin, M.G.:  Higher dimensional Thompson groups.  Geom. Dedicata  {\bf 108}, 163-192 (2004)

\bibitem{Brookfield} 
Brookfield, G.:  
  Cancellation in primely generated refinement monoids.  
 Algebra Universalis  {\bf 46}(3),  343--371 (2001)
   

\bibitem{CartanEilenberg}  Cartan, H., Eilenberg, S.:  Homological algebra.  Princeton University Press, Princeton, NJ (1956)


\bibitem{Chen}  Chen, X.W.:   Irreducible representations of Leavitt path algebras.   Forum Math. (to appear; published online 12/14/2012)


   \bibitem{Cohn66}
 Cohn, P. M.: 
   Some remarks on the invariant basis property.  
Topology  {\bf 5}, 215--228 (1966) 
  

\bibitem{Cuntz} 
Cuntz, J.:  
  Simple $C^*$-algebras generated by isometries.  
  Comm. Math. Phys.   {\bf 57}(2),  173--185 (1977)
       

\bibitem{Cuntzsimple}   Cuntz, J.:    The structure of multiplication and addition in simple C$^*$-algebras.    Math. Scand. \textbf{40}, 215--233  (1977)

    \bibitem{Cuntz2} Cuntz, J:   $K$-theory for certain C$^*$-algebras.  Ann. Math.  \textbf{113}, 181--197  (1981)


\bibitem{CuntzKrieger}
 Cuntz, J.,   Krieger, W.: 
  A class of {$C\sp{\ast} $}-algebras and topological Markov
             chains.  
Invent. Math. 
  {\bf 56}(3), 251--268 (1980)
  
  \bibitem{DMP}  Dicks, W., Mart\'{i}nez-P\'{e}rez, C.:  Isomorphisms of Brin-Higman-Thompson groups.  Israel J. Math. {\bf 199}, 189-218 (2014)

\bibitem{DixmierPaper}  Dixmier, J.:  Sur les C$^*$-algebres.  Bull. Soc. Math.  France \textbf{88},  95--112 (1960)


\bibitem{DixmierBook} Dixmier, J.:  Les C*-algebres et leurs representations.   Gauthier-Villars, Paris (1969)

\bibitem{FLR}  Fowler, N.J., Laca, M., Raeburn, I.:  The C$^*$-algebras of infinite graphs.  Proc. Amer. Math. Soc. {\bf 8}, 2319--2327 (2000)

     
  \bibitem{Franks} 
  Franks, J.:  
    Flow equivalence of subshifts of finite type.  
  Ergodic Theory Dynam. Systems  {\bf 4}(1), 53--66 (1984) 
      


\bibitem{GRTW-K6}  Gabe, J., Ruiz, E., Tomforde, M., Whalen, T.:  $K$-theory for Leavitt path algebras:  computation and classification.  ArXiV: 1407.5094v1 (2014) 

  
\bibitem{GoodearlLimits} 
 Goodearl, K. R.:  
   Leavitt path algebras and direct limits.  
in:   Rings, Modules and Representations, 165--188.  
  Contemporary Math. Series {\bf 480}, Amer. Math. Soc., Providence, RI   (2009)  
 

\bibitem{GW}   Goodearl, K. R.,  Warfield, R. B., Jr.:   Algebras over zero-dimensional rings.  Math. Ann. {\bf 223}(2), 157-168 (1976)

  
  \bibitem{HazGradedGroth}  
  Hazrat, R.:  The graded Grothendieck group and the classification of Leavitt path algebras.  Math. Ann. {\bf 355}(1), 273 -- 325 (2013)
  
  
    
  \bibitem{Hazrat2013}  
 Hazrat, R.:  
    A note on the isomorphism conjectures for Leavitt path algebras.  
J. Algebra {\bf 375}, 33--40 (2013)



     
\bibitem{Jacobson1950}  Jacobson, N.:   Some remarks on one-sided inverses.  Proc Amer. Math. Soc. {\bf 1},  352-355 (1950) 

\bibitem{JP}  Jeong, J.A., Park, G. H.: Graph C$^*$-algebras with real rank zero.  J. Funct. An. \textbf{188}, 216-226  (2002)


\bibitem{KrauseLenagan}   Krause, G. R.,  Lenagan, T.H.:   Growth of algebras and Gelfand-Kirillov dimension. Revised edition. Graduate Studies in Mathematics {\bf 22}. Amer.  Math. Soc., Providence, RI (2000) 


  \bibitem{KumjianPask}  Kumjian, A., Pask, D.:   
Higher rank graph C$^*$-algebras. New York J. Math.  {\bf 6}, 1--20 (2000)


  
  \bibitem{KPR98}  Kumjian, A.,  Pask, D.,   Raeburn, I.:   
Cuntz-Krieger algebras of directed graphs. Pacific J. Math.  {\bf 184}(1), 161--174  (1998)


    
  \bibitem{KPRR97}   Kumjian, A.,  Pask, D., Raeburn, I.,  Renault, J.:  
Graphs, groupoids, and Cuntz-Krieger algebras. 
J. Funct. Anal. {\bf 144}(2), 505--541 (1997)  
  
    

\bibitem{Leav62}  
Leavitt, W. G.:  
The module type of a ring.  Trans. Amer. Math. Soc.  {\bf 103}, 113--130 (1962)
    

\bibitem{Leav64} 
Leavitt, W. G.:  
 The module type of homomorphic images.  
  Duke Math. J.  {\bf 32}, 305--311 (1965) 


  \bibitem{LindMarcus}  Lind, D., Marcus, B.:   An introduction to symbolic dynamics and coding. Cambridge University Press, Cambridge  (1995)  
  
 

\bibitem{NTV}  N\v{a}st\v{a}sescu, C., Torrecillas, B.,  Van Oystaeyen, F.:   IBN for graded rings.  Comm. Alg. {\bf 28}(3), 1351-1360 (2000)


\bibitem{NvO}    N\v{a}st\v{a}sescu, C.,  
Van Oystaeyen, F.:   Methods of graded rings.  Lecture Notes in Mathematics {\bf 1836},   Springer-Verlag, Berlin (2004)   

\bibitem{Nicholson}  Nicholson, W. K.:   Lifting idempotents and exchange rings. Trans. Amer. Math. Soc. {\bf 229}, 269-178 (1977)

\bibitem{Pardo2011}  Pardo, E.:   The isomorphism problem for Higman-Thompson groups. J. Algebra {\bf 344}, 172--183 (2011)

\bibitem{PardoMadVet}  Pardo, E.:  (private communication)

\bibitem{PaschkeSalinas}  Paschke, W. L.,  Salinas, N.:   Matrix algebras over $\mathcal{O}_n$. Michigan Math. J.  {\bf 26}(1),  3--12  (1979)


\bibitem{Phillips}  Phillips, N.C.:  A classification theorem for nuclear purely infinite simple C$^*$-algebras.  Doc. Math. \textbf{5}, 49-114  (2000)


\bibitem{RaeburnBook}  Raeburn, I.: Graph algebras. CBMS Regional Conference Series in Mathematics, vol. 103.  American Mathematical Society, Providence, RI  (2005) 

\bibitem{Rangasimples}  Rangaswamy, K.M.:  Notes on simple modules over Leavitt path algebras.   arXiv:  1401.6589  

\bibitem{Rordam}   R{\o}rdam, M.:  A short proof of Elliott's theorem.  C.R. Math. Rep. Acad. Sci. Canada {\bf 16}, 31--36 (1994)  

\bibitem{RowenBook}   Rowen, L.H.:  Ring Theory.  Volume 1.  Pure and Applied Mathematics {\bf 127}.  Academic Press, Boston, MA  (1988) 

  
  \bibitem{RuizTomf}  Ruiz, E.,  Tomforde, M.:  Classification of unital simple Leavitt path algebras of infinite graphs. J. Algebra {\bf 384}, 45--83 (2013)
  
\bibitem{Sch}  Schofield, A.H.:  Representations of Rings Over Skew Fields.   London Math. Soc., Lecture Notes Series, Vol. 92.  Cambridge University Press, Cambridge, UK (1985)  

\bibitem{ChrisArXiV}  Smith, C.:   A description of the Parry-Sullivan Number of a graph using circuits.   ArXiV: 0903.2000v1 (2009)

\bibitem{Smith2014}  Smith, S.P.:  The space of Penrose tilings and the non-commutative curve with homogeneous coordinate ring $k\langle x,y \rangle / (y^2)$.  J. Noncomm. Geom. {\bf 8}, 541--586 (2014)

\bibitem{Smith}  Smith, S.P.:  Category equivalences involving graded modules over path algebras of quivers.  Advances in Mathematics {\bf 230}, 1780--1810 (2012) 



\bibitem{Smithnote}   Smith, S.P.:  (private communication)

\bibitem{SVdB}   Stafford, J.T., Van den Bergh, M.:   Noncommutative curves and noncommutative surfaces.  Bull. Amer. Math. Soc. {\bf 38}, 171--216 (2001) 

\bibitem{TomfordeUniqueness}    Tomforde, M.:  Uniqueness theorems and ideal structure for Leavitt path algebras. J. Algebra {\bf 318}(1), 270--299  (2007)

\bibitem{Tom}  Tomforde, M.:  Leavitt path algebras with coefficients in a commutative ring.  
  J. Pure Appl. Algebra 
{\bf 215}(4), 471-484 (2011)
      

\bibitem{MarksHistoryBook} Tomforde, M.:   Graph C$^*$-algebras: Theory, Technique, and Examples.  (in preparation)

  
  \bibitem{Watatani82}   Watatani, Y.:   Graph theory for C$^*$-algebras.  in:  Operator algebras and applications, Part I (Kingston, Ont., 1980), 195--197.   Proc. Sympos. Pure Math. {\bf 38}, Amer. Math. Soc., Providence, RI (1982)  
  
  
  
    
    
    
    
    
    
    
    
    

    


    
















    
    
\end{thebibliography}
\end{document}